%
%
%
%

\input amstex
\documentstyle{amsppt}

\input pictex
\NoBlackBoxes

\topmatter
\title Classification of knotted tori in the 2-metastable dimension \endtitle
\author M. Cencelj, D. Repov\v s and M. Skopenkov \endauthor

\address Institute of Mathematics, Physics and Mechanics,
and Faculty of Education, University of
Ljubljana, P. O. Box 2964, 1001 Ljubljana, Slovenia.
e-mail:
matija.cencelj\@guest.arnes.si \endaddress
\address Faculty of Mathematics and Physics, and Faculty of Education, University of
Ljubljana, P. O. Box 2964, 1001 Ljubljana, Slovenia.
e-mail:
dusan.repovs\@guest.arnes.si \endaddress
\address Institute for information transmission problems RAS,
Bolshoy Karetny per. 19, Moscow, 127994, Russia
and
King Abdullah University of Science and Technology, P.O. Box 2187, 4700 Thuwal, 23955-6900, Saudi Arabia.
e-mail:
skopenkov\@rambler.ru \endaddress

\thanks Cencelj and Repov\v s were supported in part by Slovenian
Research Agency grants P1-0292-0101 and J1-4144-0101.
Skopenkov was supported in part by RFBR grant 12-01-00748-а, President of the Russian Federation grant MK-3965.2012.1,``Dynasty'' foundation and Simons--IUM fellowship.
\endthanks



\abstract This paper is devoted to the classical Knotting Problem: for a given manifold $N$ and number $m$
describe the set of isotopy classes of embeddings $N\to S^m$.
We study the specific case of {\it knotted tori}, i. e. the embeddings $S^p\times S^q\to
S^m$. The classification of knotted tori up to isotopy in the {\it metastable} dimension range $m\ge p+\frac{3}{2}q+2$, $p\le q$, was given
by A. Haefliger, E. Zeeman and A.B. Skopenkov. We consider
the dimensions below the metastable range, and give an explicit criterion for the finiteness of this set of isotopy classes in the
{\it 2-metastable} dimension:

\smallskip
{\bf Theorem.} Assume that $p+\frac{4}{3}q+2 < m <
p+\frac{3}{2}q+2$ and $m > 2p+q+2$. Then the set of isotopy classes of smooth
embeddings $S^p\times S^q\to S^m$ is infinite if and
only if either $q+1$ or $p+q+1$ is divisible by $4$.
\smallskip

Our approach to the classification is based on an analogue of the
U.~Koschorke exact sequence from the theory of link maps. Our
sequence involves a new {\it $\beta$-invariant} of knotted tori.
The exactness is proved using embedded surgery and the
N.~Habegger--U.~Kaiser techniques of studying the complement.
\endabstract
\endtopmatter

\document

\head 1. Introduction \endhead

This paper is devoted to the classical Knotting Problem: {\it for a given manifold $N$ and a number $m$
describe the set of isotopy classes of embeddings $N\to S^m$}.
For a recent survey see \cite{ReSk99, Sko07L}. This subject was actively studied in
the sixties \cite{Hud63, Hae66C} 
and there has been a renewed interest for it in the last years
\cite{CeRe03, CRS04, CrSk08, Sko09}. In this paper the results announced in \cite{CRS07,CRS08} are proved.

This problem generalizes the subject of classical knot theory. In contrast with
knots in $\Bbb R^3$, a complete answer can sometimes be obtained in higher dimensions. We work in the smooth category except when explicitly indicated otherwise. Let us list some known results.

\smallskip
{\bf Knots.} The classification of {\it knots} $S^q\to S^m$ in codimension at least $3$, i. e., for $m>q+2$, has been reduced to a homotopy problem \cite{Hae66A,Le65}. In particular, a complete {\it rational} classification is known:

\proclaim{Theorem~1.1} \cite{Hae66A, Cor.~6.7} Assume
that $q+2 < m < \frac{3}{2}q+2$. Then up to isotopy the set of
smooth embeddings $S^q\to S^m$ is infinite if and only if $q+1$ is
divisible by $4$.
\endproclaim

{\bf Links.} The classification of {\it links} $S^p\sqcup S^q\to S^m$ is the next
natural problem after knots. In codimension at least $3$ there is an
exact sequence involving the set of links up to isotopy and certain
homotopy groups \cite{Hae66C}. In the certain dimension range, called {\it
2-metastable}, there is an explicit description of the isotopy
classes of links $S^p\sqcup S^q\to S^m$ modulo knots $S^p\to S^m$
and $S^q\to S^m$ in terms of homotopy groups of spheres and Stiefel
manifolds \cite{Hae66A}. A short proof of this result was given in \cite{Sko09}.

\centerline{

\beginpicture 
\setcoordinatesystem units <0.600mm,0.600mm>
\setplotarea x from -60 to 60, y from -40 to 20
\put {Figure 1} at 0 -40.000


\setcoordinatesystem units <0.600mm,0.600mm> point at 30 0
\setplotarea x from -30 to 30, y from -35 to 20


\setsolid \setplotsymbol ({\fiverm.}) \circulararc 336.650 degrees
from -7.975 -3.967 center at -0.000 -10.000 \circulararc 336.650
degrees from 7.975 -1.033 center at -0.000 5.000 \put {$S^0 \times
S^1 \to S^3$} at 0 -26.000 \put {a} at 0 -34.000

\setcoordinatesystem units <0.600mm,0.600mm> point at -30 0
\setplotarea x from -30 to 30, y from -35 to 20


\setsolid
\setplotsymbol ({\fiverm.})
\circulararc 3.423 degrees from -8.354 5.471 center at -11.020 -16.690
\circulararc 34.433 degrees from 8.944 -6.708 center at -11.020 -16.690
\circulararc 6.253 degrees from -6.422 10.241 center at -11.020 -16.690
\circulararc 36.084 degrees from 13.416 -4.472 center at -11.020 -16.690
\circulararc 36.084 degrees from -10.581 -9.383 center at -8.944 17.889
\circulararc 6.253 degrees from 12.080 0.441 center at -8.944 17.889
\circulararc 34.433 degrees from -10.282 -4.392 center at -8.944 17.889
\circulararc 3.423 degrees from 8.915 4.500 center at -8.944 17.889
\circulararc 3.423 degrees from -0.560 -9.970 center at 19.964 -1.198
\circulararc 34.433 degrees from 1.337 11.100 center at 19.964 -1.198
\circulararc 6.253 degrees from -5.658 -10.682 center at 19.964 -1.198
\circulararc 36.084 degrees from -2.835 13.855 center at 19.964 -1.198
\circulararc 33.373 degrees from 12.083 -13.748 center at 4.472 -8.944
\circulararc 77.104 degrees from 13.470 -8.770 center at 11.801 -10.531
\circulararc 39.888 degrees from 9.342 -12.449 center at 4.472 -8.944
\circulararc 42.759 degrees from 12.083 -13.748 center at 9.053 -16.600
\circulararc 180.000 degrees from -4.472 -13.416 center at 4.472 -8.944
\circulararc 180.000 degrees from -9.383 10.581 center at -9.982 0.599
\circulararc 180.000 degrees from 13.855 2.835 center at 5.510 8.345
\circulararc 180.000 degrees from 0.000 -11.180 center at 4.472 -8.944
\circulararc 180.000 degrees from -9.682 5.590 center at -9.982 0.599
\circulararc 180.000 degrees from 9.682 5.590 center at 5.510 8.345
\put {$B^2$} at 18.566 -10.728
\put {$S^1 \times S^1 \to S^3$} at 0 -26.000
\put {b} at 0 -34.000
\endpicture
 }
\bigskip

{\bf Knotted tori.} In this paper we study the classification of {\it knotted tori}, i.e.
smooth embeddings $S^p\times S^q\to S^m$. This theory generalizes the theory of
2-component links of the same dimension; see Figure~1. The
investigation of knotted tori is a natural next step (after link theory) towards the classification of embeddings of arbitrary
manifolds \cite{Sko07F, Web67}, by the handle decomposition theorem. This subject is also interesting because of many interesting examples~\cite{Hud63, MiRe71, Sko02}. A systematic investigation of this subject was begun in papers of A.B.~Skopenkov \cite{Sko02, Sko07F, Sko08}.

There exists an explicit description of the set of isotopy classes of knotted tori in the {\it metastable} dimension $m\ge p+\frac{3}{2}q+2$, $p\le q$ \cite{Hae62T, Sko02} (up to small indeterminacy for $m<\frac{3}{2}p+\frac{3}{2}q+2$); see Figure~2. 
This dimension restriction is a natural limit for classical
classification methods for a $(p-1)$-connected $(p+q)$-manifold.
In spite of many interesting approaches \cite{Br68, GW99, Wa70}
little is known below the metastable
dimension: all known explicit classification results concern knots and links (listed above), knotted tori in dimension $m=p+\frac{3}{2}q+\frac{3}{2}$
\cite{Sko08}, 3-manifolds in $\Bbb R^6$ \cite{Sko08Z} and 4-manifolds in $\Bbb
R^7$ \cite{CrSk08}.

Now let us state the main ``practical'' result of the paper announced in \cite{CRS07}. It is an explicit criterion for finiteness of the knotted tori set in the {\it 2-metastable} range (see the dashed domain in Figure~2, in which
the number $p$ is fixed, whereas the numbers $q$ and $m$ vary):

\proclaim{Theorem 1.2} Assume that $p+\frac{4}{3}q+2 < m <p+\frac{3}{2}q+2$ and $m > 2p+q+2$. Then the set of isotopy classes of smooth
embeddings $S^p\times S^q\to S^m$  is infinite if and
only if either $q+1$ or $p+q+1$ is divisible by $4$.
\endproclaim

\proclaim{Example 1} The set of isotopy classes of knotted tori $S^1\times S^5\to
S^{10}$ is finite. All dimensions
in the range $p=1$, $1\le q\le 13$, $m>p+\frac{4}{3}q+2$, such that the set of isotopy classes knotted
tori is infinite, are shown in Figure~3. Information in the figure is obtained from Theorem~1.2, the results in \cite{Sko07L, Th. 3.10 and 2.9}, Proposition~5.6 and Remark~6.3 below.
\endproclaim

Our approach to the classification of embeddings is based on an analogue of the Koschorke exact sequence (Theorem~2.1 below), involving a new {\it
$\beta$-invariant of almost embeddings $S^p\times S^q\to S^m$}.
The exactness is proved using the Habegger--Kaiser techniques of studying the complement of an almost embedding.

{\bf Organization of the paper.} In \S2 we state our main
``theoretical'' result --- Theorem~2.1. In \S3 we define the new {\it
$\beta$-invariant}. In \S4, which is the main part of the paper,
we prove the completeness of the $\beta$-invariant. In \S5 we
prove Theorem~2.1, using the results of \S4. In \S6 we
deduce Theorem~1.2 from Theorem~2.1. In \S7 we state some open
problems. In Appendix we consider the embedded surgery on
self-intersection sets which is used in \S4.

In the publication \cite{Sko10} Theorem~1.2 is generalized to arbitrary $m>2p+q+2$
by reducing the classification of knotted tori to the classification of links. A finiteness criterion for links is required for this purpose; it is obtained in~\cite{CFS11}.

\medskip

\flushpar
$
\matrix
\text{
\beginpicture 
\setcoordinatesystem units <0.700mm,0.700mm>
\setplotarea x from -10 to 80, y from -15 to 105


\setsolid
\setplotsymbol ({\fiverm.})
\arrow <2mm> [.2,.4] from 0.000 0.000 to 0.000 100.000
\arrow <2mm> [.2,.4] from 0.000 0.000 to 75.000 0.000

\plot 0.000 20.000 40.000 100.000 /
\plot 30.000 35.000 70.000 55.000 /
\plot 10.000 30.000 30.000 35.000 /

\setplotsymbol ({\rm.})
\plot 0.000 20.000 70.000 90.000 /

\setshadesymbol ({\rm.})
\setshadegrid span <0.05in>
\vshade 10 30 30 <z,z,,> 30 35 50 <z,z,,> 70 55 90 /

\put {$\bullet$} at 30.000 35.000
\put {0} at 15.000 94.000
\put {1-metastable} at 45.000 81.000
\put {2-metastable} at 51.000 55.000
\put {Unknown} at 53.000 25.000
\put {$m$} at -2.5000 100.000
\put {$0$} at -2.5000 -2.5000
\put {$q$} at 75.000 -2.5000
\put {$p$ fixed} at 35.000 -5.000
\put {Figure 2} at 35.000 -12.500
\endpicture
 } &
\text{
\beginpicture 
\setcoordinatesystem units <0.700mm,0.700mm>
\setplotarea x from -10 to 80, y from -15 to 105


\setsolid
\setplotsymbol ({\fiverm.})
\plot 60.000 0.000 60.000 95.000 /
\plot 65.000 0.000 65.000 95.000 /
\plot 5.000 0.000 5.000 95.000 /
\plot 10.000 -0.000 10.000 95.000 /
\plot 50.000 -0.000 50.000 95.000 /
\plot 55.000 0.000 55.000 95.000 /
\plot 15.000 0.000 15.000 95.000 /
\plot 20.000 -0.000 20.000 95.000 /
\plot 40.000 0.000 40.000 95.000 /
\plot 45.000 0.000 45.000 95.000 /
\plot 25.000 0.000 25.000 95.000 /
\plot 0.000 95.000 65.000 95.000 /
\plot 35.000 0.000 35.000 95.000 /
\plot 30.000 0.000 30.000 95.000 /
\plot -0.000 5.000 65.000 5.000 /
\plot -0.000 85.000 65.000 85.000 /
\plot 0.000 90.000 65.000 90.000 /
\plot 0.000 10.000 65.000 10.000 /
\plot -0.000 15.000 65.000 15.000 /
\plot 0.000 75.000 65.000 75.000 /
\plot 0.000 80.000 65.000 80.000 /
\plot 0.000 20.000 65.000 20.000 /
\plot 0.000 25.000 65.000 25.000 /
\plot 0.000 65.000 65.000 65.000 /
\plot 0.000 70.000 65.000 70.000 /
\plot -0.000 30.000 65.000 30.000 /
\plot 0.000 35.000 65.000 35.000 /
\plot 0.000 55.000 65.000 55.000 /
\plot 0.000 60.000 65.000 60.000 /
\plot 0.000 40.000 65.000 40.000 /
\plot 0.000 45.000 65.000 45.000 /
\plot 0.000 50.000 65.000 50.000 /
\arrow <2mm> [.2,.4] from 0.000 0.000 to 0.000 100.000
\arrow <2mm> [.2,.4] from 0.000 0.000 to 75.000 0.000

\setshadesymbol ({\rm.}) \setshadegrid span <0.03in>
\shaderectangleson

\putrectangle corners at 5.000 5.000 and 10.000 10.000
\putrectangle corners at 10.000 10.000 and 15.000 15.000
\putrectangle corners at 15.000 25.000 and 20.000 30.000
\putrectangle corners at 20.000 30.000 and 25.000 35.000
\putrectangle corners at 25.000 30.000 and 30.000 35.000
\putrectangle corners at 25.000 35.000 and 30.000 40.000
\putrectangle corners at 30.000 40.000 and 35.000 45.000
\putrectangle corners at 25.000 45.000 and 30.000 50.000
\putrectangle corners at 30.000 50.000 and 35.000 55.000
\putrectangle corners at 45.000 60.000 and 50.000 65.000
\putrectangle corners at 45.000 65.000 and 50.000 70.000
\putrectangle corners at 50.000 65.000 and 55.000 70.000
\putrectangle corners at 50.000 70.000 and 55.000 75.000
\putrectangle corners at 40.000 70.000 and 45.000 75.000
\putrectangle corners at 35.000 65.000 and 40.000 70.000
\putrectangle corners at 45.000 85.000 and 50.000 90.000
\putrectangle corners at 50.000 90.000 and 55.000 95.000

\shaderectanglesoff

\setplotsymbol ({\bf.}) \plot 60.000 80.000 65.000 80.000 / \plot
60.000 70.000 60.000 80.000 / \plot 55.000 70.000 60.000 70.000 /
\plot 55.000 65.000 55.000 70.000 / \plot 50.000 65.000 55.000
65.000 / \plot 50.000 60.000 50.000 65.000 / \plot 45.000 60.000
50.000 60.000 / \plot 45.000 50.000 45.000 60.000 / \plot 40.000
50.000 45.000 50.000 / \plot 40.000 45.000 40.000 50.000 / \plot
35.000 45.000 40.000 45.000 / \plot 35.000 40.000 35.000 45.000 /
\plot 30.000 40.000 35.000 40.000 / \plot 30.000 30.000 30.000
40.000 / \plot 25.000 30.000 30.000 30.000 / \plot 25.000 25.000
25.000 30.000 / \plot 20.000 25.000 25.000 25.000 / \plot 20.000
20.000 20.000 25.000 / \plot 15.000 20.000 20.000 20.000 / \plot
15.000 10.000 15.000 20.000 / \plot 10.000 10.000 15.000 10.000 /
\plot 10.000 5.000 10.000 10.000 / \plot 5.000 5.000 10.000 5.000
/ \plot 5.000 0.000 5.000 5.000 / \plot 0.000 0.000 5.000 0.000 /
\put {$m$} at -2.5000 100.000 \put {20} at -2.5000 77.500 \put
{15} at -2.5000 52.500 \put {10} at -2.5000 27.500 \put {5} at
-2.5000 2.500 \put {$q$} at 75.000 -2.5000 \put {10} at 47.500
-2.5000 \put {5} at 22.500 -2.5000 \put {1} at 2.500 -2.5000 \put
{$p=1$} at 35.000 -5.000 \put {Figure 3} at 35.000 -12.500
\endpicture
 }
\endmatrix
$

\bigskip

\head 2. The main idea \endhead

In this section we state our main theoretical result. We present
an 
exact sequence, which reduces the
classification of embeddings $S^p\times S^q\to S^m$ to an easier
classification of all {\it almost embeddings} $S^p\times S^q\to S^m$.

{\bf Isotopies and concordances.} An embedding $f:X\times I\to S^m\times I$ is {\it a concordance} if
$X\times 0=f^{-1}(S^m\times 0)$ and $X\times 1=f^{-1}(S^m\times 1)$.
A concordance is an {\it isotopy} if $f(X\times t)\subset S^m\times t$ for each $t\in I$. A concordance or isotopy is {\it ambient} if $X=S^m$.
We tacitly use the well-known facts that in codimension at least $3$ {\it existence of concordance implies existence of an isotopy} and {\it any concordance or isotopy extends to an ambient one} \cite{Hud69}.

{\bf Almost embeddings.} Informally, an {\it almost
embedding} is a map admitting only ``local'' self-intersections; see
Figure~4. To give a formal definition,
fix a base point $*\in S^p$ and a codimension $0$ ball $B^{p+q}\subset S^p\times S^q$ such
that $B^{p+q}\cap (*\times S^q)=\emptyset$;
see Figure~1b. A
map $f:S^p\times S^q\to S^m$ is an {\it almost embedding},
if the following two conditions hold:

\smallskip
\flushpar (i) $f$ is an embedding outside $B^{p+q}$; and

\flushpar (ii) $fB^{p+q}\cap f(S^p\times S^q-B^{p+q})=\emptyset$.
\smallskip

\flushpar An {\it almost concordance} is defined analogously, except that the
ball $B^{p+q}$ is replaced by $B^{p+q}\times I$.

\bigskip
\flushpar$\matrix
\text{
\beginpicture 
\setcoordinatesystem units <0.80mm,0.80mm>
\setplotarea x from -40 to 40, y from -55 to 25


\setsolid \setplotsymbol ({\fiverm.}) \circulararc 360 degrees
from 35.000 -15.000 center at 20.000 -15.000 \circulararc 360
degrees from 35.000 5.000 center at 20.000 5.000 \circulararc 360
degrees from -12.500 -25.000 center at -20.000 -25.000
\circulararc 360 degrees from -10.000 0.000 center at -20.000
0.000 \setplotsymbol ({\rm.}) \plot -26.364 6.364 -27.778 7.778 /
\plot -13.636 6.364 -12.222 7.778 / \circulararc 90.000 degrees
from -12.929 7.071 center at -20.000 0.000 \circulararc 360
degrees from -15.000 15.000 center at -20.000 15.000 \put
{$\bullet$} at 31.180 -5.000 \put {$\bullet$} at 8.820 -5.000 \put
{$\bullet$} at -20.000 10.000 \put {Almost} at -20.000 -37.000
\put {embedding} at -20.000 -42.000 \put {$f : S^0 \times S^1 \to
S^2$} at -20.000 -48.000 \put {NOT almost} at 20.000 -37.000 \put
{embedding} at 20.000 -42.000 \put {$S^0 \times S^1 \to S^2$} at
20.000 -48.000 \put {Figure 4} at 0.000 -56.000 \put {$fB^1$} at
-30.000 15.000
\endpicture
 } 
\endmatrix$

\smallskip

\flushpar$\matrix
\text{
\beginpicture 
\setcoordinatesystem units <0.600mm,0.600mm>
\setplotarea x from -50 to 50, y from -45 to 40


\setsolid \setplotsymbol ({\fiverm.}) \plot 24.464 16.365 24.464
0.000 / \plot 21.317 17.333 21.317 0.000 / \plot 21.317 -17.333
21.317 0.000 / \plot 24.464 -16.365 24.464 0.000 / \circulararc
44.054 degrees from 21.317 -17.333 center at 40.335 -19.033
\circulararc 51.141 degrees from 24.464 -16.365 center at 40.335
-19.033 \circulararc 51.141 degrees from 35.536 16.365 center at
19.665 19.033 \circulararc 44.054 degrees from 38.683 17.333
center at 19.665 19.033 \plot 35.536 16.365 35.536 0.000 / \plot
38.683 17.333 38.683 0.000 / \plot 38.683 -17.333 38.683 0.000 /
\plot 35.536 -16.365 35.536 0.000 / \circulararc 137.364 degrees
from 37.453 -37.908 center at 38.660 -30.000 \circulararc 137.364
degrees from 15.095 -25.000 center at 21.340 -30.000 \circulararc
137.364 degrees from 37.905 -34.943 center at 38.660 -30.000
\circulararc 137.364 degrees from 17.437 -26.875 center at 21.340
-30.000 \circulararc 16.210 degrees from 32.298 -36.353 center at
40.335 -19.033 \circulararc 11.774 degrees from 34.710 -34.112
center at 40.335 -19.033 \circulararc 62.915 degrees from 22.095
-34.943 center at 19.665 -19.033 \circulararc 60.263 degrees from
22.547 -37.908 center at 19.665 -19.033 \circulararc 60.263
degrees from 44.905 -25.000 center at 30.000 -36.934 \circulararc
62.915 degrees from 42.563 -26.875 center at 30.000 -36.934
\circulararc 16.210 degrees from 19.019 -21.314 center at 30.000
-36.934 \circulararc 11.774 degrees from 19.753 -24.523 center at
30.000 -36.934 \circulararc 11.774 degrees from 40.247 24.523
center at 30.000 36.934 \circulararc 16.210 degrees from 40.981
21.314 center at 30.000 36.934 \circulararc 62.915 degrees from
17.437 26.875 center at 30.000 36.934 \circulararc 60.263 degrees
from 15.095 25.000 center at 30.000 36.934 \circulararc 60.263
degrees from 37.453 37.908 center at 40.335 19.033 \circulararc
62.915 degrees from 37.905 34.943 center at 40.335 19.033
\circulararc 11.774 degrees from 25.290 34.112 center at 19.665
19.033 \circulararc 16.210 degrees from 27.702 36.353 center at
19.665 19.033 \circulararc 137.364 degrees from 42.563 26.875
center at 38.660 30.000 \circulararc 137.364 degrees from 22.095
34.943 center at 21.340 30.000 \circulararc 137.364 degrees from
44.905 25.000 center at 38.660 30.000 \circulararc 137.364 degrees
from 22.547 37.908 center at 21.340 30.000 \circulararc 137.364
degrees from -37.453 12.092 center at -30.000 15.000 \circulararc
137.364 degrees from -37.453 37.908 center at -38.660 30.000
\circulararc 137.364 degrees from -15.095 25.000 center at -21.340
30.000 \circulararc 137.364 degrees from -34.658 13.182 center at
-30.000 15.000 \circulararc 137.364 degrees from -37.905 34.943
center at -38.660 30.000 \circulararc 137.364 degrees from -17.437
26.875 center at -21.340 30.000 \circulararc 60.263 degrees from
-22.547 12.092 center at -40.335 19.033 \circulararc 62.915
degrees from -25.342 13.182 center at -40.335 19.033 \circulararc
16.210 degrees from -32.298 36.353 center at -40.335 19.033
\circulararc 11.774 degrees from -34.710 34.112 center at -40.335
19.033 \circulararc 11.774 degrees from -35.536 16.365 center at
-19.665 19.033 \circulararc 16.210 degrees from -38.683 17.333
center at -19.665 19.033 \circulararc 62.915 degrees from -22.095
34.943 center at -19.665 19.033 \circulararc 60.263 degrees from
-22.547 37.908 center at -19.665 19.033 \circulararc 60.263
degrees from -44.905 25.000 center at -30.000 36.934 \circulararc
62.915 degrees from -42.563 26.875 center at -30.000 36.934
\circulararc 16.210 degrees from -19.019 21.314 center at -30.000
36.934 \circulararc 11.774 degrees from -19.753 24.523 center at
-30.000 36.934 \circulararc 11.774 degrees from -40.247 -24.523
center at -30.000 -36.934 \circulararc 16.210 degrees from -40.981
-21.314 center at -30.000 -36.934 \circulararc 62.915 degrees from
-17.437 -26.875 center at -30.000 -36.934 \circulararc 60.263
degrees from -15.095 -25.000 center at -30.000 -36.934
\circulararc 60.263 degrees from -37.453 -37.908 center at -40.335
-19.033 \circulararc 62.915 degrees from -37.905 -34.943 center at
-40.335 -19.033 \circulararc 16.210 degrees from -21.317 -17.333
center at -40.335 -19.033 \circulararc 11.774 degrees from -24.464
-16.365 center at -40.335 -19.033 \circulararc 11.774 degrees from
-25.290 -34.112 center at -19.665 -19.033 \circulararc 16.210
degrees from -27.702 -36.353 center at -19.665 -19.033
\circulararc 62.915 degrees from -34.658 -13.182 center at -19.665
-19.033 \circulararc 60.263 degrees from -37.453 -12.092 center at
-19.665 -19.033 \circulararc 137.364 degrees from -42.563 -26.875
center at -38.660 -30.000 \circulararc 137.364 degrees from
-22.095 -34.943 center at -21.340 -30.000 \circulararc 137.364
degrees from -25.342 -13.182 center at -30.000 -15.000
\circulararc 137.364 degrees from -44.905 -25.000 center at
-38.660 -30.000 \circulararc 137.364 degrees from -22.547 -37.908
center at -21.340 -30.000 \circulararc 137.364 degrees from
-22.547 -12.092 center at -30.000 -15.000 \setdashes \plot -45.000
-10.000 -36.245 -10.000 / \plot -23.755 -10.000 -15.000 -10.000 /
\plot -23.755 10.000 -15.000 10.000 / \plot -45.000 10.000 -36.245
10.000 / \put {$+$} at -30.000 0.000 \put {$=$} at 0.000 0.000
\put {Figure 5} at 0.000 -45.000
\endpicture
 }& \text{
\beginpicture 
\setcoordinatesystem units <0.700mm,0.700mm>
\setplotarea x from -45 to 40, y from -40 to 35


\setsolid \setplotsymbol ({\fiverm.}) \arrow <2mm> [.2,.4] from
20.000 0.000 to 30.000 0.000 \arrow <2mm> [.2,.4] from -15.000
0.000 to -5.000 0.000 \plot 35.000 30.000 35.000 -30.000 / \plot
14.571 0.000 0.774 -14.001 / \plot 0.697 13.874 14.571 0.000 /
\plot -39.324 13.895 -32.773 7.344 / \circulararc 45.000 degrees
from -39.324 -13.895 center at -37.690 -15.528 \plot -30.519 5.090
-21.464 -3.964 / \plot -40.000 -15.528 -40.000 -30.000 / \plot
-39.902 15.421 -39.902 29.894 / \circulararc 12.223 degrees from
-30.524 -8.555 center at -17.402 -0.106 \circulararc 59.403
degrees from -28.438 10.929 center at -17.402 -0.106 \plot -27.357
-2.132 -39.226 -13.801 / \plot -21.367 3.858 -23.277 1.948 /
\circulararc 180.000 degrees from -21.367 3.858 center at -24.902
7.394 \circulararc 45.000 degrees from -39.902 15.421 center at
-37.592 15.421 \circulararc 180.000 degrees from -28.438 -11.142
center at -24.902 -7.606 \circulararc 45.000 degrees from 0.774
-14.001 center at 2.408 -15.634 \plot 0.000 -15.528 0.000 -30.000
/ \plot 0.000 15.528 0.000 30.000 / \circulararc 45.000 degrees
from 0.000 15.528 center at 2.310 15.528 \put {$\bullet$} at
14.571 0.000 \put {Figure 6} at 0.000 -38.000
\endpicture
 }
\endmatrix$
\bigskip

{\bf Commutative group structure.} The parametric connected
sum operation gives a natural commutative group structure on the set of
embeddings $S^p\times S^q\to S^m$ up to concordance; see Figure~5.
This structure is well-defined for $m > 2p+q+2$ \cite{Sko08}.
We shall give the formal definition in \S5.

\smallskip
{\bf Action of knots on knotted tori.}
For $m>p+q+2$ the set of embeddings
$S^{p+q}\to S^m$ up to concordance is a group with respect to the connected sum operation \cite{Hae66A}.
The same operation gives an action of this group on the set of embeddings
$S^p\times S^q\to S^m$ up to concordance. In \S5 we shall prove that this action is injective for $m>2p+q+2$.
Notice that the set of orbits of the action is in a one-to-one correspondence with the set of embeddings $S^p\times S^q\to S^m$
up to concordance smooth outside a finite set; see Figure~6.

\smallskip
{\bf Notation:}

\flushpar (a) ${E}^m(S^p\times S^q)/E^m(S^{p+q})$ denotes the group of all 
embeddings $S^p\times S^q\to S^m$ up to concordance and connected sums with embeddings
$S^{p+q}\to S^m$.

\flushpar (b) $\overline{E}^m(S^p\times S^q)$ denotes the group of all
almost embeddings $S^p\times S^q\to S^m$ up to almost concordance
(with the parametric connected sum group structure).

\flushpar (c) $\Omega^m_{p,q}:=\pi_{p+2q-m+1}(V_{N+m-p-q-1,N})$,
where $V_{i,j}$ is the {\it
Stiefel manifold} of $j$-frames at the origin of $\Bbb{R}^i$ and $N$ is a large number. Equivalently, for $m\ge p +\frac {4} {3} q+2$ the group
$\Omega_{p,q}^m$ is the {\it normal bordism group} $\Omega_{2p+3q-2m+2}(P^\infty,(m-p-q-1)\lambda)$.
Many of the groups  $\Omega_{p,q}^m$ are known \cite{Kos88L, Pae56}.

Let us state the main ``theoretical'' result of the paper.

\proclaim{Theorem~2.1}
For every $m\ge p+\frac{4}{3}q+2$
and $m > 2p+q+2$ there exists an exact sequence
$$
\multline
{E}^m(S^p\times S^q)/E^m(S^{p+q})\to
\overline{E}^m(S^p\times S^q)\overset{\beta}\to{\to}
\Omega_{p,q}^m\to\\
\to {E}^{m-1}(S^p\times S^{q-1})/E^{m-1}(S^{p+q-1})\to \overline{E} ^{m-1} (S^p\times S^{q-1})\to\Omega _ {p, q-1} ^{m-1}\to\dots
\endmultline
$$
\endproclaim

This theorem has a list of immediate corollaries. First, it
provides estimates of the order or the rank of the group ${E}^m(S^p\times S^q)/E^m(S^{p+q})$ given such estimates for the group
$\overline{E}^m(S^p\times S^q)$. Theorem~2.1 also easily implies the
Haefliger formula for the group of links $S^q\sqcup S^q\to S^m$ in
the 2-metastable range \cite{Hae66C}. A short proof of this
classical result (together with Theorem~2.1 for $p=0$)
was given in \cite{Sko09}. Theorem~2.1 for $p=0$ is analogous to
the Koschorke exact sequence: compare
Theorems~3.1 and~3.5 in \cite{Sko09}, and Theorem~3.1 in \cite{Kos90}.

\smallskip
{\bf The beta-invariant.} The map
$\beta:\overline{E}^m(S^p\times S^q)\to\Omega^m_{p,q}$ in Theorem~2.1 is
a new invariant and it is the main tool of the present paper. It
generalizes both

\flushpar (a) the normal bordism $\beta$-invariant of link maps
\cite{HaKa98, Kos88, Ki90}; and

\flushpar (b) the $\beta$-invariant of knotted tori \cite{Sko07F, Sko08},
cf.~also \cite{Hae66C, Hud63}.

The idea of this invariant is the following. For any almost embedding
$f:S^p\times S^q\to S^m$ we have by definition $fB^{p+q}\cap f(*\times
S^q)=\emptyset$. The $\beta$-invariant measures the ``linking'' of
$f(*\times S^q)$ and (a certain retract of) $fB^{p+q}$.

The main idea of our paper is that from this point of view the
study of almost embeddings $S^p\times S^q\to S^m$ is similar to
the study of link maps $S^q\sqcup S^{p+q}\to S^m$. So we can apply
the full strength of the Habegger--Kaiser technique from link map theory to our problem.
Now let us concentrate on what is done in addition to the technique
of the Habegger--Kaiser paper \cite{HaKa98}.

\smallskip

{\bf Sketch of the proof of Theorem 2.1.} Let us outline the proof
of exactness at the term $\overline{E}^m(S^p\times S^q)$. We need to prove
that any almost embedding $f:S^p\times S^q\to S^m$ such that
$\beta(f)=0$ is almost concordant to an embedding.

It suffices to construct a ball $B^m\subset S^m$ such that
$f^{-1}B^m=B^{p+q}$. Indeed, then the knot $f:\partial B^{p+q}\to \partial B^m$ is trivial by smoothing theory.
Thus one can deform the restriction $f:B^{p+q}\to B^m$ to an embedding
and obtain the desired embedding $S^p\times S^q\to S^m$.

To construct the ball $B^m$ it suffices to span the meridians $f(*\times S^q)$ and $f(S^p\times *)$ by two discs $D^{q+1}$ and $D^{p+1}$ (called {\it webs}),
whose interiors are disjoint and do not intersect the image $f(S^p\times S^q)$. Then the required ball $B^m$ is the complement
of a small tubular neighborhood of $D^{q+1}\cup D^{p+1}$ in
$S^m$.

The existence of the web $D^{p+1}$ is granted by general position and the inequality $m
> 2p+q+2$.
Let us show how to construct the web $D^{q+1}$ under some additional assumptions. This is the most difficult
step, which requires the conditions $\beta(f)=0$ and $m
\ge p+\frac{4}{3}q+2$. By the results of Habegger--Kaiser \cite{HaKa98} we may
assume that $f\left|_{*\times S^q}\right.$ is null-homotopic
outside $fB^{p+q}$. So we can span $f(*\times S^q)$ by a (not
necessarily embedded) disc $D^{q+1}$ in $S^m-fB^{p+q}$. Then we can
remove the self-intersection of the disc using the
embedding theorem.

This completes the proof under the assumption that the interior of $D^{q+1}$ does not intersect also $f(S^p\times S^q-B^{p+q})$ and the disc is smooth. The proof without this assumptions uses an appropriate relative version of this argument; see \S4. Let us emphasize that a {\it relative} version of the above construction is used in the formal proof of Theorem~2.1.

\head 3. The beta-invariant \endhead

In this section we shall give a detailed construction of the
$\beta$-invariant of almost embeddings $S^p\times S^q\to S^m$.

{\bf Idea of the $\beta$-invariant.} The visualization
of this idea is based on an analogy of low-dimensional {\it almost embeddings} $f\sqcup g:S^1\sqcup S^0\to S^2$.

Fix an arc $B^1\subset S^1$. A map $f\sqcup g:S^1\sqcup S^0\to S^2$ is an {\it almost embedding} if it is an embedding outside the arc $B^1$ and $fB^1\cap(f(S^1-B^1)\cup gS^0)=\emptyset$. An {\it almost isotopy} $f_t\sqcup g_t:S^1\sqcup S^0\to S^2$ is a homotopy in the class of almost embeddings.

A simple almost isotopy invariant of an almost
embedding $f\sqcup g:S^1\sqcup S^0\to S^2$ is the linking number $\operatorname{lk}(f,g)$ which assumes
its values in $\Bbb Z_2$. This invariant is incomplete: e.g., the almost embedding in Figure~7 is not almost isotopic to the ``unlink'' but $\operatorname{lk}(f,g)=0$.

\smallskip

\beginpicture 
\setcoordinatesystem units <0.600mm,0.600mm>
\setplotarea x from -50 to 45, y from -35 to 35


\setsolid
\setplotsymbol ({\fiverm.})
\arrow <2mm> [.2,.4] from 5.000 0.000 to 15.000 0.000
\circulararc 360 degrees from -0.000 0.000 center at -20.000 0.000
\circulararc 263.621 degrees from 22.546 -8.333 center at 30.000 -15.000
\plot 22.546 8.333 37.454 -8.333 /
\plot 22.546 -8.333 37.454 8.333 /
\setplotsymbol ({\rm.})
\circulararc 180.000 degrees from -20.000 20.000 center at -20.000 0.000
\circulararc 263.621 degrees from 37.454 8.333 center at 30.000 15.000
\plot 30.000 -0.000 37.454 8.333 /
\plot 22.546 8.333 30.000 -0.000 /
\put {$fC_{\Delta}$} at 30.000 30.000
\put {$\Delta$} at 35.000 0.000
\put {$gS^0$} at 30.000 -20.000
\put {$gS^0$} at 30.000 20.000
\put {$f$} at 10.000 5.000
\put {$C_{\Delta}$} at -45.000 0.000
\put {$f^{-1} \Delta$} at -20.000 -25.000
\put {$f^{-1} \Delta$} at -20.000 25.000
\put {$\bullet$} at -20.000 -20.000
\put {$\bullet$} at -20.000 20.000
\put {$\bullet$} at 30.000 -15.000
\put {$\bullet$} at 30.000 15.000
\put {$\bullet$} at 30.000 0.000
\put {Figure 7} at 10 -33.000
\endpicture

\bigskip

In the situation of Figure~7 the following {\it $\beta$-invariant} of almost isotopy is useful. Take a double point $\Delta$ of
the map $f:S^1\to S^2$. Generically $f^{-1}\Delta$ consists of two points.
Join these two points by an arc $C_\Delta\subset S^1$. The image
$fC_\Delta$ is a cycle, and the number
$\beta(f,g)=\sum_{\Delta}\operatorname{lk}(fC_\Delta,g)\pmod2$,
where the summation is over all double points of $f$, is an almost
isotopy invariant. It is well-defined only if
$\operatorname{lk}(f,g)=0$. For example, for the almost embedding in Figure~7
we have $\beta(f,g)=1$ which proves that the shown almost embedding is indeed not almost isotopic to an ``unlink''.

{\bf Construction.} To realize this idea in higher dimensions we construct:

\flushpar (i) an analogue of the cycle $fC_\Delta$, see the
definition of the cycle $\tilde f$ below; and

\flushpar (ii) a generalization of the linking number
$\operatorname{lk}(fC_\Delta,f(*\times S^q))$, see the definition
of~ $\beta(f)$ below.

\definition{Definition of the double point data} \cite{Kos88L}
Let $f:S^p\times S^q\to S^m$ be an almost embedding.
Make $f$ a general position immersion by an almost isotopy. Consider the
diagram
$$
\CD
\tilde\Delta @>{\tilde i}>> B^{p+q}\\
@V{2:1}VV @VfVV\\
\Delta @>i>> S^m
\endCD
$$
Here $\tilde\Delta=\operatorname{Cl}\{\, (x,y)\in B^{p+q}\times B^{p+q} \,|\,
x\ne y, fx=fy \,\}$ and $\Delta=\tilde\Delta/\Bbb Z_2$ are the
{\it double point manifolds}, where $\Bbb Z_2$ acts on $\tilde\Delta$ by interchanging of factors. The immersions $\tilde i:\tilde
\Delta\to B^{p+q}$ and $i:\Delta\to S^m$ are given by the formulas
$\tilde i(x,y)=x$ and $i\{x,y\}=fx$. Denote by
$\Sigma(f)=\operatorname{Im}\tilde i$ the {\it singular set} of the map
$f$.
\enddefinition

\definition{Definition of the cycle $\tilde f:T(\lambda_\Delta)\to S^m$}
\cite{HaKa98} Denote by $T(\lambda_\Delta)$ the mapping cone
of the double covering $\tilde\Delta\to\Delta$. Take a map $\bar
i:C\tilde\Delta\to B$ extending the map $\tilde i:\tilde\Delta\to B$. Define $\tilde
f:T(\lambda_\Delta)\to S^m$ to be the quotient map of the
composition $f\bar i:C\tilde\Delta\to S^m$.
\enddefinition

Now let us define the generalized linking number. This is the announced step (ii) of the construction of the invariant. Extend
$f\left|_{*\times S^q}\right.$ to a general position immersion
$\bar f:D^{q+1}\to S^m$ ({\it web}). Informally, the desired ``linking number''
is the bordism class of the ``intersection''
$\operatorname{Im}\tilde f\cap\operatorname{Im}\bar f$ with a
natural ``skew framing''. It assumes values in the {\it normal
bordism group} $\Omega^m_{p,q}$ generalizing the group of framed
links.

The construction of all framings below is obvious, so the reader may ignore steps (2) and (3) in all the definitions below. In fact these steps are not used in the paper except in the proof of Proposition~4.4 essentially borrowed from~\cite{HaKa98}.

\definition{Definition of the normal bordism group
$\Omega_s( P^\infty, l\lambda)$ \cite{Kos88L} }
Let $l$ and $s$ be integer numbers.
{\it An
$l\lambda$-manifold} is the triple consisting of

\flushpar(1) a manifold $M$ (``link'');

\flushpar(2) a line bundle $\lambda_M$ on $M$; and

\flushpar(3) a stable isomorphism $\bar g_M:\nu(M)\cong
l\lambda_M:=\underbrace{\lambda_M\oplus\dots\oplus \lambda_M}_{l}$
(``skew framing''), where \topsmash{$\nu(M)$ is the stable normal
bundle of~$M$.}

\flushpar {\it The normal bordism group} $\Omega_s( P^\infty,
l\lambda)$ is the set of $s$-dimensional $l\lambda$-manifolds up
to bordism (with analogous bundle structure). The disjoint union
commutative group structure is evidently defined on this set.
\enddefinition

Hereafter set $s=2p+3q-2m+2$, $l=m-p-q-1$ and $n=p+q$. By
\cite{Kos88L} we have
$\Omega_s(P^\infty,l\lambda)\cong\pi_{l+s}(V_{N+l,N})=\Omega^m_{p,q}$ for $s<l$ and a large integer $N$.
Denote by $N(X,Y)$ the normal bundle of a manifold $X$ immersed into a manifold $Y$ and by
$\epsilon$ the trivial line bundle.

\definition{Definition of the double point $(m-n)\lambda$-manifold}
\cite{Kos88L}
The {\it double point $(m-n)\lambda$-manifold} is a triple $(\Delta,\lambda_\Delta,\bar g_\Delta)$, where:

\flushpar(1) $\Delta$ is the double point manifold;

\flushpar(2) $\lambda_\Delta$ is the line bundle associated with the covering
$\tilde \Delta\to \Delta$; and

\flushpar(3) $\bar g_\Delta:N(\Delta,S^m)\cong
(m-n)\lambda_\Delta\oplus\epsilon^{m-n}$ is constructed as
follows.

For each point $\{x,y\}\in \Delta$ we have canonical isomorphisms
$$
N(\Delta,S^m)_{\{x,y\}}\cong N(\tilde \Delta,B)_{(x,y)}\oplus
N(\tilde \Delta,B)_{(y,x)}\cong N(B,S^m)_y\oplus N(B,S^m)_x.
$$
Let the vectors $\{e^k_x\}_{k=1}^{m-n}$ form a trivialization of
$N(B,S^m)$ at the point $x\in B$. Then the vectors $\{ e^k_x,
e^k_y\}$ form a ``skew framing'' of the manifold $\Delta$. Interchanging of the
points $x$ and $y$ in a pair~$(x,y)\in\tilde\Delta$ implies interchanging of the vectors $e^k_x$
and $e^k_y$. Thus the bundle $N(\Delta,S^m)$ can be decomposed into the sum
of all line bundles $\langle e^k_x+e^k_y\rangle\cong\epsilon$ and
$\langle e^k_x-e^k_y\rangle\cong\lambda_\Delta$ for $k=1,\dots, m-n$. This decomposition defines the
required isomorphism $\bar g_\Delta:N(\Delta,S^m)\cong
(m-n)\lambda_\Delta\oplus \epsilon^{m-n}$.
\enddefinition

Notice that the bundle $\lambda_\Delta$ can be identified with a subset of the mapping cone $T(\lambda_\Delta)$
of the covering $\tilde\Delta\to\Delta$.
Denote the restriction of the cycle $\tilde f:T(\lambda_\Delta)\to S^m$ to the subset still by $\tilde f:\lambda_\Delta\to S^m$. Denote by $\check f:\lambda_\Delta\to S^m$ a general position smooth map (not necessarily an immersion ) sufficiently close to $\tilde f:\lambda_\Delta\to S^m$.

\definition{Definition of the beta-invariant $\beta(f)$}
The beta-invariant $\beta(f)$ of an almost embedding $f:S^p\times S^q\to S^m$ is the bordism class of the $l\lambda$-manifold $(\beta,\lambda_\beta,\bar g_\beta)$ defined as follows:

\flushpar(1) {\it The manifold $\beta=\beta(f)$.} Put
$$
\beta=\{\, (x,y)\in D^{q+1}\times \lambda_\Delta\,:\, \bar fx=\check fy
\,\}.
$$

\flushpar(2) {\it The line bundle $\lambda_\beta$ on the manifold $\beta$.}
Denote by $\operatornamewithlimits{pr}:\beta\to \Delta$ the
obvious composition $\beta\to D^{q+1}\times \lambda_\Delta\to
0\times \Delta=\Delta$. Put
$\lambda_\beta=\operatornamewithlimits{pr^*}(\lambda_\Delta)$.

\flushpar(3) {\it The stable isomorphism $\bar g_\beta:\nu(\beta)\cong
l\lambda_\beta$.} Restrict the stable isomorphism $\bar
g_\Delta:\nu(\Delta)\cong (l+1)\lambda_\Delta$ given by the definition of the
double point $(m-n)\lambda$-manifold, step (3), to the bundle
$<e^1_x-e^1_y>^\perp$. We get an isomorphism
$\nu(\lambda_\Delta)\cong l\lambda_\Delta\oplus \epsilon^{l+1}$.
Identify $\nu(\lambda_\Delta)$ and $\nu(D^{q+1}\times\lambda_\Delta)$.
By restricting the previous isomorphism to the manifold $\beta$ we get an isomorphism
$g_1:\nu(D^{q+1}\times\lambda_\Delta)\left|_{\beta}\right. \cong
l\lambda_\beta\oplus \epsilon^{l+1}$. Take a trivialization of the
normal bundle $N(D^{q+1},S^m)$. The trivialization gives an
isomorphism $g_2:\nu(\beta) \cong
\nu(D^{q+1}\times\lambda_\Delta)\left|_{\beta}\right.\oplus \epsilon^{m-q-1}$. Put $\bar g_\beta=(g_1\oplus \operatorname{id})\circ
g_2$.
\enddefinition

\proclaim{Proposition~3.1} There is a
well-defined map $\overline{E}^m(S^p\times S^q)\to \Omega^m_{p,q}$
given by the formula $f\mapsto\beta(f)$.
\endproclaim

\demo{\bf Proof} We need to check the following:

(1) {\it The class of $\beta(f)$ does not depend on the choices in
the construction.} Indeed, we made the following four choices. In
the definition of the cycle $\tilde f$ we took an extension $\bar
i:C\tilde \Delta\to D^{p+q}$. In the definition of the double point
$(m-n)\lambda$-manifold, step (3), we took a trivialization of the
bundle $N(D^{p+q},S^m)$. In the definition of $\beta(f)$, step
(1), we took an extension $\bar f:D^{q+1}\to S^m$. And in step (3)
of the same definition we took a trivialization of
$N(D^{q+1},S^m)$. Clearly, these extensions and trivializations
are unique up to homotopy, hence the class of $\beta(f)$ is
well-defined.

(2) {\it If $f_1$ and $f_2$ are almost concordant, then
$\beta(f_1)=\beta(f_2)$.} Indeed, let $f:S^p\times S^q\times I\to
S^m\times I$ be a general position almost concordance between $f_1$
and $f_2$. One can construct analogously as above the
$\beta$-invariant $\beta(f)$ of this almost concordance. It will be a
bordism between $\beta(f_1)$ and $\beta(f_2)$. Thus
$\beta(f_1)=\beta(f_2)$. \qed \enddemo

Now we give a relative version of the above construction; see
Figure~8. A map $f:X\to Y$ is {\it proper} if $f^{-1}\partial Y=\partial X$.

\definition{Definition of a proper almost embedding $S^p\times D^q\to D^m$}
Fix a codimension 0 ball $B^{p+q}\subset S^p\times \operatorname{Int} D^q$ such
that $B^{p+q}\cap(*\times D^q)=\emptyset$.
A proper map $f:S^p\times
D^q\to D^m$ is said to be a {\it proper almost embedding} if the
following conditions hold:

\smallskip

\flushpar (i) $f$ is an embedding outside $B^{p+q}$; and

\flushpar (ii) $fB^{p+q}\cap f(S^p\times D^q-B^{p+q})=\emptyset$.

\smallskip

\flushpar A {\it proper almost concordance} is defined analogously, except that the ball $B^{p+q}$ is replaced by $B^{p+q}\times I$.
The {\it standard} embedding $S^p\times D^q\to D^m$ is the composition $S^p\times D^q\subset D^{p+1}\times D^q\cong D^{p+q+1}\subset D^m$.
\enddefinition

Hereafter $f:S^p\times D^q\to D^m$ will be assumed to be a proper almost embedding, unless stated otherwise.

\smallskip

\beginpicture 
\setcoordinatesystem units <0.700mm,0.700mm>
\setplotarea x from -30 to 30, y from -35 to 30


\setsolid \setplotsymbol ({\fiverm.}) \circulararc 13.818 degrees
from -2.502 13.124 center at 12.205 -20.000 \plot -7.352 10.526
-8.325 8.059 / \plot -5.500 11.635 -7.165 8.001 / \plot -3.839
12.443 -5.800 8.288 / \circulararc 120.000 degrees from -10.000
8.660 center at -7.500 12.990 \circulararc 120.000 degrees from
-15.000 -0.000 center at -7.500 12.990 \circulararc 60.000 degrees
from -12.500 8.660 center at -8.750 2.165 \circulararc 120.000
degrees from -8.750 6.495 center at -11.250 6.495 \circulararc
60.000 degrees from -12.500 -0.000 center at -16.250 6.495
\circulararc 120.000 degrees from -16.250 2.165 center at -13.750
2.165 \circulararc 120.000 degrees from 6.250 10.825 center at
5.000 12.990 \circulararc 60.000 degrees from -1.250 10.825 center
at 2.500 17.321 \circulararc 120.000 degrees from -1.250 15.155
center at 0.000 12.990 \circulararc 60.000 degrees from 6.250
15.155 center at 2.500 8.660 \circulararc 360 degrees from 25.000
0.000 center at 0.000 0.000 \circulararc 42.261 degrees from
-7.221 -22.420 center at -3.770 -21.110 \circulararc 42.261
degrees from -2.167 -23.454 center at 0.917 -21.424 \circulararc
42.261 degrees from -11.929 -20.310 center at -8.276 -19.783
\circulararc 42.261 degrees from 8.007 -22.151 center at 9.935
-19.003 \circulararc 42.261 degrees from 2.992 -23.363 center at
5.559 -20.711 \circulararc 42.261 degrees from 24.400 -5.443
center at 20.711 -5.559 \circulararc 42.261 degrees from 22.630
-10.625 center at 19.003 -9.935 \circulararc 42.261 degrees from
22.630 10.625 center at 19.783 8.276 \circulararc 42.261 degrees
from 25.000 0.000 center at 21.424 -0.917 \circulararc 42.261
degrees from 24.400 5.443 center at 21.110 3.770 \circulararc
42.467 degrees from 1.717 16.119 center at -0.431 15.019
\circulararc 106.543 degrees from -1.251 10.825 center at -1.366
10.264 \circulararc 102.259 degrees from 0.915 9.990 center at
0.804 9.131 \circulararc 102.259 degrees from 5.059 9.058 center
at 4.196 9.131 \circulararc 106.543 degrees from 6.937 10.214
center at 6.366 10.264 \circulararc 42.467 degrees from 4.589
17.281 center at 5.431 15.019 \circulararc 72.093 degrees from
8.982 13.426 center at 8.513 12.210 \circulararc 39.567 degrees
from 7.691 16.124 center at 8.082 13.431 \circulararc 39.567
degrees from -1.065 15.258 center at -3.082 13.431 \circulararc
39.567 degrees from -12.619 10.392 center at -10.091 9.384
\circulararc 39.567 degrees from -16.247 2.376 center at -15.673
-0.284 \circulararc 72.093 degrees from -15.035 0.020 center at
-14.830 -1.268 \circulararc 42.467 degrees from -15.601 5.216
center at -15.723 2.806 \circulararc 106.543 degrees from -12.501
-0.001 center at -12.072 -0.382 \circulararc 102.259 degrees from
-10.694 1.457 center at -10.006 0.931 \circulararc 102.259 degrees
from -7.815 4.580 center at -8.310 3.869 \circulararc 106.543
degrees from -7.877 6.785 center at -8.206 6.315 \circulararc
42.467 degrees from -15.171 8.284 center at -12.792 7.883
\circulararc 24.761 degrees from -9.093 -1.925 center at -11.752
2.006 \circulararc 46.460 degrees from -3.879 -1.566 center at
-5.004 -0.186 \circulararc 39.363 degrees from 1.677 1.125 center
at -0.364 2.636 \put {$f(S^p \times D^q)$} at 0.000 -5.000 \put
{$\overline{f}D^{q+1}$} at -10.000 15.000 \put {$D^m$} at 23.000
20.000 \put {Figure 8} at 0.000 -30.000
\endpicture

\bigskip

\definition{Definition of a proper almost embedding $D^{p+q}\to D^m$}
A {\it proper almost embedding} $D^{p+q}\to D^m$ is a proper map whose restriction to the boundary is an embedding.
\enddefinition

\definition{Definition of the web $\bar f$}
Fix the standard equatorial decomposition $\partial
D^{q+1}=D^q_+\cup_{\partial D^q_+=\partial D^q_-=\partial D^q} D^q_-$. A {\it web} of a proper almost embedding
$f:S^p\times D^q\to D^m$ is a map $\bar
f:D^{q+1}\to D^m$ satisfying the following two conditions, see Figure~8:

\smallskip
\flushpar (i) $\bar f\left|_{D^q_+}\right.=f\left|_{*\times
D^q}\right.$; and

\flushpar (ii) $\bar fD^q_-\subset \partial D^m$.
\enddefinition

\definition{Definition of the relative beta-invariants $\beta(f)$ and
$\beta(f,g)$} The definition of the (relative)
{\it $\beta$-invariant $\beta(f)$} of a proper almost embedding
$f:S^p\times D^q\to D^m$ is completely analogous to the definition
of the invariant $\beta(f)$ above, except that the map $f\left|_{*\times S^q}\right.$
is replaced by $f\left|_{*\times D^q}\right.$.

Given a proper map
$g:D^q\to D^m$ missing the image $f(S^p\times D^q)$, define {\it the
$\beta$-invariant $\beta(f,g)$} analogously to $\beta(f)$, only
replace the map $f\left|_{*\times D^q}\right.$ by the map $g$.
\enddefinition

An obvious consequence of the definitions is:

\proclaim{Proposition~3.2} The map $g\mapsto \beta(f,g)$ induces a
homomorphism of groups $\pi_q(D^m-\operatorname{Im}f,\partial
D^m-\operatorname{Im}f)\to\Omega^m_{p,q}$. Moreover, if a proper map $g:D^q\to D^m$
is sufficiently close to $f\left|_{*\times D^q}\right.$, then
$\beta(f,g)=\beta(f)$.
\endproclaim

\head 4. Completeness of the beta-invariant \endhead

In this section we prove the completeness of the
relative $\beta$-invariant:

\proclaim{Theorem~4.1} Assume that $m\ge p+\frac{4}{3}q+2$. Then
any proper almost embedding $f:S^p\times D^q\to D^m$ such that
$\beta(f)=0$ is properly almost concordant to a connected sum (relatively the boundary) of the
standard embedding $S^p\times D^q\to D^m$ and a proper almost embedding $D^{p+q}\to D^m$.
\endproclaim

First let us state our central lemma which describes the
homotopy groups of the complement of a proper almost embedding; cf.~\cite{HaKa98, Corollary~4.4}:

\proclaim{Complement Lemma~4.2}
Assume that
$m\ge p+\frac{4}{3}q+2$. Then any proper almost embedding $f:S^p\times D^q\to D^m$ is
properly almost concordant to a general position proper almost embedding $f':S^p\times D^q\to D^m$ such
that
$$
\pi_{q}(D^m-\operatorname{Im}f',\partial D^m-\operatorname{Im}f') \cong\Omega^m_{p,q}.
$$
This isomorphism is given by the formula $g\mapsto\beta(f',g)$.
\endproclaim

In particular, if $\beta(f)=0$ then any proper map $g:D^q\to D^m-\operatorname{Im}f'$ sufficiently close to
$f'\left|_{*\times D^q}\right.$ is properly null-homotopic (outside $\operatorname{Im}f'$). This remark follows from Proposition~3.2 above and forms the basis of the following argument.

\demo{Proof of 4.1 modulo Lemma 4.2} The proof consists of $3$
steps:

(1) {\it Construction of a web,
whose interior misses
$\operatorname{Im}f$.} Take a proper almost embedding $f:S^p\times D^q\to D^m$.
Let $f':S^p\times D^q\to D^m$ be a proper almost embedding (properly almost concordant to $f$)
such that the isomorphism of Lemma~4.2 holds.

Without loss of generality assume that
$f'(S^p\times D^q)$ is orthogonal to $\partial D^m$. Since the restriction of the normal
bundle $N(f'(S^p\times D^q),D^m)$ to the disc $f'(*\times D^q)$ is
trivial, there exists a unit vector field on
$f'(*\times D^q)$ orthogonal to $f'(S^p\times D^q)$.
Attach a collar neighborhood to $f'(*\times D^q)$
along this vector field.
One boundary component of this collar
forms a proper map $g:D^q\to D^m-\operatorname{Im}f'$.

By Lemma~4.2 (and the paragraph after its statement) the map $g:D^q\to D^m-\operatorname {Im} f '$ is properly null-homotopic.
Hence one can attach a disc (possibly with self-intersections, missing $\operatorname{Im}f'$ and with boundary in $\operatorname {Im}g\cup \partial D^m$)
to the boundary component of the collar. The union of the disc and the collar is the image of the desired web $\bar f:D^{q+1}\to D^m$.

(2) {\it Removing self-intersection of the web.} By \cite{Ha67, Theorem~2.1} there exists
a piecewise smooth homeomorphism $h:D^m\to D^m$ such that $f'_{PL}:=hf'$
and $\bar f_{PL}:=h\bar f$ are piecewise linear. Denote by $M^m$
the complement of a regular neighborhood of $\operatorname{Im}f'_{PL}$
in $D^m$. Then the pair $(M^m,\partial M^m\cap\partial D^m)$ is
sufficiently highly connected (see Proposition~4.6 below). The
following theorem allows to remove the self-intersection of the
web:

\proclaim{Embedding Theorem moving a part of the boundary~4.3}
Let $(M^m,M^{m-1})$ be a pair of
piecewise linear manifolds such that $M^{m-1}\subset \partial M^m$. Suppose this pair is
$(2q-m+3)$-connected and $m \ge q+4$. Let $\bar
f_{PL}:(D^{q+1},D^{q}_-)\to (M^m,M^{m-1})$ be a piecewise linear map
which embeds $D^{q}_+$ into $\partial M^m-\operatorname{Int} M^{m-1}$. Then
$\bar f_{PL}$ is homotopic $\operatorname{rel}D^{q}_+$ to a piecewise linear embedding
$\bar f_{Emb}:(D^{q+1},D^{q}_-)\to (M^m,M^{m-1})$.
\endproclaim

This theorem is proved completely analogously to the result \cite{Hud69,
Theorem~9.2.1}.
Since all the obstructions to smoothing the embedding $h ^ {-1} \bar f _ {Emb} $ belong to trivial groups $H^k(D^{q+1},D^{q}_-;C^{m-q}_{q-k})$ it follows
by smoothing theory \cite{Ha67} that the piecewise linear embedding $h^{-1}\bar f_{Emb}$ is properly homotopic $\operatorname{rel}D^{q}_+$
to a web $\bar f':D^{q+1}\to D^m$ which is a smooth embedding.

(3) {\it Decomposition of $f'$ into a connected sum.} Let $B^m$ be the complement in
$D^m$ of the union of tubular neighborhoods of $f'(S^p\times D^q-\operatorname{Int}B^{p+q})$ and $\bar f'D^{q+1}$.
Then $B^m$ is a smooth ball such that $(f')^{-1}B^m=B^{p+q}$.
Denote the restriction $f':B^{p+q}\to B^m$ by $g':B^{p+q}\to B^m$. It is easy to see that
$f'$ is properly almost concordant to a connected sum (relatively the boundary) of $g':B^{p+q}\to B^m$ and an embedding $S^p\times D^q\to D^m$.
The latter is properly ambient concordant to the standard one (by easy Proposition~5.7(a) below). Thus Theorem~4.1 modulo Lemma~4.2 follows.
\qed
\enddemo

\medskip

The rest of \S4 is devoted to the proof of Lemma~4.2. Our argument is parallel
to \cite{HaKa98, \S\S3--4}.

\smallskip
{\bf Notation:} Hereafter we shall identify the disc $D^m$ with the upper half-sphere of $S^m$
and we shall fix the decomposition $S^m=D^m\cup CS^{m-1}$.
Denote by $\{X,Y\}$ the set of stable homotopy classes of maps $X\to Y$.

\proclaim{Example}
$\{ (D^q; S^{q-1})\,,\,(D^m-\operatorname{Im}f,\partial)\}
\cong\{S^q,S^m-f(S^p\times D^q)\}$ (because a pair
$(X,Y)$ is stably homotopy equivalent to the space $X\cup CY$).
\endproclaim

The following proposition shows how to express the $\beta$-invariant via
the homotopy class of a map $g:D^q\to D^m-\operatorname {Im} f$ in the group $\pi_q(D^m-\operatorname{Im}f,\partial)$.
This can be considered as an alternative definition of the
$\beta$-invariant $\beta(f,g)$.

\proclaim{Proposition~4.4} (cf. \cite{HaKa98, Proposition~3.2})
Under the composition
$$
\multline
\pi_q(D^m-\operatorname{Im}f,\partial)
\overset{\Sigma^\infty}\to{\to}
\{S^q,S^m-\operatorname{Im}f\}
\overset{SW}\to{\to}\\
\overset{SW}\to{\to}
\{\operatorname{Im}f,S^{m-q-1}\}
\overset{\tilde f_*}\to{\to}
\{T(\lambda_\Delta),S^{m-q-1}\}
\overset{PT}\to{\to}
\Omega_s(P^{\infty},l\lambda)
\endmultline
$$
the homotopy class
of a proper map $g:D^q\to D^m-\operatorname {Im} f$
is sent to $\beta(f,g)$.
\endproclaim

Here the first arrow is the iterated suspension map. The second map is the
Spanier-Whitehead duality. The third arrow is induced by the map
$\tilde f$ defined in \S3. The fourth map is given by the
Pontryagin--Thom construction (see \cite{Kos88L} for details). This
proposition is proved by a direct verification. In fact, it
suffices to prove it for a map $f:B^{p+q}\to S^m$, which can be done
analogously to \cite{HaKa98, Proposition~3.2}.

This assertion suggests to find the homotopy type of
$\operatorname{Im}f$:

\proclaim{Proposition~4.5} (cf. \cite{HaKa98, \S4}) Denote by $C$ the
mapping cone of $f:\Sigma(f)\to f\Sigma(f)$. Then
$\operatorname{Im}f\simeq C\vee S^p$.
\endproclaim

The proof of this proposition immediately follows from the observation that both of these spaces can be
obtained from $\operatorname{Cyl}(\Sigma(f)\to f\Sigma(f))
\cup_{\Sigma(f)\subset S^p\times D^q}(S^p\times D^q)$ by
appropriate contractions; see Figure~9.

\beginpicture 
\setcoordinatesystem units <0.600mm,0.600mm>
\setplotarea x from -60 to 60, y from -30 to 30


\arrow <2mm> [.2,.4] from -5.000 0.000 to 5.000 0.000
\put {$f$} at 0.000 5.000
\put {Figure 9} at 0.000 -25.000

\setcoordinatesystem units <0.600mm,0.600mm> point at 30 0
\setplotarea x from -30 to 30, y from -30 to 30


\setsolid
\setplotsymbol ({\fiverm.})
\circulararc 360 degrees from 20.000 0.000 center at 0.000 0.000
\setplotsymbol ({\rm.})
\plot -16.454 9.500 -18.187 10.500 /
\plot 16.454 9.500 18.187 10.500 /
\circulararc 120.000 degrees from 17.321 10.000 center at -0.000 0.000
\setdashes
\setplotsymbol ({\fiverm.})
\plot 0.000 0.000 -14.142 14.142 /
\plot 0.000 0.000 0.000 20.000 /
\plot 0.000 0.000 14.142 14.142 /
\put {$B^{p+q}$} at -20.000 18.000
\put {$\Sigma(f)$} at 0.000 25.000
\put {Cyl} at 0.000 -5.000
\put {$\bullet$} at -14.142 14.142
\put {$\bullet$} at 14.142 14.142
\put {$\bullet$} at 0.000 20.000
\put {$\bullet$} at 0.000 0.000

\setcoordinatesystem units <0.600mm,0.600mm> point at -30 0
\setplotarea x from -30 to 30, y from -30 to 30


\setsolid
\setplotsymbol ({\fiverm.})
\plot -4.330 7.500 4.330 -7.500 /
\plot 4.330 7.500 -4.330 -7.500 /
\plot -8.660 0.000 8.660 0.000 /
\circulararc 240.000 degrees from -8.660 -0.000 center at -8.660 -5.000
\circulararc 240.000 degrees from 4.330 -7.500 center at 8.660 -5.000
\circulararc 240.000 degrees from 4.330 7.500 center at 0.000 10.000
\put {Im $f$} at 0.000 -15.000
\put {$\bullet$} at 0.000 -0.000

\endpicture

\bigskip

Let us begin the study of the homotopy type of the pair
$(D^m-\operatorname{Im}f,\partial)$.

\proclaim{Proposition~4.6} (cf. \cite{HaKa98, Lemma~4.2}) For a general position immersion $f$ the pair
$(D^m-\operatorname{Im}f,\partial)$ is $c$-connected, where
$c=\min\{ m-p-2,2m-2p-2q-3\}$.
\endproclaim

\demo{Proof} In codimension at least 3 the pair
$(D^m-\operatorname{Im}f,\partial)$ is generically 1-connected. By the
homology excision theorem $H_i(D^m-\operatorname{Im}f,\partial)
\cong H_i(S^m-\operatorname{Im}f)$. By the Alexander duality
$H_i(S^m-\operatorname{Im}f)\cong H^{m-i-1}(\operatorname{Im}f)$.
By Proposition~4.5 we have $H^{m-i-1}(\operatorname{Im}f)=0$ for $i\le c$,
because generically $\operatorname{dim}C=2p+2q-m+1$. So by the
Hurewicz theorem the pair $(D^m-\operatorname{Im}f,\partial)$ is
$c$-connected. \qed
\enddemo

\demo{Proof of~Lemma~4.2} It suffices to show that for an
appropriate choice of proper almost embedding $f'$ all maps in Proposition~4.4 are bijective. We need to modify the initial proper almost embedding $f$ only in step (4) below.

(1) {\it The first map is bijective} by Proposition~4.6 and the
Suspension theorem, because the assumption $m\ge p+\frac{4}{3}q+2$
implies $q/2\le c$.

(2) {\it The second map is bijective} by the Spanier-Whitehead
duality.

(3) {\it The third map is bijective.} By Proposition~4.5
$\{\operatorname{Im}f,S^{m-q-1}\}\allowmathbreak
\cong\{C,S^{m-q-1}\}$, hence it remains to check that
$\{C,S^{m-q-1}\}\cong\{T(\lambda_\Delta),S^{m-q-1}\}$. Denote by
$\lambda_T$ the restriction of $\lambda_\Delta$ to the triple
points, by $C_T$ the image of $T(\lambda_T)$ under the map
$T(\lambda_\Delta)\to C$, and by $R$ the mapping cone of the map
$T(\lambda_T)\to C_T$. Then $R$ is a deformation retract of the
mapping cone of the map $T(\lambda_\Delta)\to C$. Consider the Puppe
exact sequence of the pair
$(\operatorname{Cyl}(T(\lambda_\Delta)\to C), T(\lambda_\Delta))$:
$$
\{R,S^{m-q-1}\}\to\{C,S^{m-q-1}\}
\to\{T(\lambda_\Delta),S^{m-q-1}\}\to\{R,S^{m-q}\}.
$$
Since $R$ has dimension at most $3q+3p-2m+2$, it follows by the assumption $m\ge
p+\frac{4}{3}q+2$ that $\{C,S^{m-q-1}\}\cong
\{T(\lambda_\Delta),S^{m-q-1}\}$.
Clearly, the obtained isomorphism
$\{\operatorname {Im} f, S^{m-q-1} \} \cong \{T (\lambda_\Delta), S ^ {m-q-1} \} $
is induced by the map $\tilde f:T (\lambda_\Delta)\to\operatorname {Im} f$.

(4) {\it Making the fourth map bijective.} One can see that the fourth
map is bijective if the map $\Delta\to P^\infty$ classifying the
bundle $\lambda_\Delta$ is $(s+1)$-connected. Indeed, by the
Thom-Pontryagin construction $\{T(\lambda_\Delta),S^{m-q-1}\}\cong
\Omega_s(\Delta,l\lambda_\Delta)$ and by homotopy lifting the
induced map
$\Omega_s(\Delta,l\lambda_\Delta)\to\Omega_s(P^\infty,l\lambda)$
is an isomorphism. (The definition of the group
$\Omega_s(\Delta,l\lambda_\Delta)$ and the details can be found in
\cite{HaKa98,\S3} and \cite{Kos88L}).

So it remains to make the classifying map $\Delta\to P^\infty$ $(s+1)$-connected by a proper almost concordance of $f$. This is made by the following theorem (applied here for $s=2p+3q-2m+2$ and $n=p+q$):

\proclaim{Surgery Theorem 4.7} (cf. \cite{HaKa98, Theorem~4.5}) Let
$M^m$ be an $(s+1)$-connected manifold
and let $f:B^n\to M^m$ be a proper immersion such that
the restriction $f\left|_{\partial B^n}\right.$ is an embedding. Assume that $2s \le 2n-m-2$ and $0\le s
\le m-n-3$. Then by a regular
homotopy  $\operatorname{rel}\partial B^n$
of the immersion $f:B^n\to M^m$
the classifying map
$\Delta\to P^\infty$ of the covering $\tilde\Delta\to\Delta$ can be made $(s+1)$-connected.
\endproclaim

The proof of Theorem~4.7 is completely analogous to the
proof of \cite{HaKa98, Theorem~4.5}. We present it in Appendix for the convenience of the reader.
\qed
\enddemo

Thus we have proved,
modulo Theorem~4.7, the completeness of the $\beta$-invariant.

\head 5. The exact sequence \endhead

In this section we deduce Theorem~2.1 from the completeness of the
$\beta$-invariant. Formally, Theorem~2.1 follows from assertions~4.1,
5.1, 5.2 and~5.5.

\proclaim{Theorem~5.1} Assume that $m\ge p+\frac{4}{3}q+2$. Then
for each $x\in\Omega^m_{p,q}$ there exists a proper almost
embedding $\omega_x:S^p\times D^q\to D^m$ such that $\beta(\omega_x)=x$.
\endproclaim

\demo{Proof of~5.1 modulo Lemma~4.2} The construction of
$\omega_x$ proceeds in $3$ steps:

(1) {\it Construction of a map $f'\sqcup g:S^p\times D^q\sqcup
D^q\to D^m$ such that $\beta(f',g)=x$.} Start with the standard
embedding $f:S^p\times D^q\to D^m$. By Lemma~4.2 from \S4
the map $f$ is properly almost concordant to a proper almost
embedding $f'S^p\times D^q\to D^m$, such that
$\beta(f',{\cdot}):\pi_q(D^m-\operatorname{Im}f',\partial)\to
\Omega^m_{p,q}$ is an isomorphism.
Take a proper map $g:D^q\to
D^m-\operatorname{Im}f'$ such that $\beta(f',g)=x$.

(2) {\it Homotopical extension of $g:D^q\to D^m-\operatorname{Im}f'$ to a proper embedding $g':S^p\times D^q\to
D^m-\operatorname{Im}f'$.} Since
$(D^m-\operatorname{Im}f',\partial)$ is sufficiently highly connected (see
Proposition~4.6),  it follows by the embedding theorem \cite{Hud69,
Theorem~8.2.1} that one can remove the self-intersection of $g:D^q\to
D^m-\operatorname {Im} f '$. By
Hirsch theory one can make $g:D^q\to
D^m-\operatorname {Im} f '$ a smooth embedding with a
trivial normal bundle. So one can extend $g:D^q\to D^m-\operatorname{Im}f'$ to a proper embedding
$g':S^p\times D^q \to D^m-\operatorname{Im}f'$ such that the image
$\operatorname{Im}g'$ is contained in a tubular neighborhood of
$\operatorname{Im}g$.

(3) {\it Making the $S^p$-parametric connected sum of $f'$ and $g'$.}
Fix a point $*\in \partial D^{q}$.
By construction it follows that
$f'\left|_{\partial}\right.$ is concordant to the standard embedding.

Thus the sphere $f'(S^p\times *)$ can be spanned by a framed disc
$D^{p+1}\subset \partial D^m$ ({\it web}) such that $\operatorname{Int} D^{p+1}\cap\operatorname{Im}f'=\emptyset$
and the first $q$ vector fields of the framing of
$\partial D^{p+1}$ are tangent to $f'(S^p\times S^q)$. Generically the web does not intersect
also $\operatorname{Im}g$, and hence neither $\operatorname{Im}g'$.
Clearly, the sphere $g'(S^p\times *)$ can also be spanned by
a web $\bar D^{p+1}\subset \partial D^m$ whose interior does not intersect
$\operatorname{Im}(f'\sqcup g')$ and $D^{p+1}$.

Join the centers of the webs $D^{p+1}$ and $\bar D^{p+1}$ by an
arc $I$ in $\partial D^m$. Generically $I$ intersects
$\operatorname{Im}(f'\sqcup g')\cup D^{p+1}\cup \bar D^{p+1}$ only
at the boundary $\partial I$. Let $\bar D^m$ be the union of small tubular neighborhoods of
$D^{p+1}$, $I$ and $\bar D^{p+1}$ in the ball $D^m$. Clearly,
the intersection $\operatorname{Im}(f'\sqcup g')\cap \bar D^m$ is
standard. By performing the $S^p$-parametric connected sum of $f':S^p\times D^q\to
D^m$
and $g':S^p\times D^q\to
D^m$ in $\bar D^{m}$, we obtain the required proper almost
embedding $\omega_x:S^p\times D^q\to D^m$. \qed\enddemo

{\bf Notation:} $\overline{E}^m(S^p\times D^q,S^p\times S^{q-1})$ denotes the group
of proper almost embeddings $S^p\times D^q\to D^m$ up to proper
almost concordance and connected sums with proper almost embeddings $D^{p+q}\to D^m$
(with the parametric connected sum group structure).

{\it Remark.} The relative $\beta$-invariant is a
map $\overline{E}^m(S^p\times D^q,S^p\times S^{q-1})\to \Omega^m_{p,q}$.
Theorems~4.1 and~5.1 together assert that this map is bijective for
$m\ge p+\frac{4}{3}q+2$.

\flushpar
\beginpicture 
\setcoordinatesystem units <0.850mm,0.850mm>
\setplotarea x from -80 to 70, y from -30 to 30
\put {Figure 10} at 0.000 -29.000


\arrow <2mm> [.2,.4] from -43.000 0.000 to -37.000 0.000
\put {$e$} at -40.000 5.000
\arrow <2mm> [.2,.4] from -3.000 0.000 to 3.000 0.000
\put {$h$} at 0.000 5.000
\arrow <2mm> [.2,.4] from 32.000 0.000 to 38.000 0.000
\put {$p$} at 35.000 5.000

\setcoordinatesystem units <0.850mm,0.850mm> point at 60 0
\setplotarea x from -20 to 20, y from -25 to 25


\setsolid
\setplotsymbol ({\fiverm.})
\circulararc 120.000 degrees from -1.121 6.424 center at -3.536 7.071
\circulararc 60.000 degrees from -4.182 9.486 center at -2.241 2.241
\circulararc 120.000 degrees from -9.486 4.183 center at -7.071 3.536
\circulararc 60.000 degrees from -6.424 1.121 center at -8.365 8.365
\circulararc 57.312 degrees from -2.433 -13.738 center at -0.428 -12.776
\circulararc 57.312 degrees from 6.452 -12.370 center at 7.455 -10.385
\circulararc 57.312 degrees from 1.825 -13.832 center at 3.447 -12.310
\circulararc 57.312 degrees from 14.497 -3.852 center at 12.310 -3.447
\circulararc 57.312 degrees from 12.301 -8.584 center at 10.385 -7.455
\circulararc 57.312 degrees from 14.984 0.702 center at 12.776 0.428
\circulararc 360 degrees from 15.000 0.000 center at 0.000 0.000
\circulararc 60.000 degrees from 6.424 -1.121 center at 8.365 -8.365
\circulararc 120.000 degrees from 1.121 -6.424 center at 3.536 -7.071
\circulararc 60.000 degrees from 4.183 -9.486 center at 2.241 -2.241
\circulararc 120.000 degrees from 9.486 -4.183 center at 7.071 -3.536

\put {$S^1 \sqcup S^1 \to S^2$} at 0.000 20.000
\put {$\frac{{E}^2(S^0\times S^1)}{E^2(S^1)}$} at 0.000 -20.000

\setcoordinatesystem units <0.850mm,0.850mm> point at 20 0
\setplotarea x from -20 to 20, y from -25 to 25


\setsolid
\setplotsymbol ({\fiverm.})
\circulararc 57.312 degrees from -2.433 -13.738 center at -0.428 -12.776
\circulararc 57.312 degrees from 6.452 -12.370 center at 7.455 -10.385
\circulararc 57.312 degrees from 1.825 -13.832 center at 3.447 -12.310
\circulararc 57.312 degrees from 14.497 -3.852 center at 12.310 -3.447
\circulararc 57.312 degrees from 12.301 -8.584 center at 10.385 -7.455
\circulararc 57.312 degrees from 14.984 0.702 center at 12.776 0.428
\circulararc 360 degrees from 15.000 0.000 center at 0.000 0.000
\plot -2.500 5.159 -4.183 -1.121 /
\plot 1.121 4.183 -5.314 2.458 /
\circulararc 60.000 degrees from 6.424 -1.121 center at 8.365 -8.365
\circulararc 120.000 degrees from 1.121 -6.424 center at 3.536 -7.071
\circulararc 60.000 degrees from 4.183 -9.486 center at 2.241 -2.241
\circulararc 120.000 degrees from 9.486 -4.183 center at 7.071 -3.536
\circulararc 60.000 degrees from -4.022 6.674 center at -3.052 3.052
\circulararc 120.000 degrees from -6.674 4.022 center at -5.467 3.699
\circulararc 120.000 degrees from -2.492 5.143 center at -3.699 5.467
\circulararc 120.000 degrees from 4.183 1.121 center at 1.768 1.768
\circulararc 60.000 degrees from -1.121 -4.183 center at -3.062 3.062
\circulararc 120.000 degrees from -4.183 -1.121 center at -1.768 -1.768
\setdashes
\circulararc 120.004 degrees from 8.365 2.241 center at 3.536 3.536
\circulararc 59.990 degrees from -2.241 -8.365 center at -6.125 6.125
\circulararc 120.004 degrees from -8.365 -2.242 center at -3.536 -3.536
\circulararc 60.001 degrees from 2.242 8.365 center at 6.124 -6.124
\put {$\bullet$} at -3.062 3.062
\put {$S^1 \sqcup S^1 \to S^2$} at 0.000 20.000
\put {\eightpoint $\overline{E}^2(S^0\times S^1)$} at 0.000 -20.000

\setcoordinatesystem units <0.680mm,0.680mm> point at -22 0
\setplotarea x from -20 to 20, y from -30 to 30


\setsolid \setplotsymbol ({\fiverm.}) \circulararc 93.135 degrees
from 7.262 -13.125 center at 0.000 -20.000 \circulararc 70.529
degrees from 5.000 0.000 center at -10.000 0.000 \circulararc
180.000 degrees from -5.000 -0.000 center at -0.000 0.000
\circulararc 70.529 degrees from 5.000 14.142 center at 10.000
0.000 \circulararc 360 degrees from 15.000 0.000 center at 0.000
0.000 \put {$\bullet$} at -5.000 14.142 \put {$\bullet$} at 5.000
14.142 \put {$\bullet$} at -7.262 -13.125 \put {$\bullet$} at
7.262 -13.125 \put {$D^1 \sqcup D^1 \to D^2$} at 0.000 25.000 \put
{\eightpoint $\tsize\overline{E}^2(S^0\times D^1,S^0\times S^{0})$} at 0.000 -25.000

\setcoordinatesystem units <0.680mm,0.680mm> point at -67 0
\setplotarea x from -20 to 20, y from -30 to 30


\setsolid
\setplotsymbol ({\fiverm.})
\circulararc 360 degrees from 15.000 0.000 center at 0.000 0.000
\put {$\bullet$} at -5.000 14.142
\put {$\bullet$} at 5.000 14.142
\put {$\bullet$} at -7.262 -13.125
\put {$\bullet$} at 7.262 -13.125
\put {$S^0 \sqcup S^0 \to S^1$} at 0.000 25.000
\put {$\frac{{E}^1(S^0\times S^0)}{E^{1}(S^{0})}$} at 0.000 -25.000

\endpicture

\bigskip

\proclaim{Lemma~5.2} (cf. \cite{Ne84, Sko09, Th. 3.1}, see Figure~10) For every $m
> 2p+q+2$ there is an exact sequence
$$
\dfrac{{E}^m(S^p\times S^q)}{E^{m}(S^{p+q})}
\overset{e}\to{\to} \overline{E}^m(S^p\times S^q)\overset{h}\to{\to}
\overline{E}^m(S^p\times D^q,S^p\times S^{q-1})\overset{p}\to{\to}
\dfrac{{E}^{m-1}(S^p\times S^{q-1})}{E^{m-1}(S^{p+q-1})}\to\dots
$$
\endproclaim

In the rest of \S5 we shall prove the Lemma~5.2.
We begin with the construction of the commutative group structure
on the set of knotted tori
(and also almost embeddings $S^p\times S^q\to S^m$).
This construction is equivalent to the one from \cite{Sko08, \S2}.

\definition{Definition of the web $D^{p+1}$} Fix a point $*\in S^q$.
Without loss of generality assume that $*\times S^q\cap B^{p+q}=\emptyset$.
A {\it web} of an almost embedding $f:S^p\times S^q\to S^m$ is a
framed disc $D^{p+1}\subset S^m$ satisfying the following $3$ conditions:

\smallskip

\flushpar (i) $\partial D^{p+1}=f(S^p\times *)$;

\flushpar (ii) $\operatorname{Int}
D^{p+1}\cap\operatorname{Im}f=\emptyset$; and

\flushpar (iii) the first $q$ vector fields of the framing of $\partial D^{p+1}$ coincide with the obvious framing of the ``meridian'' $f(S^p\times *)$ in the ``torus'' $f(S^p\times S^q)$.
%
\flushpar {\it Webs}
of an almost concordance $f:S^p\times S^q\times I\to S^m\times I$ and
of a proper almost embedding $f:S^p\times D^q\to D^m$ are
defined analogously (the marked point $*\in D^q$ is the center of the disk
$D^q$).
\enddefinition

The following assertion is equivalent to \cite{Sko08, Standardization Lemma 2.1}.

\proclaim{Proposition~5.3} \cite{Sko08, Standardization Lemma 2.1} If $m>2p+q+1$ then for any
almost embedding $f:S^p\times S^q\to S^m$ there exist a web.
If $m>2p+q+2$, then for any almost concordance between the almost embeddings $f_1,f_2:S^p\times S^q\to
S^m$ there exists a web extending given webs of $f_1$ and $f_2$.
\endproclaim

\demo{Proof} \cite{Sko08} The bundle $N(f|_{S^p\times *};S^m)$
is stably trivial and $m-p-q\ge p$, hence this bundle is trivial.
Take a $(m-p-q)$-framing $\xi$ of this bundle.

Take the section formed by the first vectors of $\xi$. Since
$m\ge2p+q+2\ge2p+2$, it follows that the embedding $f|_{S^p\times *}$ is unknotted
in $S^m$. So there is an embedding $\bar f:D^{p+1}\subset S^m$
satisfying property (i) from the definition of the web above. Since
$m\ge2p+q+2$, we may also assume property (ii) by general position.

By deleting the first vector from $\xi$ we obtain a
$(m-p-q-1)$-framing $\xi_1$ of the boundary $\bar f(\partial D^{p+1})$ orthogonal to the disc
$\bar f(D^{p+1})$. Denote by $\eta$ the standard normal $q$-framing
of $f(S^p\times *)$ in $f(S^p\times S^q)$. Then $(\xi_1,\eta)$ is a
normal $(m-p-1)$-framing on $\bar f(\partial D^{p+1})$ orthogonal to
$\bar f(D^{p+1})$. Since $p<m-p-q-1$, the map
$\pi_p(SO_{m-p-q-1})\to\pi_p(SO_{m-p-1})$ is epimorphic. Hence we
can change $\xi_1$ (and thus $\xi$) so that $(\xi_1,\eta)$ extends
to a normal framing on $\bar f(D^{p+1})$. Clearly, the obtained framing
satisfies property (iii).

The second assertion of the proposition is proved analogously
\cite{Sko08, Proof of Standardization Lemma 2.1 in \S3}
\qed\enddemo

The following definition is equivalent to the one in \cite{Sko08}.

\definition{Definition of the $S^p$-parametric connected sum} (see Figure~5)
Let $f_1,f_2:S^p\times S^q\to S^m$
be a pair of almost embeddings. Without loss of generality we may assume
that $\operatorname{Im}f_1\subset D^m_+$, $\operatorname{Im}f_2\subset D^m_-$.
Take webs $D^{p+1}\subset D^m_+$ and $\bar D^{p+1}\subset D^m_-$ of these almost embeddings.
Join the centers of the webs $D^{p+1}$ and $\bar D^{p+1}$ by an
arc $I$ in $S^m$. Generically $I$ intersects
$\operatorname{Im}f_1\cup\operatorname{Im}f_2\cup D^{p+1}\cup \bar D^{p+1}$ only
at the boundary $\partial I$. Let $\bar D^m$ be the union of small tubular neighborhoods of
$D^{p+1}$, $I$ and $\bar D^{p+1}$ in $S^m$. Clearly,
the intersection $(\operatorname{Im}f_1 \cup\operatorname{Im}f_2)\cap \bar D^m$ is
standard. By the {\it sum} of $f_1$ and $f_2$ we mean the almost embedding $f_1+f_2:S^p\times S^q\to S^m$
obtained by the $S^p$-parametric connected sum of $f_1$ and $f_2$ inside $\bar D^m$.
\enddefinition

\definition{Definition of the additive inverse} Let $f:S^p\times S^q\to S^m$
be an almost embedding. By the {\it inverse} of $f$ we mean the almost
embedding given by the formula $(-f)(x,y)=\sigma_mf(x,\sigma_qy)$,
where $\sigma_k$ is the reflection of $S^k$ across the hyperplane $x_1=0$.
\enddefinition

\definition{Definition of the neutral element}
Denote by $0$ the {\it standard embedding}
$S^q\times S^p\to D^{q+1}\times D^{p+1}\cong D^{p+q+2}\subset D^m\subset S^m$.
\enddefinition

The following important result is checked directly (see the details in \cite{Sko08}).

\proclaim{Group Structure Theorem~5.4} \cite{Sko08}
A commutative group structure on the set of almost embeddings $S^p\times S^q\to S^m$ up to almost concordance
is well-defined for $m> 2p+q+2$ and $q> 1$ by the above construction.
\endproclaim

The analogues of this theorem for proper almost embeddings and smooth embeddings are also true and can be proved analogously.

\proclaim{Proposition~5.5}
For $m>2p+q+2$ the relative $\beta$-invariant $\beta:\overline{E}^m(S^p\times D^q,S^p\times S^{q-1})\to \Omega^m_{p,q}$ is a homomorphism.
\endproclaim

\demo{Proof} Consider a triple
of proper almost embeddings $f_1$, $f_2$, and $f_3=f_1+f_2$,
where ``$+$'' denotes the $S^p$-parametric connected sum
relatively to the boundary. One can see
that up to homotopy $\tilde f_3=\tilde f_1\vee\tilde f_2$ and
$\bar f_3=\bar f_1 \natural \bar f_2$, where ``$\natural$'' denotes the connected sum relatively to the boundary. So $\beta(f_3)=\beta(f_1)\sqcup\beta(f_2)$, hence
$\beta(f_1+f_2)=\beta(f_1)+\beta(f_2)$.
\qed\enddemo

Now we are going to prove that the action of embeddings $S^{p+q}\to S^m$ on the set of embeddings $S^p\times S^q\to S^m$
is injective.

\smallskip
{\bf Notations:}

\flushpar(a) $E^m(S^{p+q})$ is the group of all embeddings $S^{p+q}\to S^m$ up to concordance;

\flushpar(b) $E^m(S^p\times S^q)$ is the group of all embeddings $S^p\times S^q\to S^m$ up to concordance
(with the parametric connected sum group structure).

\flushpar(c) $\kappa^*:E^{m}(S^{p+q})\to E^{m}(S^p\times S^q)$ is a map taking an embedding $g:S^{p+q}\to S^m$ to the connected sum of $g$ and the standard embedding $S^p\times S^q\to S^m$. (The connected sum is made along an arc $I$ joining the images of the embeddings; these images are assumed to be separated by a hyperplane).

\proclaim{Proposition~5.6} The
map $\kappa^*:E^{m}(S^{p+q})\to E^{m}(S^p\times S^q)$ is injective for $m> 2p+q+2$.
\endproclaim

In fact, this proposition immediately implies the case ``$p+q+1$ divisible by $4$'' of
Theorem~1.2, by Theorem~1.1 in \S1.

\demo{Proof of Proposition 5.6}
It suffices to construct a left inverse
$\bar\kappa^*:E^{m}(S^p\times S^q)\to E^{m}(S^{p+q})$ of the map $\kappa^*$.

The map $\bar\kappa^*:E^{m}(S^p\times S^q)\to E^{m}(S^{p+q})$ is defined as follows.
Take an embedding $f:S^p\times S^q\to S^m$. By Proposition~5.3 there exists a web $D^{p+1}$
of this embedding.
Perform embedded surgery on $S^p\times S^q$ along the framed disc $D^{p+1}$.
Let $\bar\kappa^*(f)$ be the isotopy class of the embedding $S^{p+q}\to S^m$ obtained by the surgery.

The element $\bar\kappa^*(f)$ is well-defined by the second assertion of Proposition~5.3.
We have $\bar\kappa^*\kappa^*=\mathop{id}$,
because $\bar\kappa^*\kappa^*(f)=\bar\kappa^*(f\# s)=f\#\bar\kappa^*(s)=f\# 0=f$ for any $f\in E^m(S^{p+q})$.
\qed\enddemo

\proclaim{Proposition~5.7} (a) For each $m>2p+q+2$ all proper embeddings
$S^p\times D^q\to D^m$ are properly concordant.

(b) \cite{Sko08, Triviality criterion} For $m>2p+q+2$ an embedding $S^p\times S^{q-1}\to S^{m-1}$
is concordant to the standard embedding if and only if it extends to a proper embedding $S^p\times D^q\to D^m$.
\endproclaim

\demo{Proof} (a) Take a proper embedding $f:S^p\times D^q\to D^m$. By an analogue
of Proposition~5.3 there exists a web $D^{p+1}$ of $f$. Let $\bar D^m$ be the tubular neighborhood of
$D^{p+1}$. Clearly, the restriction $f:f^{-1}\bar D^m\to \bar D^m$
is concordant to the standard embedding $S^p\times D^q\to D^m$.
It remains to prove that $f$ is concordant to this restriction.
Let $c:S^p\times D^q\times I\to D^m\times I$ be the identical concordance (equal to the embedding $f:S^p\times D^q\to D^m$ for each $t\in I$).
Let $h_m:D^m\times I\to D^m\times I$ be a diffeomorphism fixed on $D^m\times 0$
and taking $\bar D^m\times 1$ to $D^m\times 1$. Let $h_{p+q}:S^p\times D^q\times I\to S^p\times D^q\times I$
be a diffeomorphism fixed on $S^p\times S^q\times 0$ and taking $S^p\times D^q\times 1$ to $f^{-1}\bar D^m$.
Then the composition
$h_m c h_{p+q}$ is a proper concordance between $f:S^p\times D^q\to D^m$ and the restriction $f:f^{-1}\bar D^m\to \bar D^m$.

(b) follows directly from (a).
\qed\enddemo

The analogue of this proposition for proper almost embeddings is also true, except that the ball $B^{p+q}$ in the definition of a proper almost embedding should be replaced by a codimension $0$ ball $\bar B^{p+q}\subset S^p\times D^q$ meeting the boundary at a common face. A {\it face} of the ball $\bar B^{p+q}$ is a ball contained in $\partial \bar B^{p+q}$.

\demo{Proof of Lemma 5.2} (cf. \cite{Sko09, Proof of Theorem 3.1}) (1) {\it
Construction of the homomorphisms.} Let $e$ be the obvious map.
Let $p$ be the ``restriction to the boundary'' map. The homomorphism $h$
is the ``cutting'' map defined as follows. Take an almost embedding
$f:S^p\times S^q\to S^m$. By Proposition~5.3 there exists a web
$D^{p+1}\subset S^m$. Let $\bar D^m$ be a tubular neighborhood of
$D^{p+1}$. Set $h(f)$ to be the restriction of $f$ to a map
$(S^p\times S^q-f^{-1}\operatorname{Int} \bar D^m)\to
S^m-\operatorname{Int} \bar D^m$.

(2) {\it Exactness at $E^m(S^p\times S^q)/E^m(S^{p+q})$.}  The sequence is exact at
$E^m(S^p\times S^q)/E^m(S^{p+q})$ because an embedding $f:S^p\times S^q\to S^m$
extends
to a proper almost embedding $S^{p}\times D^{q+1}\to D^{m+1}$ if and
only if $f$ is almost concordant to the standard embedding (by the analogue of Proposition~5.7(b)).

(3) {\it Exactness at $\overline{E}^m(S^p\times D^q,S^p\times S^{q-1})$.} The inclusion $\operatorname{Im}h\subset\operatorname{Ker}p$
follows by Proposition~5.7(b).
To prove $\operatorname{Ker}p\subset\operatorname{Im}h$
take a proper almost embedding $f:S^p\times D^q\to D^m$ such that $p(f)=0$. Let $f\left|_\partial\right.:S^p\times \partial D^q\to \partial D^m$
be the restriction of $f$ to the boundary.
By definition, there exists a smooth embedding $g:S^{p+q-1}\to S^{m-1}$  such that the connected sum $f\left|_\partial\right.+g$ is concordant to the standard embedding. Extend $g:S^{p+q-1}\to S^{m-1}$ to a proper almost embedding $g':D^{p+q}\to D^m$.
Let $f+g'$ be the connected sum of $f$ and $g'$ relatively to the boundary. By Proposition~5.7(b)
the map $f+g':S^p\times D^q\to D^m$ extends to an almost embedding $f':S^p\times S^q\to S^m$. Thus $f=h(f')$.

(4) {\it Exactness at $\overline{E}^m(S^p\times S^q)$.} The inclusion
$\operatorname{Im}e\subset \operatorname{Ker}h$ follows by
Proposition~5.7(a). To prove $\operatorname{Ker}h\subset\operatorname{Im}e$, take an almost embedding $f:S^p\times S^q\to S^m$
such that $h(f)=0$. Then, by definition,
there exist a proper almost concordance $c$ between $h(f)$ and a connected sum of
the standard embedding $S^p\times D^q\to D^m$ and a proper almost embedding $g:D^{p+q}\to D^m$.
Since the restriction of $c$ to the boundary is a smooth concordance, by Proposition~5.6
it follows that the restriction of $g$ to the boundary is unknotted. Thus we may assume that $g$ is a smooth
embedding. By the concordance extension theorem 
the restriction of $c$ to the boundary extends to an ambient
concordance of the disc $S^m-D^m$. So $c$ can be extended to an almost concordance
of $f$ without adding new self-intersections. The latter is an
almost concordance between $f$ and an embedding $f':S^p\times S^q\to S^m$.
Hence $f=e(f')$. \qed\enddemo

So the proof of Theorem~2.1 is completed. \qed

\head 6. The finiteness criterion \endhead

In this section we deduce Theorem~1.2 from Theorem~2.1. Thus we need a
classification of almost embeddings $S^p\times S^q\to S^m$. We summarize this classification and the above results
in the following theorem.

\smallskip
{\bf Notation:}
$E^m(D^p\times S^q)$ denotes the group of smooth
embeddings $D^p\times S^q\to S^m$ up to concordance (with the $D^p$-parametric connected sum
group structure).

\proclaim{Theorem~6.1}
For $m>p+\frac{4}{3}q+2$ and $m > 2p+q+2$  there exist the
following exact sequences

\smallskip
\flushpar $ \matrix
\format\l&\quad\l&\,\,\c&\,\,\l&\,\,\c&\,\,\l&\,\,\c&\,\,\l&\,\,\c&\,\,\l\\
(1) &0 &\to &E^m(S^{p+q}) &\to &{E}^m(S^p\times S^q)&\to &\frac{{E}^m(S^p\times S^q)}{E^{m}(S^{p+q})}
\\
(2) &\overline{E}^{m+1}(S^p\times S^{q+1}) &\to &\Omega^{m+1}_{p,q+1} &\to &\frac{{E}^m(S^p\times S^q)}{E^{m}(S^{p+q})}&\to
&\overline{E}^m(S^p\times S^q)
\\
(3) &E^{m+1}(D^p\times S^{q+1})&\to&\pi_{p+q}(S^{m-q-1}) &\to &\overline{E}^m(S^p\times S^q) &\to
&E^m(D^p\times S^q)
\\
(4) &E^{m+1}(S^{q+1})&\to&\pi_{q}(V_{m-q,p})&\to &E^m(D^p\times S^q)&\to &E^m(S^{q}).
\endmatrix
$
\endproclaim


The groups in the second column of~6.1 are
well-known {\it rationally}:

\proclaim{Theorem~6.2} Assume that $p+\frac{4}{3}q+2\le m <
p+\frac{3}{2}q+2$, $m > 2p+q+2$ and $m>n+2$. Then

\smallskip
\flushpar$(1)$ 
$E^{m}(S^n)$ is infinite if and only if $m\le\frac{3}{2}n+\frac{3}{2}$,
$n+1$ is divisible by $4$.

\flushpar$(2)$ $\Omega_{p,q+1}^{m+1}$ is infinite if and only if
$m=p+\frac{3}{2}q+\frac{3}{2}$, and $q+1$ is divisible by $4$.

\flushpar$(3)$ $\pi_{p+q}(S^{m-q-1})$ is infinite if and only if
$m=\frac{1}{2}p+\frac{3}{2}q+\frac{3}{2}$, and $p+q+1$ is divisible by~$4$.

\flushpar$(4)$ $\pi_q(V_{m-q,p})$ is infinite if and only if
$p\ge1$, $\frac{3}{2}q+\frac{3}{2}\le m\le p+\frac{3}{2}q+\frac{1}{2}$,
and $q+1$ is divisible by $4$.
\endproclaim

We argue by the following plan.
First we prove Theorem~1.2 modulo Theorem~6.1 and~6.2. Then we prove
Theorems~6.1 and~6.2 themselves, using some known results.

\demo{Proof of Theorem~1.2 modulo 6.1 and 6.2} (1) {\it Case when $q+1, p+q+1$ are not divisible by $4$.}
Recall that if $X\to Y\to Z$ is an exact
sequence with finite $X$ and $Z$, then $Y$ is also finite.
Applying this 4 times to the last $3$ columns of Theorem~6.1
starting from the bottom, we are done, because by Theorem~6.2 the
groups in the second column of~6.1 are finite when $q+1, p+q+1$ are not divisible by $4$.

(2) {\it Case when $p+q+1$ is divisible by $4$.} By Theorems~6.1(1) and~6.2(1) it follows directly that
the group $E^m(S^p\times S^q)$ is in this case infinite.

\smallskip

(3) {\it Case when $q+1$ is divisible by $4$, $m\le
\frac{3}{2}q+\frac{3}{2}$.}
By Theorem~6.2(1)
the group $E^m(S^q)$ in this case is infinite. Take an infinite
order element $x$. The obstruction to existence of a $(p+1)$-framing
on the embedding $x:S^{q}\to S^m$ belongs to the group $\pi_{q-1}(V_{m-q,p+1})$. By Theorem~6.2(4)
this group in our case is finite. So
for some positive integer $N$ the embedding $Nx$ extends to a
smooth embedding $H:S^p\times S^q\to S^m$. Since the restriction of the embedding to the sphere $*\times S^q$ has an infinite order it follows that $H$ itself has an
infinite order.

\smallskip
(4) {\it Case when $q+1$ is divisible by $4$, $\frac{3}{2}q+\frac{3}{2}< m\le p+\frac{3}{2}q+\frac{1}{2}$, $p\ge 1$.} It suffices to construct an embedding $T:S^p\times S^q\to S^m$ having an infinite order in
$E^m(S^p\times S^q)$.

{\it Construction of the embedding $T$.} By Theorem~6.2(4) the group
$\pi_q(V_{m-q,p})$ is in this case infinite. Take an infinite
order element $x$ of this group. Consider the map
$\tau:\pi_q(V_{m-q,p})\to E^m(D^p\times S^q)$ from~6.1(4). This map
takes the element $x$ to the canonical $p$-frame $D^p\times S^q\to
S^m$ of the standard sphere $S^q\subset S^m$. The complete
obstruction to extension of this $p$-frame to a $(p+1)$-frame
belongs to $\pi_{q-1}(S^{m-p-q-1})$. The latter group is in
our case finite. So for some positive integer $N$ the element $N\tau(x)$ can be
extended to a smooth embedding $S^p\times S^q\to S^m$, which is
the desired torus $T$.

{\it Proof that $T$ has an infinite order.} It suffices to prove
that the element $\tau(x)\in E^m(D^p\times S^q)$, which is the
restriction of $T$ to $D^p\times S^q$, has an infinite order.
Suppose to the contrary. Then $N\tau(x)=0$ for some positive integer $N$. So
by Theorem~6.1(4) $Nx$ belongs to the image of the map
$E^{m+1}(S^{q+1})\to\pi_{q}(V_{m-q,p})$. But the group $E^{m+1}(S^{q+1})$ is in our case finite. Thus $x$ has finite order in contrast to our choice above. This contradiction proves that $T$ has infinite order.

\smallskip
(5) {\it Case when $q+1$ is divisible by $4$, $m=p+\frac{3}{2}q+\frac{3}{2}$.}

{\it Construction of the embedding $W$.} By Theorem~6.2(2)
the group $\Omega^{m+1}_{p,q+1}$ is in this case infinite. Take an
infinite order element $x$ of this group. Let $W:S^p\times S^q\to S^m$
be an embedding realizing the image of $x$ under the map
$\Omega^{m+1}_{p,q+1}\to{{E}^m(S^p\times S^q)}/{E^{m}(S^{p+q})}$
from Theorem~6.1(2).

{\it Proof that $W$ has an infinite order.} Consider the exact
sequence~6.1(2). It
suffices to prove that $\overline{E}^{m+1}(S^p\times S^{q+1})$ is in
our case finite. Since $q+1$ is divisible $4$, it follows by Theorem~6.2 that
$\pi_{q+1}(V_{m-q,p})$ and $E^{m+1}(S^{q+1})$ are finite. An easy computation shows
that $\pi_{p+q+1}(S^{m-q-1})$ is also finite in our case.
So by Theorems~6.1(3)--(4) it follows that $\overline{E}^{m+1}(S^p\times S^{q+1})$ is finite.
\qed\enddemo

In the rest of the section we shall prove Theorems~6.1 and~6.2.

Assertions~(1) and~(2) of Theorem~6.1 are reformulations of Proposition~5.6 and~Theorem~2.1 which were proved in~\S5.
Theorem~6.1~(4) is proved by a direct verification analogously to
\cite{Hae66A, Corollary~5.9}.
We sketch the proof below for the convenience of the reader.
Theorem~6.1(3) is proved in \cite{Sko08, Restriction Lemma~5.2}
for $p\ge 1$, $m>2p+q+2$ and $m\ge \frac{1}{2}p+\frac{3}{2}q+2$. The additional restriction
$m\ge \frac{1}{2}p+\frac{3}{2}q+2$ in the proof of this assertion in \cite{Sko08}
can be dropped. We sketch an alternative proof to keep the paper self-contained.

\demo{Sketch of the proof of assertion (3) in Theorem~6.1} (a) {\it Definition of the
groups $\tilde E^m(S^p\times S^q)$ and $\tilde{E}^m(S^p\times D^q,S^{p}\times S^{q-1})$.} A map
$f:S^p\times S^q\to S^m$ is said to be a {\it weak almost
embedding}, if it is an embedding outside the fixed ball $B^{p+q}\subset
S^p\times S^q$ (the intersection $fB^{p+q}\cap f(S^p\times S^q-B^{p+q})$ is allowed to be nonempty). A {\it weak almost concordance} is defined
analogously. Denote by $\tilde E^m(S^p\times S^q)$ the group of weak almost
embeddings up to weak almost concordance.
Identify the groups $\tilde E^m(S^p\times S^q)$ and $E^m(D^p\times S^q)$.  Clearly, these groups are isomorphic for $m>2p+q+2$.

Fix a ball $\bar B^{p+q}\subset (S^p-*)\times D^q$ meeting the boundary at its face.
A proper map $f:S^p\times D^q\to D^m$ is said to
be a {\it proper weak almost embedding}, if the following two
conditions hold:

\flushpar(i) $f$ is an embedding outside $\bar B^{p+q}$; and

\flushpar(ii) $f(S^p\times\partial D^q\cap \bar B^{p+q})\cap
f(S^p\times\partial D^q - \bar B^{p+q})=\emptyset$.

\flushpar A {\it proper weak almost concordance} is defined
analogously. Denote by $\tilde{E}^m(S^p\times D^q,S^{p}\times S^{q-1})$
the group of proper weak almost embeddings up to proper weak almost concordance.

(b) {\it For every $m>2p+q+2$ there exist an exact sequence:}
$$
\overline{E}^m(S^{p}\times S^{q}) \overset{e}\to{\to}
\tilde{E}^m(S^p\times S^q)\overset{h}\to{\to}
\tilde{E}^m(S^p\times D^q,S^{p}\times S^{q-1})\overset{p}\to{\to}
\overline{E}^{m-1}(S^{p}\times S^{q-1})\to\dots
$$
Here $e$, $h$ and $p$ are the obvious forgetful, cutting and
restriction homomorphisms, respectively. This assertion is proved
completely analogously as Lemma~5.2.

(c) {\it Definition of the homomorphism
$\lambda:\tilde{E}^m(S^p\times D^q,S^{p}\times S^{q-1})\to\pi_{p+q-1}(S^{m-q-1})$.} Take a proper
weak almost embedding $f:S^p\times D^q\to D^m$. By definition
$f\partial \bar B^{p+q}\cap f(*\times D^q)=\emptyset$. Notice that
$D^m-f(*\times D^q)\simeq S^{m-q-1}$. Let $\lambda(f)$ be the
homotopy class of the restriction $f:\partial \bar B^{p+q}\to D^m-f(*\times D^q)$.

(d) {\it $\lambda$ is injective.} Take a proper weak almost
embedding $f:S^p\times D^q\to D^m$ such that $\lambda(f)=0$. Then
$f\left|_{\partial \bar B^{p+q}}\right.$ extends to a map $g:\bar B^{p+q}\to D^m$
missing $f(*\times D^q)$. Since $f\left|_{S^p\times D^q-\bar B^{p+q}}\right.$
is an embedding it follows that
the intersection of $g\bar B^{p+q}$ with $f(S^p\times D^q-\bar B^{p+q})$ can be removed
by a homotopy relatively to the boundary.
Thus we may assume that $g$ misses
$f(S^p\times D^q-\bar B^{p+q})$. Perform a proper weak almost concordance which
replaces $f\left|_{\bar B^{p+q}}\right.$ by $g$. By an analog of Proposition~5.7(a)
the obtained map is properly weakly almost concordant to the standard
embedding $S^p\times D^q\to D^m$.

(e) {\it $\lambda$ is surjective.} Take an element $x\in
\pi_{p+q-1}(S^{m-q-1})$. Take the standard embedding
$f:S^{p}\times D^q\to D^m$. Realize the element $x$ by a map
$g:S^{p+q-1}\to \partial D^m-f(*\times D^q)$.
Since $f\left|_{S^p\times D^q-\bar B^{p+q}}\right.$
is an embedding it follows that
the intersection of $g\bar B^{p+q}$ with $f(S^p\times D^q-\bar B^{p+q})$ can be removed
by a homotopy relatively to the boundary.
Thus we may assume that $g$ misses
$f(S^p\times D^q-\bar B^{p+q})$.
Extend the map $g:S^{p+q-1}\to \partial D^m$ to a
proper map $g':D^{p+q}\to D^m$. Let $\mu_x$ be the connected sum
(relatively the boundary) of $g'$ and $f$.
Clearly, $\lambda(\mu_x)=x$. This completes the proof of
assertion~(3). \qed\enddemo

\demo{Sketch of the proof of assertion (4) in Theorem~6.1}
(a) {\it Definition of homomorphisms.} The map $i^*:{E}^m(D^p\times S^q) \to E^{m}(S^q)$
is restriction-induced. Here $0\times S^q$ is identified with $S^q$ in obvious way.

The map $Ob:E^{m}(S^q)\to \pi_{q-1}(V_{m-q,p})$ is the complete obstruction to
the existence of a $p$-framing on an embedding $S^q\to S^m$.
This obstruction is defined as follows. Take an embedding $f:S^q\to S^m$.
Take a (unique up to homotopy)
$(m-q)$-framing
of the disc $fD^q_+$. Take a (unique up to homotopy) $p$-framing of the disc $fD^q_-$. Thus the sphere $fS^{q-1}$ is equipped both with the $p$-framing and the $(m-q)$-framing. Using the $(m-q)$-framing identify each fiber of the normal bundle to $fD^q_+$ with the space $\Bbb{R}^{m-q}$. To each point $x\in S^{q-1}$ assign the $p$-framing at the point $fx$. This leads to a map $S^{q-1}\to V_{m-q,p}$. By definition $Ob(f)\in \pi_{q-1}(V_{m-q,p})$ is the homotopy class of this map.

The map $\tau:\pi_{q}(V_{m-q,p}) \to {E}^m(D^p\times S^q)$ is defined as follows.
Represent $f\in \pi_{q}(V_{m-q,p})$ as a smooth map $f:D^p\times S^q\to D^{m-q}$ linear in each fiber $D^p\times y$, $y\in S^q$.
Define $\tau(f)$ to be the composition $D^p\times S^q\to D^{m-q}\times S^q\to S^m$ of the embedding $f\times \operatorname{pr}_2$ and the standard
embedding $s$, i.e., $\tau(f)(x,y)=s(f(x,y),y)$ for each $x\in S^p$, $y\in S^q$.

(b){\it Exactness at ${E}^m(D^p\times S^q)$ and $E^{m}(S^q)$} is checked directly.

(c){\it Proof of the exactness at $\pi_{q}(V_{m-q,p})$.} Let $f:S^{q+1}\to S^{m+1}$ be an embedding. Then $f$ is isotopic to a {\it standardized} embedding $f':S^{q+1}\to S^m$, i.e., satisfying the conditions:

\noindent $\bullet$ $f':D^{q+1}_-\to D^{m+1}_-$ is the restriction of the inclusion $S^{q+1}\to S^{m+1}$;

\noindent $\bullet$ $f'(\mathop{Int}D^{q+1}_+)\subset \mathop{Int}D^{m+1}_+$.

Take a $p$-framing of $f'(D^{q+1}_+)$.
Represent this framing as an embedding $g:D^p\times D^{q+1}_+\to D^{m+1}_+$ linear in each fiber $D^p\times y$, $y\in D^{q+1}_+$. By definition $\tau Ob (f')=g\left|_{D^p\times \partial D^{q+1}_+}\right.$.
Thus $\tau Ob(f'):D^p\times \partial D^{q+1}_+\to \partial D^{m+1}_+$ extends to the proper embedding $g:D^p\times D^{q+1}_+\to D^{m+1}_+$.
So $\tau Ob (f')$ is isotopic to the standard embedding $D^p\times S^q\to S^m$. Thus $\mathop{Im} Ob\subset \mathop{Ker} \tau$. Analogously $\mathop{Im} Ob\supset \mathop{Ker}\tau$.
\enddemo

Theorem~6.2 can be easily reduced to known results. Assertion (1) is proved in
\cite{Hae66A, Corollary~6.7}. Assertion (3) follows from the Serre theorem.

\demo{Proof of assertion (4) in Theorem~6.2}
The assumptions $m> 2p+q+2$ and $m<p+\frac{3}{2}q+2$
together imply that $m\le 2q$. We are going to prove assertion~(4)
with assumptions $m> 2p+q+2$ and $m<p+\frac{3}{2}q+2$ replaced by
the only assumption $m\le 2q$. We use induction over $p$.

(a) {\it Case $q+1$ not divisible by 4.} Since $m\le 2q$, it follows
that $\pi_q(V_{m-q,1})\cong \pi_q(S^{m-q-1})$ is finite.
Using the homotopy exact sequence of
the ``restriction'' bundle $S^{m-p-q}\to V_{m-q,p}\to
V_{m-q,p-1} $ tensored by $\Bbb Q$, we get inductively
that $\pi_q(V_{m-q,p})$ is finite.

(b) {\it Case $q+1$ divisible by 4, and either $m<\frac{3}{2}q+\frac{3}{2}$ or
$m>p+\frac{3}{2}q+\frac{1}{2}$.} In this case the groups $\pi_q(S^{m-q-i})$
are still finite for each $i=1,2,\dots,p$. Similarly as above we get that $\pi_q(V_{m-q,p})$ is finite.

(c) {\it Case $q+1$ divisible by 4, and
$\frac{3}{2}q+\frac{3}{2}\le m\le p+\frac{3}{2}q+\frac{1}{2}$.}
Take $i$ such that $m=i+\frac{3}{2}q+\frac{1}{2}$. Consider
the exact homotopy sequence above for $p=i$. Analogously as above
it can be shown that for $q+1$ divisible by $4$ and $m\le 2q$
the group $\pi_{q+1}(V_{m-q,i-1})$ is finite. Thus the group
$\pi_q(V_{m-q,i})$ is infinite. By induction $\pi_q(V_{m-q,p})$
is also infinite.
\qed\enddemo

\demo{Proof of the assertion (2) in Theorem~6.2}
Use the notation $s=2p+3q-2m+2$ and $l=m-p-q-1$. Then the group in question is $\pi_{s+l+1}(V_{M+l,M})$. Our restriction $p+\frac{4}{3}q+2\le m < p+\frac{3}{2}q+2$ is equivalent to the restriction $-1\le s\le l-3$.

(a) {\it Case $s=-1$}. By tables in~\cite{Pae56} the group $\pi_{l}(V_{M+l,M})$
is infinite if and only if $l$ is divisible by $2$. Together
with condition $s=-1$ this is equivalent to the conditions
$m=p+\frac{3}{2}q+\frac{3}{2}$, $q+1$ is divisible by 4.

(b) {\it Case $0\le s\le l-3$}.
Let us prove by induction over $s$ that the
group $\pi_{s+l+1}(V_{M+l,M})$ is finite. The base $s=0$ follows from tables in~\cite{Pae56}. For $s>0$
consider the homotopy exact sequence of the ``restriction'' bundle
$ S^{l}\to V_{M+l,M}\to V_{M+l,M-1} $ tensored by
$\Bbb Q$. In this sequence $\pi_{s+l+1}(S^l)$ is finite because $0\le s\le l-3$. By the inductive hypothesis the group
$\pi_{s+l+1}(V_{M+l,M-1})\cong\pi_{(s-1)+(l+1)+1}(V_{M+l+1,M})$ is finite. Hence the group $\pi_{s+l+1}(V_{M+l,M})$ is finite.
\qed\enddemo

\proclaim{Remark~6.3} Analogously using the bundle $SO_{m-q-p-1}\to SO_{m-q}\to \pi_{q}(V_{m-q,p+1})$ one can prove the following:
{\it suppose that $m\ge p+\frac32 q +2$; then the group $\pi_q(V_{m-q,p+1})$
is infinite if and only if either $m=2q+1$, $q$ odd, or $m=p+2q+1$, $q$ even} \cite{CFS11}.
\endproclaim

Thus the proof of Theorem~1.2 is completed. \qed

\head 7. Concluding remarks \endhead

Let us give a counterexample which shows that Theorem~2.1 in \S2 and Lemma~4.2 in \S4 cannot be proved using the ordinary Whitney trick
or the approach of \cite{HaQu74}:

\proclaim{Example} For some $p$, $q$, $m$ satisfying the inequalities of Theorem~2.1
 there exists a proper almost embedding $f:S^p\times D^q\to D^m$ such that $\beta(f)=0$ but $f$ admits no
webs. In particular, the isomorphism of Lemma~4.2 does not hold for $f'=f$.
\endproclaim

\demo{Proof} Take $m=p+\frac{3}{2}q+\frac{3}{2}$ and choose $p,q$ so that
$l=m-p-q-1$ is odd. Take a generic proper almost embedding $f:S^p\times
D^q\to D^m$ such that $\Delta$ is connected whereas $\tilde\Delta$ is
not connected (for example, start with the standard embedding and
perform the finger Whitney moves). Our proof of~Lemma~4.2 in \S4
shows that, in fact $\pi_q(D^m-\operatorname{Im}F,\partial)\cong
\Omega_s(\Delta;l\lambda_\Delta)$. The latter group is isomorphic to
$\Bbb Z$ for $s=0$ and $l$ odd by \cite{HaKa98, end of \S4}. On
the other hand, $\Omega_0(P^\infty;l\lambda)=\Bbb Z_2$. So the map
$\beta(f,{\cdot}):\Bbb Z\to\Bbb Z_2$ is not injective.

Then there exists a proper map $g:D^q\to D^m-\operatorname{Im}f$
such that $\beta(f,g)=0$, but $g$ is not null-homotopic. Perform
the construction from the proof of~Theorem~5.1 modulo Lemma~4.2,
steps (2) and (3). We get a new proper almost embedding $f$ such
that $\beta(f)=0$. On the other hand, the map $g:D^q\to
D^m-\operatorname{Im}f$ is close to $f\left|_{*\times D^q}\right.$
but it is not null-homotopic. Thus there exists no web
for~ $f$. \qed\enddemo

{\bf Further investigation.} There are several directions to study
knotted tori:

\flushpar (i) {\it Explicit classification results.} How many
embeddings $S^1\times S^5\to S^{10}$ are there up to isotopy?

\flushpar (ii) {\it Weakening the dimension restrictions.} Is it possible to
drop the restrictions $m>2p+q+2$ or $m>p+\frac{4}{3}q+2$ in
Theorem~1.2? \cite{cf. KeMi63, Hae66A}

\flushpar (iii) {\it Arbitrary manifolds.} It would be interesting to
generalize the $\beta$-invariant and Theorem~2.1 to embeddings of
arbitrary manifolds \cite{ReSk99, Sko08}.

\flushpar (iv) {\it Rational classification of embeddings.} For a given
manifold $N$ and a number $m$ determine whether the set of embeddings $N\to S^m$
up to isotopy is finite.

\head Appendix. Surgery on the double point manifold \endhead

Here we perform a surgery on the double point manifold $\Delta$ to make
the classifying map $\Delta\to P^{\infty}$ of the covering $\tilde\Delta\to\Delta$ sufficiently highly connected. This is required for the proof of Complement Lemma~4.2 stated in \S4.
Our exposition is completely analogous to \cite{HaKa98, Appendix~A} although more general, explicit and detailed.

Let $f:D^n\to M^m$ be a general position proper immersion such
that $f\left|_{\partial D^n}\right.:\partial D^n\to \partial M^m$ is an embedding. The
embedding theorem \cite{Hud69}, cf. Theorem~4.3 above, allows us to
remove the self-intersection of $f$ by a homotopy
$\operatorname{rel}
\partial D^n$ under certain conditions. In the dimension range, where
the embedding theorem is not true, we give an approach to ''simplify'' the double points of $f$.

In this section we prove that the classifying map
$\Delta\to P^\infty$ can be made $(s+1)$-connected by a homotopy $\operatorname{rel} \partial D^n$
of the map $f:D^n\to M^m$, provided that

\flushpar (i) $M^m$ is $(s+1)$-connected;

\flushpar (ii) $2s\le 2n-m-2$; and

\flushpar (iii) $0\le s\le m-n-3$.

This is a restatement of Theorem~4.7 above.

\demo{Proof of Theorem~4.7} The map $\Delta\to P^\infty$ is made
$(s+1)$-connected in 2 steps:

\flushpar Step 1. Making $\Delta$ connected and $\pi_1(\Delta)\to\pi_1(P^\infty)$
surjective (i.e., $\tilde \Delta$ connected).

\flushpar Step 2. Killing the elements of
$\operatorname{Ker}(\pi_i(\Delta)\to\pi_i(P^\infty))$ for each $i=1,\dots, s$.

These two steps are sufficient for the proof of theorem, because the map
$\pi_{i+1}(\Delta)\to\pi_{i+1}(P^\infty)=0$ is surjective for
$1\le i\le s$.

In both steps 1 and 2 we make the following Whitney-Haefliger trick,
performing a surgery on $\Delta$.
First let us construct {\it the Habegger-Kaiser standard model}
for doing surgery on a framed $i$-sphere of double points of an
$n$-disc immersed into $S^m$.

{\bf Standard model for doing surgery.} \cite{HaKa98} We will make
use of the ``model'' manifold $\Bbb R^m=\Bbb R\times \Bbb
R^{i+1}\times \Bbb R^{2n-m-i} \times \Bbb R^{m-n-1}\times \Bbb
R^{m-n-1}$ and of two embeddings, $g_+$ and $g_-$ of $\Bbb
R^n=\Bbb R^{i+1}\times \Bbb R^{2n-m-i}\times \Bbb R^{n-m-1}$ into
$\Bbb R^m$ intersecting transversally along $0\times S^i\times
\Bbb R^{2n-m-i}\times 0\times 0$. For example, one may take
$g_-(x,y,z)=(|x|^2-1,x,y,0,z)$ and $g_+(x,y,z)=(1-|x|^2,x,y,z,0)$.
The sphere $S^i$ bounds a ball $D^{i+1}\subset\Bbb
R^{i+1}\subset\Bbb R^n$. Furthermore, the sphere
$S^{i+1}=D^{i+1}_+\cup D^{i+1}_-$, where $D^{i+1}_\pm=g_\pm(
D^{i+1}\times 0\times 0)$, bounds a ball $D^{i+2}\subset\Bbb R\times\Bbb
R^{i+1}\subset\Bbb R^m$ with corners along $S^i$. Pushing one of
the two caps $D ^ {i+2}_\pm $ across $D^{i+2}$ does the required surgery. More
precisely, the double points of the resulting regular homotopy form
the trace of this surgery.

Now we are going to make some preparations for doing surgery,
which are a bit different, in Steps 1 and 2.

{\bf Step 1: making $\Delta$ and $\tilde \Delta$ connected.}
If $\Delta=\emptyset$ then we first create
a nonempty self-intersection (for example, by the Whitney finger moves).
Assume $\Delta\ne\emptyset$. Take a pair of
points $(a,b), (c,d)$ belonging to distinct components of $\tilde
\Delta$. One can assume that they are outside the triple point
set. Consider the spheres $S^0=\{ \{a,b\},\{c,d\}\}$, $S^0_+=\{
(a,b),(c,d)\}$, and $S^0_-=\{ (b,a),(d,c)\}$, and let $\eta$ be
the trivial normal bundle $N(S^0,\Delta)$. The surgery on $S^0$
(see Completion of the proof below) will connect distinct
components of $\tilde \Delta$, because
$\operatorname{dim}\Delta=2n-m\ge 2s+2\ge 2$.

{\bf Step 2: killing the elements of
$\operatorname{Ker}(\pi_i(\Delta)\to\pi_i(P^\infty))$.} Assume
that $g:S^i\to \Delta$ represents an element of
the kernel $\operatorname{Ker}(\pi_i(\Delta)\to\pi_i(P^\infty)$, $1\le i\le
s$. Since $2i\le\operatorname{dim}\Delta-1$ (because $2s\le
2n-m-2$), it follows that $g$ can be assumed to be an embedding.
By general position the triple point set has dimension $\le 3n-2m$. Since $s\le m-n-3$, it follows that
$i+3n-2m\le\operatorname{dim}\Delta-1$, so generically
$\operatorname{Im}g$ does not contain triple points.

Since the composition $S^i\to \Delta\to P^\infty$ is trivial, it
follows that $(m-n)\lambda_\Delta\left|_{S^i}\right.$ is trivial.
Then the normal bundle $\eta=N(S^i,\Delta)$ is stably trivial and hence trivial. Since $S^i\to P^{\infty}$ is trivial, it
follows that $S^i$ is also trivially covered in $\tilde \Delta$.
Denote by $S^i_+$ and $S^i_-$ the two copies of $S^i$ in $\tilde
\Delta$.

{\bf Completion of the proof: surgery on $S^i$.}
First let us span the spheres $S^i_+$ and $S^i_-$ by two disjoint
discs $D^{i+1}_+$ and $D^{i+1}_-$ in the ball $D^n$. To do this take a trivialization of the normal bundle
$N(\tilde \Delta, D^n)$ (e.g., constructed in the definition of the double point $(m-n)\lambda$-manifold in \S3). Push the spheres $S^i_+$
and $S^i_-$ along the first vector field of the trivialization.
Since $s\le m-n-3$ and $2s\le 2n-m-1$ it follows that
$3s\le n-4$. Hence $i+1+\operatorname{dim}\Delta \le
n-1$ and $2(i+1)\le n-1$. Then by general position the pushed spheres can be spanned by discs missing the image of $\tilde\Delta$.

Consider the obvious decomposition
$N(D^{i+1}_+,D^n)\left|_{S^i_+}\right.=\eta\oplus\epsilon^{m-n-1}$,
where $\eta=N(S^i,\Delta)$.
This decomposition gives an $(m-n-1)$-framing of the sphere  $S^i_+$. We wish to extend this $(m-n-1)$-framing over $D^{i+1}_+$. The complete obstruction lies in
$\pi_i(V_{n-i-1,m-n-1})$. The latter group vanishes for $2i\le 2n-m-1$, which follows from $2s\le 2n-m-2$. Thus we obtain a decomposition
$N(D^{i+1}_+,D^n)=\eta_+\oplus\epsilon^{m-n-1}$, where $\eta_+$ is the extension of the bundle $\eta$ over $D^{i+1}_+$ complimentary to the $(m-n-1)$-framing. Define the bundle $\eta_-$ analogously.

Next we extend the embedding of $S^{i+1}=D^{i+1}_+\cup D^{i+1}_-$
to an embedding of $D^{i+2}$ into $S^m$. To do this push $D^{i+1}_+$
and $D^{i+1}_-$ along the first vector field of a trivialization of the
bundle $N(D^n,S^m)$. In this way we obtain an embedding of a collar
neighborhood of the sphere $S^{i+1}$ into $S^m$.
Since $n+s+3\le m$ and $s\le n-1$ it follows that
$i+2+n\le m-1$ and $2(i+2)\le m-1$.
Thus by general position the collar can be extended
an embedded disc $D^{i+2}$ in $S^m$, whose interior
does not intersect $fD^n$.

Finally, consider the following $(m-n-1)$-framing
of the sphere $S^{i+1}$.
On the disc $D^{i+1}_+$ take the $(m-n-1)$-framing complementary
to $\eta_+$. On the disc $D^{i+1}_-$ take
the $(m-n-1)$-framing obtained from the trivialization
of $N(D^n,S^m)$ by forgetting the first vector field.
By the construction above it follows that these two
framings coincide on $S^i$. Thus we obtain an $(m-n-1)$-framing of $S^{i+1}$.
Let us extend it to $D^{i+2}$.
The complete obstruction lies in $\pi_{i+1}(V_{m-i-2,m-n-1})$.
The latter group vanishes for $2i+2\le n-1$ because
$3s\le n-4$.
Let $\eta'$ be the complementary bundle to
the obtained $(m-n-1)$-framing over $D^{i+2}$.
We have a splitting $\eta'=\eta_-\oplus\epsilon^{m-n-1}$ on $D^{i+1}_-$.
Extending it to $D^{i+2}$ we get
a splitting
$N(D^{i+2},S^m)=\eta''\oplus\epsilon^{m-n-1}\oplus\epsilon^{m-n-1}$
for some bundle $\eta''$.

Thus the relevant framing information along $D^{i+2}$ agrees with that of the standard model. So there is a diffeomorphism between a neighborhood of $D^{i+2}$ and
the Euclidean $m$-space taking a restriction of our immersion
to the standard model. So one can
perform the surgery and kill the spheroid $g:S^i\to \Delta$. The theorem is thus proved. \qed
\enddemo

\smallskip
{\bf Acknowledgements.}

The authors are grateful to A. Skopenkov for continuous attention
to our work and also to P.~ Akhmetiev, U.~ Kaiser, U.~ Koschorke, G.~ Laures, S.~ Melikhov, A.~Mischenko, V.~
Nezhinsky and E.~Shchepin  for useful discussions.


\Refs \widestnumber\key{RSS95A}

\ref \key Br68 \by W. Browder
\paper Embedding smooth manifolds
\yr 1968 \jour Proc. Int. Congr. Math. Moscow 1966
\pages 712--719 \endref

\ref \key CeRe03 \by M. Cencelj and D. Repov\v s \paper On
embeddings of tori in Euclidean spaces \jour Acta Math. Sinica
(Engl. Ser.) \vol 21:2 \yr 2005 \pages 435--438 \endref

\ref \key CRS04 \by M. Cencelj, D. Repov\v s and A. Skopenkov \paper
On the Browder-Levine-Novikov embedding theorems \lang in Russian
\jour Trudy Math. Inst. Ross. Akad. Nauk \vol 247 \yr 2004 \pages
280--290 \transl English transl. \jour Proc. Steklov Inst.
Math. \vol 247 \yr 2004 \pages 259--268 \endref

\ref \key CRS07 \by M. Cencelj, D. Repov\v s and M. Skopenkov \paper
Homotopy type of the complement to an immersion and the classification of embeddings of tori
\lang in Russian \jour Uspekhi Mat. Nauk \vol 62:5 \yr 2007 \pages 165--166
\transl English transl. \jour Russian Math. Surveys \vol 62:5 \yr 2007 \pages 985--987
\finalinfo arXiv:0803.4285v1 [math.GT] \endref

\ref \key CRS08 \by M. Cencelj, D. Repov\v s and M. Skopenkov \paper A new invariant of higher-dimensional embeddings
\paperinfo preprint \yr 2008 \finalinfo arXiv:0811.2745v1 [math.GT]
\endref


\ref \key CFS11 \by D. Crowley, S. Ferry and M. Skopenkov \paper The rational classification of links in codimension
$>2$ \jour Forum Math. \finalinfo to appear, arXiv:1106.1455v1[math.AT]
\endref

\ref \key CrSk08 \by D. Crowley and A. Skopenkov \paper A classification of smooth embeddings of 4-manifolds
in 7-space, II \jour
Int. J. Math. \vol 22:6 \yr 2011 \pages 731--757 \finalinfo arXiv:0808.1795v1 [math.GT]
\endref

\ref \key GW99 \by T. Goodwillie and M. Weiss \paper Embeddings from the point of view of immersion theory, II
\jour Geometry and Topology \vol 3 \yr 1999 \pages 103–-118
\endref

\ref \key HaKa98 \by N. Habegger and U. Kaiser \paper Link
homotopy in the 2-metastable range \jour Topology \vol 37:1 \yr
1998 \pages 75--94
\endref

\ref \key Hae62T \by A. Haefliger \paper Differentiable links
\jour Topology \vol 1 \yr 1962 \pages 241--244 \endref

\ref \key Hae66A \by A. Haefliger \paper Differentiable embeddings
of $S^n$ in $S^{n+q}$ for $q>2$ \jour Ann. Math. (2) \vol 83:3 \yr
1966 \pages 402--436 \endref

\ref \key Hae66C \by A.~Haefliger \paper Enlacements de spheres en
codimension superieure a 2 \jour Comment. Math. Helv. \vol 41 \yr
1966-67 \pages 51--72 \lang in French
\endref

\ref \key Ha67 \by A. Haefliger \pages 221--240
\paper Lissage des immersions-I
\yr 1967 \vol 6 \jour Topology
\endref

\ref \key HaQu74  \by A. Hatcher and F. Quinn \paper Bordism
invariants of intersections of submanifolds \jour Trans. Amer.
Math. Soc. \vol 200 \yr 1974 \pages 327--344 \endref

\ref \key Hud63 \by J. F. P. Hudson \paper Knotted tori \jour
Topology \vol 2 \yr 1963 \pages 11-22 \endref

\ref \key Hud69 \by J. F. P. Hudson \book Piecewise Linear
Topology \publ Benjamin \publaddr New York-Amsterdam \yr 1969
\endref




\ref \key Ki90 \by P.A. Kirk \paper Link homotopy with one codimension two component \jour Trans. Amer. Math. Soc. \vol 319:2 \yr 1990 \pages 663--688
\endref


\ref \key Kos88 \by U. Koschorke \paper Link maps and the geometry
of their invariants \yr 1988 \vol 61:4 \jour Manuscr. Math.
\pages 383--415 \endref

\ref \key Kos88L \by U. Koschorke \paper Multiple point invariants
of link maps \jour Lect. Notes Math., Springer-Verlag \vol 1350
\yr 1988 \pages 44--86 \endref


\ref \key Kos90 \by U. Koschorke \paper On link maps and their
homotopy classification \yr 1990 \vol 286:4 \jour Math. Ann.
\pages 753--782 \endref

\ref \key Le65 \by J. Levine \paper A classification of differentiable knots \jour Ann. Math. (2) \vol 82 \yr 1965 \pages 15--50 \endref


\ref \key MiRe71 \by  R. J. Milgram and E. Rees \paper On the
normal bundle to an embedding \jour Topology \vol 10 \yr 1971
\pages 299-308 \endref

\ref \key Ne84 \by V. Nezhinsky \paper A suspension sequence in link theory \jour Izv.
Akad. Nauk \vol 48:1 \yr 1984 \pages 126--143 \lang in Russian
\endref


\ref \key Pae56 \by G.F. Paechter \paper The groups $\pi_r(V_{mn})$
\jour Quart. J. Math. Oxford, Ser.2 \moreref \paper I    \vol 7:28
\yr 1956 \pages 249--268 \moreref \paper II   \vol 9:33  \yr 1958
\pages 8--27 \moreref \paper III  \vol 10:37 \yr 1959 \pages
17--37 \moreref \paper IV   \vol 10:40 \yr 1959 \pages 241--260
\moreref \paper V    \vol 11:41 \yr 1960 \pages 1--16
\endref

\ref \key ReSk99 \by D. Repov\v s and A. Skopenkov \paper New
results on embeddings of polyhedra and manifolds into Euclidean
spaces \jour Uspekhi Mat. Nauk \vol 54:6 \yr 1999 \pages 61--109
\lang in Russian \transl English transl. \jour Russ. Math. Surv.
\vol 54:6 \yr 1999 \pages 1149--1196 \endref

\ref \key Sko02 \by A. Skopenkov \paper On the Haefliger-Hirsh-Wu
invariants for embeddings and immersions \jour Comment. Math.
Helv. \vol 77 \yr 2002 \pages 78--124\endref

\ref \key Sko07F \by A. Skopenkov \paper A new invariant and parametric connected sum of embeddings \jour
Fund. Math. \vol 197 \yr 2007 \pages 253--269 \finalinfo arXiv:math/0509621[math.GT]
\endref

\ref \key Sko07L \by A. Skopenkov \paper Embedding and knotting of higher-dimensional manifolds in Euclidean spaces, in: Surveys in Contemporary Mathematics \ed N. Young and Y. Choi \jour London
Math. Soc. Lect. Notes \vol 347 \yr 2007 \pages 248--342 \finalinfo arXiv:math/0604045[math.GT] \endref


\ref \key Sko08Z \by A. Skopenkov \paper A classification of smooth embeddings of 3-manifolds in 6-space
\jour Math. Zeitschrift \vol 260 \yr 2008 \pages 647--672 \finalinfo arXiv:math/0603429v5[math.GT]
\endref

\ref \key Sko08 \by A. Skopenkov \paper Classification of embeddings below the metastable dimension
\finalinfo submitted, arXiv:math/0607422[math.GT] \endref


\ref \key Sko09 \by M. Skopenkov \paper Suspension theorems for links and link maps
\jour Proc. Amer. Math. Soc. \vol 137:1 \yr 2009 \pages 359--369 \finalinfo arXiv:math/0610320v2[math.GT]
\endref

\ref \key Sko10 \by M. Skopenkov \paper When the set of embeddings is finite? \finalinfo submitted, 	 arXiv:1106.1878v1[math.GT] \endref

\ref \key Wa70 \by C. T. C. Wall \book Surgery on Compact Manifolds \publ Academic Press \publaddr London \yr 1970
\endref

\ref \key Web67 \by C. Weber \paper Plongements de poly\`edres
dans le domain m\'etastable \jour Comment. Math. Helv. \vol 42 \yr
1967 \pages 1--27 \lang in French \endref

\endRefs
\enddocument
\end


\title{Классификация вложений торов в 2-метастабильной размерности}
\author[D.~Repov\v s]{Д.~Реповш}
\address{Institute of Mathematics, Physics and Mechanics,
and University of Ljubljana, Slovenia}
\email{dusan.repovs@guest.arnes.si}

\author[М.~Skopenkov]{М.~Скопенков}
\address{Институт проблем передачи информации им.~А.\,А.~Харкевича РАН и
King Abdullah University of Science and Technology,
Saudi Arabia}
\email{skopenkov@rambler.ru}

\author[М.~Cencelj]{М.~Ценцель}
\address{Institute of Mathematics, Physics and Mechanics,
and University of Ljubljana, Slovenia}
\email{matija.cencelj@guest.arnes.si}

\date{24.12.2011}
\udk{515.164.6}

\maketitle

\begin{fulltext}

\begin{abstract}
Статья посвящена классической проблеме заузливания: для данного многообразия $N$ и числа $m$
описать множество изотопических классов вложений $N\to S^m$. Изучается частный случай {\it заузленных торов}, то есть вложений $S^p\times S^q\to S^m$. Изотопическая классификация заузленных торов в {\it метастабильной} размерности $m\ge p +\frac {3} {2} q+2$, $p\le q$, была получена А.~Хефлигером, Е.~Зиманом и А.~Б.~Скопенковым. В работе рассматриваются размерности ниже метастабильной и 
приводится явный критерий конечности множества изотопических классов заузленных торов в {\it 2-метастабильной} размерности:

\smallskip
{\bf Теорема.} Пусть $p +\frac {4} {3} q+2 <m <
p +\frac {3} {2} q+2$ и $m> 2p+q+2$. Тогда множество изотопических классов гладких
вложений $S^p\times S^q\to S^m$  бесконечно, если и
только если хотя бы одно из чисел $q+1$ или $p+q+1$ делится на $4$.
\smallskip

Предлагаемый подход к классификации вложений основан на аналоге точной последовательности У.~Кошорке из теории сингулярных зацеплений. В данную
последовательность входит новый {\it $\beta$-инвариант} заузленных торов.
В доказательстве точности используется вложенная хирургия и метод Н.~Хабеггера--У.~Кайзера исследования дополнения.
\end{abstract}

\begin{keywords}
заузленный тор, зацепление, сингулярное зацепление, вложение, хирургия 
\end{keywords}

\markright{Классификация вложений торов}

\footnotetext[0]{Реповш и Ценцель частично поддержаны Slovenian Research Agency (гранты \No~P1-0292-0101 и J1-4144-0101). Скопенков частично поддержан грантом РФФИ \No~12-01-00748-а, грантом Президента РФ МК-3965.2012.1, фондом ``Династия'' и фондом Саймонса.
}


\section{Введение}\label{s1}

Статья посвящена классической проблеме заузливания: \emph{для данного многообразия $N$ и числа $m$
описать множество изотопических классов вложений $N\to S^m$.} Актуальные обзоры можно найти в статьях \cite{ReSk99, Sko07L}. Данная тематика активно изучалась в шестидесятые годы \cite{Hae62T, Hae66A, Hae66C, Hud63}, 
и в настоящее время интерес к ней возобновился
\cite{CeRe03, CRS04, CrSk08, Sko09}. В настоящей статье доказываются результаты, анонсированные в работах~\cite{CRS07,CRS08}.

Данная проблема обобщает предмет классической теории узлов. В отличие от 
узлов в $\mathbb{R}^3$,
в более высоких размерностях иногда удается получить
полный ответ. Мы работаем в гладкой категории, кроме тех случаев, когда явно оговаривается обратное. Перечислим некоторые известные результаты.

\subsection{Узлы}\label{s1.1}
Классификация {\it узлов} $S^q\to S^m$ в коразмерности по крайней мере $3$, то есть при $m> q+2$, сводится к некоторой гомотопической задаче \cite{Hae66A, Le65}. В частности, известна полная {\it рациональная} классификация:

\begin{theorem}\label{t0}
\cite[следствие~6.7]{Hae66A} Пусть
$q+2<m<\frac{3}{2}q+2$. Тогда множество
изотопических классов гладких вложений $S^q\to S^m$ бесконечно, если и только если $q+1$
делится на $4$.
\end{theorem}

\subsection{Зацепления} \label{s1.2}
Классификация {\it зацеплений} $S^p\sqcup S^q\to S^m$ --- следующая естественная задача после классификации узлов. При условии коразмерности по крайней мере $3$ существует
точная последовательность, содержащая множество изотопических классов зацеплений и некоторые гомотопические группы \cite{Hae66C}. В некотором диапазоне размерностей, называемом {\it
2-метастабильным}, известно явное описание изотопических
классов зацеплений $S^p\sqcup S^q\to S^m$
с точностью до узлов $S^p\to S^m$ и $S^q\to S^m$ в терминах гомотопических групп сфер и многообразий Штифеля \cite {Hae66A}. Простое доказательство этого результата приводится в \cite {Sko09}.

\centerline{

\beginpicture 
\setcoordinatesystem units <0.600mm,0.600mm>
\setplotarea x from -60 to 60, y from -40 to 20
\put {Рисунок 1} at 0 -40.000


\setcoordinatesystem units <0.600mm,0.600mm> point at 30 0
\setplotarea x from -30 to 30, y from -35 to 20


\setsolid \setplotsymbol ({\fiverm.}) \circulararc 336.650 degrees
from -7.975 -3.967 center at -0.000 -10.000 \circulararc 336.650
degrees from 7.975 -1.033 center at -0.000 5.000 \put {$S^0 \times
S^1 \to S^3$} at 0 -26.000 \put {a} at 0 -34.000

\setcoordinatesystem units <0.600mm,0.600mm> point at -30 0
\setplotarea x from -30 to 30, y from -35 to 20


\setsolid
\setplotsymbol ({\fiverm.})
\circulararc 3.423 degrees from -8.354 5.471 center at -11.020 -16.690
\circulararc 34.433 degrees from 8.944 -6.708 center at -11.020 -16.690
\circulararc 6.253 degrees from -6.422 10.241 center at -11.020 -16.690
\circulararc 36.084 degrees from 13.416 -4.472 center at -11.020 -16.690
\circulararc 36.084 degrees from -10.581 -9.383 center at -8.944 17.889
\circulararc 6.253 degrees from 12.080 0.441 center at -8.944 17.889
\circulararc 34.433 degrees from -10.282 -4.392 center at -8.944 17.889
\circulararc 3.423 degrees from 8.915 4.500 center at -8.944 17.889
\circulararc 3.423 degrees from -0.560 -9.970 center at 19.964 -1.198
\circulararc 34.433 degrees from 1.337 11.100 center at 19.964 -1.198
\circulararc 6.253 degrees from -5.658 -10.682 center at 19.964 -1.198
\circulararc 36.084 degrees from -2.835 13.855 center at 19.964 -1.198
\circulararc 33.373 degrees from 12.083 -13.748 center at 4.472 -8.944
\circulararc 77.104 degrees from 13.470 -8.770 center at 11.801 -10.531
\circulararc 39.888 degrees from 9.342 -12.449 center at 4.472 -8.944
\circulararc 42.759 degrees from 12.083 -13.748 center at 9.053 -16.600
\circulararc 180.000 degrees from -4.472 -13.416 center at 4.472 -8.944
\circulararc 180.000 degrees from -9.383 10.581 center at -9.982 0.599
\circulararc 180.000 degrees from 13.855 2.835 center at 5.510 8.345
\circulararc 180.000 degrees from 0.000 -11.180 center at 4.472 -8.944
\circulararc 180.000 degrees from -9.682 5.590 center at -9.982 0.599
\circulararc 180.000 degrees from 9.682 5.590 center at 5.510 8.345
\put {$B^2$} at 18.566 -10.728
\put {$S^1 \times S^1 \to S^3$} at 0 -26.000
\put {b} at 0 -34.000
\endpicture

}

\subsection{Заузленные торы}\label{s1.3}
В данной статье изучается классификация {\it заузленных торов}, то есть гладких вложений $S^p\times S^q\to S^m$. Эта теория обобщает теорию $2$-компонентных зацеплений одинаковой размерности; см.~рисунок~1.
Исследование заузленных торов --- следующий естественный шаг после теории зацеплений по направлению к классификации вложений произвольных многообразий \cite{Sko07F,Web67}, согласно теореме о разбиении на ручки.
Данная тематика также интересна благодаря многим интересным примерам \cite{Hud63, MiRe71, Sko02}. Ее систематическое исследование началось в работах А.Б.~Скопенкова \cite{Sko02, Sko07F, Sko08}.

Известно явное описание множества изотопических классов заузленных торов в {\it метастабильной} размерности $m\ge p +\frac {3} {2} q+2$, $p\le q$
\cite{Hae62T, Sko02} (за исключением небольшой неопределенности при $m <\frac {3} {2} p +\frac {3} {2} q+2$); см.~рисунок~2. 
Данное ограничение на размерность --- естественный предел для классических
методов классификации вложений $(p-1)$-связного $(p+q)$-мерного многообразия.
Несмотря на наличие многих интересных подходов \cite{Br68, Wa70, GW99},
мало что известно о вложениях в размерности ниже метастабильной: все известные явные классификационные результаты касаются узлов и зацеплений (они перечислены выше), заузленных торов в размерности $m=p +\frac {3} {2} q +\frac {3} {2} $
\cite{Sko08}, $3$-мерных многообразий в $\Bbb R^6$ \cite{Sko08Z} и $4$-мерных многообразий в $\Bbb R^7$ \cite{CrSk08}.

Сформулируем основной ``практический'' результат работы, анонсированный в \cite {CRS07}. Это явный критерий конечности множества изотопических классов заузленных торов в {\it 2-метастабильной} размерности (см.~заштрихованную область на рисунке~2, где число $p$ фиксировано, а числа $q$ и $m$ изменяются):

\begin{theorem}\label{tI}
Пусть $p +\frac {4} {3} q+2 <m< p +\frac {3} {2} q+2$ и $m> 2p+q+2$. Тогда множество изотопических классов гладких вложений $S^p\times S^q\to S^m$ бесконечно, если и
только если хотя бы одно из чисел $q+1$ или $p+q+1$ делится на $4$.
\end{theorem}

\begin{example}
Множество изотопических классов заузленных торов $S^1\times S^5\to
S ^ {10} $ конечно.
На рисунке~3 отмечены все размерности из диапазона $p=1$,
$1\le q\le 13$ и $m>p+\frac{4}{3}q+2$, для которых множество изотопических классов заузленных торов  бесконечно. Информация на рисунке получена из
теоремы~\ref{tI}, результатов работы \cite[теоремы 3.10 и 2.9]{Sko07L}, предложения~\ref{p5.6} и замечания~\ref{r6.3} ниже.
\end{example}

Наш подход к классификации вложений основывается на аналоге точной
последовательности Кошорке (теорема~\ref{tII} ниже), содержащем новый {\it $\beta$-инвариант почти вложений $S^p\times S^q\to S^m$}.
Точность данной последовательности доказывается с помощью метода
Хабеггера--Кайзера исследования дополнения к почти вложению.

\subsection{Структура работы} В \S2 мы формулируем основной
``теоретический'' результат --- теорему~\ref{tII}. В \S3 мы приводим определение нового {\it $\beta$-инварианта}. В \S4, который является центральным разделом статьи, мы доказываем полноту $\beta$-инва\-ри\-анта. В \S5 мы
доказываем теорему~\ref{tII}, используя результаты \S4. В \S6 мы выводим теорему~\ref{tI} из теоремы~\ref{tII}. В \S7 мы формулируем некоторые открытые вопросы. В Приложении мы детально рассматриваем вложенную хирургию многообразия
самопересечения, которая используется в \S4.

В публикации \cite{Sko10} теорема~\ref{tI} обобщена на произвольные размерности $m> 2p+q+2$ с помощью сведения классификации заузленных торов к классификации зацеплений. Необходимый для этого критерий конечности множества зацеплений получен в работе~\cite{CFS11}.

\medskip

\noindent
\begin{tabular}{cc}

\beginpicture 
\setcoordinatesystem units <0.700mm,0.700mm>
\setplotarea x from -10 to 80, y from -15 to 105


\setsolid
\setplotsymbol ({\fiverm.})
\arrow <2mm> [.2,.4] from 0.000 0.000 to 0.000 100.000
\arrow <2mm> [.2,.4] from 0.000 0.000 to 75.000 0.000

\plot 0.000 20.000 40.000 100.000 /
\plot 30.000 35.000 70.000 55.000 /
\plot 10.000 30.000 30.000 35.000 /

\setplotsymbol ({\rm.})
\plot 0.000 20.000 70.000 90.000 /

\setshadesymbol ({\rm.})
\setshadegrid span <0.05in>
\vshade 10 30 30 <z,z,,> 30 35 50 <z,z,,> 70 55 90 /

\put {$\bullet$} at 30.000 35.000
\put {0} at 15.000 94.000
\put {\tiny 1-метастабильная} at 45.000 81.000
\put {\tiny 2-метастабильная} at 51.000 55.000
\put {\tiny Классификация неизвестна} at 53.000 25.000
\put {$m$} at -2.5000 100.000
\put {$0$} at -2.5000 -2.5000
\put {$q$} at 75.000 -2.5000
\put {$p$ фиксировано} at 35.000 -5.000
\put {Рисунок 2} at 35.000 -12.500
\endpicture
 &

\beginpicture 
\setcoordinatesystem units <0.700mm,0.700mm>
\setplotarea x from -10 to 80, y from -15 to 105


\setsolid
\setplotsymbol ({\fiverm.})
\plot 60.000 0.000 60.000 95.000 /
\plot 65.000 0.000 65.000 95.000 /
\plot 5.000 0.000 5.000 95.000 /
\plot 10.000 -0.000 10.000 95.000 /
\plot 50.000 -0.000 50.000 95.000 /
\plot 55.000 0.000 55.000 95.000 /
\plot 15.000 0.000 15.000 95.000 /
\plot 20.000 -0.000 20.000 95.000 /
\plot 40.000 0.000 40.000 95.000 /
\plot 45.000 0.000 45.000 95.000 /
\plot 25.000 0.000 25.000 95.000 /
\plot 0.000 95.000 65.000 95.000 /
\plot 35.000 0.000 35.000 95.000 /
\plot 30.000 0.000 30.000 95.000 /
\plot -0.000 5.000 65.000 5.000 /
\plot -0.000 85.000 65.000 85.000 /
\plot 0.000 90.000 65.000 90.000 /
\plot 0.000 10.000 65.000 10.000 /
\plot -0.000 15.000 65.000 15.000 /
\plot 0.000 75.000 65.000 75.000 /
\plot 0.000 80.000 65.000 80.000 /
\plot 0.000 20.000 65.000 20.000 /
\plot 0.000 25.000 65.000 25.000 /
\plot 0.000 65.000 65.000 65.000 /
\plot 0.000 70.000 65.000 70.000 /
\plot -0.000 30.000 65.000 30.000 /
\plot 0.000 35.000 65.000 35.000 /
\plot 0.000 55.000 65.000 55.000 /
\plot 0.000 60.000 65.000 60.000 /
\plot 0.000 40.000 65.000 40.000 /
\plot 0.000 45.000 65.000 45.000 /
\plot 0.000 50.000 65.000 50.000 /
\arrow <2mm> [.2,.4] from 0.000 0.000 to 0.000 100.000
\arrow <2mm> [.2,.4] from 0.000 0.000 to 75.000 0.000

\setshadesymbol ({\rm.}) \setshadegrid span <0.03in>
\shaderectangleson

\putrectangle corners at 5.000 5.000 and 10.000 10.000
\putrectangle corners at 10.000 10.000 and 15.000 15.000
\putrectangle corners at 15.000 25.000 and 20.000 30.000
\putrectangle corners at 20.000 30.000 and 25.000 35.000
\putrectangle corners at 25.000 30.000 and 30.000 35.000
\putrectangle corners at 25.000 35.000 and 30.000 40.000
\putrectangle corners at 30.000 40.000 and 35.000 45.000
\putrectangle corners at 25.000 45.000 and 30.000 50.000
\putrectangle corners at 30.000 50.000 and 35.000 55.000
\putrectangle corners at 45.000 60.000 and 50.000 65.000
\putrectangle corners at 45.000 65.000 and 50.000 70.000
\putrectangle corners at 50.000 65.000 and 55.000 70.000
\putrectangle corners at 50.000 70.000 and 55.000 75.000
\putrectangle corners at 40.000 70.000 and 45.000 75.000
\putrectangle corners at 35.000 65.000 and 40.000 70.000
\putrectangle corners at 45.000 85.000 and 50.000 90.000
\putrectangle corners at 50.000 90.000 and 55.000 95.000

\shaderectanglesoff

\setplotsymbol ({\bf.}) \plot 60.000 80.000 65.000 80.000 / \plot
60.000 70.000 60.000 80.000 / \plot 55.000 70.000 60.000 70.000 /
\plot 55.000 65.000 55.000 70.000 / \plot 50.000 65.000 55.000
65.000 / \plot 50.000 60.000 50.000 65.000 / \plot 45.000 60.000
50.000 60.000 / \plot 45.000 50.000 45.000 60.000 / \plot 40.000
50.000 45.000 50.000 / \plot 40.000 45.000 40.000 50.000 / \plot
35.000 45.000 40.000 45.000 / \plot 35.000 40.000 35.000 45.000 /
\plot 30.000 40.000 35.000 40.000 / \plot 30.000 30.000 30.000
40.000 / \plot 25.000 30.000 30.000 30.000 / \plot 25.000 25.000
25.000 30.000 / \plot 20.000 25.000 25.000 25.000 / \plot 20.000
20.000 20.000 25.000 / \plot 15.000 20.000 20.000 20.000 / \plot
15.000 10.000 15.000 20.000 / \plot 10.000 10.000 15.000 10.000 /
\plot 10.000 5.000 10.000 10.000 / \plot 5.000 5.000 10.000 5.000
/ \plot 5.000 0.000 5.000 5.000 / \plot 0.000 0.000 5.000 0.000 /
\put {$m$} at -2.5000 100.000 \put {20} at -2.5000 77.500 \put
{15} at -2.5000 52.500 \put {10} at -2.5000 27.500 \put {5} at
-2.5000 2.500 \put {$q$} at 75.000 -2.5000 \put {10} at 47.500
-2.5000 \put {5} at 22.500 -2.5000 \put {1} at 2.500 -2.5000 \put
{$p=1$} at 35.000 -5.000 \put {Рисунок 3} at 35.000 -12.500
\endpicture

\end{tabular}

\bigskip

\section{Основная идея}\label{s2}

В этом разделе мы формулируем основной ``теоретический'' результат работы. Это точная последовательность, которая сводит классификацию вложений $S^p\times S^q\to S^m$ к более простой
классификации {\it почти вложений} $S^p\times S^q\to S^m$.

\subsection{Изотопии и конкордантности} Вложение $f:X\times I\to S^m\times I$ называется {\it конкордантностью}, если
$X\times 0=f ^ {-1} (S^m\times 0) $ и $X\times 1=f ^ {-1} (S^m\times 1) $.
Конкордантность называется {\it изотопией}, если
$f(X\times t)\subset S^m\times t$ для каждого $t\in I$.
Конкордантность или изотопия называется объемлемой, если $X=S^m$.
В дальнейшем мы неявно используем известные утверждения, что в коразмерности по крайней мере $3$ {\it существование конкордантости влечет существование изотопии} и {\it любая конкордантность или изотопия продолжается до объемлемой} \cite{Hud69}.

\bigskip
\noindent \begin{tabular}{l}

\beginpicture 
\setcoordinatesystem units <0.80mm,0.80mm>
\setplotarea x from -40 to 40, y from -55 to 25


\setsolid \setplotsymbol ({\fiverm.}) \circulararc 360 degrees
from 35.000 -15.000 center at 20.000 -15.000 \circulararc 360
degrees from 35.000 5.000 center at 20.000 5.000 \circulararc 360
degrees from -12.500 -25.000 center at -20.000 -25.000
\circulararc 360 degrees from -10.000 0.000 center at -20.000
0.000 \setplotsymbol ({\rm.}) \plot -26.364 6.364 -27.778 7.778 /
\plot -13.636 6.364 -12.222 7.778 / \circulararc 90.000 degrees
from -12.929 7.071 center at -20.000 0.000 \circulararc 360
degrees from -15.000 15.000 center at -20.000 15.000 \put
{$\bullet$} at 31.180 -5.000 \put {$\bullet$} at 8.820 -5.000 \put
{$\bullet$} at -20.000 10.000 \put {Почти} at -20.000 -37.000
\put {вложение} at -20.000 -42.000 \put {$f : S^0 \times S^1 \to
S^2$} at -20.000 -48.000 \put {НЕ почти} at 20.000 -37.000 \put
{вложение} at 20.000 -42.000 \put {$S^0 \times S^1 \to S^2$} at
20.000 -48.000 \put {Рисунок 4} at 0.000 -56.000 \put {$fB^1$} at
-30.000 15.000
\endpicture
\end{tabular}

\smallskip

\noindent \begin{tabular}{cc}

\beginpicture 
\setcoordinatesystem units <0.600mm,0.600mm>
\setplotarea x from -50 to 50, y from -45 to 40


\setsolid \setplotsymbol ({\fiverm.}) \plot 24.464 16.365 24.464
0.000 / \plot 21.317 17.333 21.317 0.000 / \plot 21.317 -17.333
21.317 0.000 / \plot 24.464 -16.365 24.464 0.000 / \circulararc
44.054 degrees from 21.317 -17.333 center at 40.335 -19.033
\circulararc 51.141 degrees from 24.464 -16.365 center at 40.335
-19.033 \circulararc 51.141 degrees from 35.536 16.365 center at
19.665 19.033 \circulararc 44.054 degrees from 38.683 17.333
center at 19.665 19.033 \plot 35.536 16.365 35.536 0.000 / \plot
38.683 17.333 38.683 0.000 / \plot 38.683 -17.333 38.683 0.000 /
\plot 35.536 -16.365 35.536 0.000 / \circulararc 137.364 degrees
from 37.453 -37.908 center at 38.660 -30.000 \circulararc 137.364
degrees from 15.095 -25.000 center at 21.340 -30.000 \circulararc
137.364 degrees from 37.905 -34.943 center at 38.660 -30.000
\circulararc 137.364 degrees from 17.437 -26.875 center at 21.340
-30.000 \circulararc 16.210 degrees from 32.298 -36.353 center at
40.335 -19.033 \circulararc 11.774 degrees from 34.710 -34.112
center at 40.335 -19.033 \circulararc 62.915 degrees from 22.095
-34.943 center at 19.665 -19.033 \circulararc 60.263 degrees from
22.547 -37.908 center at 19.665 -19.033 \circulararc 60.263
degrees from 44.905 -25.000 center at 30.000 -36.934 \circulararc
62.915 degrees from 42.563 -26.875 center at 30.000 -36.934
\circulararc 16.210 degrees from 19.019 -21.314 center at 30.000
-36.934 \circulararc 11.774 degrees from 19.753 -24.523 center at
30.000 -36.934 \circulararc 11.774 degrees from 40.247 24.523
center at 30.000 36.934 \circulararc 16.210 degrees from 40.981
21.314 center at 30.000 36.934 \circulararc 62.915 degrees from
17.437 26.875 center at 30.000 36.934 \circulararc 60.263 degrees
from 15.095 25.000 center at 30.000 36.934 \circulararc 60.263
degrees from 37.453 37.908 center at 40.335 19.033 \circulararc
62.915 degrees from 37.905 34.943 center at 40.335 19.033
\circulararc 11.774 degrees from 25.290 34.112 center at 19.665
19.033 \circulararc 16.210 degrees from 27.702 36.353 center at
19.665 19.033 \circulararc 137.364 degrees from 42.563 26.875
center at 38.660 30.000 \circulararc 137.364 degrees from 22.095
34.943 center at 21.340 30.000 \circulararc 137.364 degrees from
44.905 25.000 center at 38.660 30.000 \circulararc 137.364 degrees
from 22.547 37.908 center at 21.340 30.000 \circulararc 137.364
degrees from -37.453 12.092 center at -30.000 15.000 \circulararc
137.364 degrees from -37.453 37.908 center at -38.660 30.000
\circulararc 137.364 degrees from -15.095 25.000 center at -21.340
30.000 \circulararc 137.364 degrees from -34.658 13.182 center at
-30.000 15.000 \circulararc 137.364 degrees from -37.905 34.943
center at -38.660 30.000 \circulararc 137.364 degrees from -17.437
26.875 center at -21.340 30.000 \circulararc 60.263 degrees from
-22.547 12.092 center at -40.335 19.033 \circulararc 62.915
degrees from -25.342 13.182 center at -40.335 19.033 \circulararc
16.210 degrees from -32.298 36.353 center at -40.335 19.033
\circulararc 11.774 degrees from -34.710 34.112 center at -40.335
19.033 \circulararc 11.774 degrees from -35.536 16.365 center at
-19.665 19.033 \circulararc 16.210 degrees from -38.683 17.333
center at -19.665 19.033 \circulararc 62.915 degrees from -22.095
34.943 center at -19.665 19.033 \circulararc 60.263 degrees from
-22.547 37.908 center at -19.665 19.033 \circulararc 60.263
degrees from -44.905 25.000 center at -30.000 36.934 \circulararc
62.915 degrees from -42.563 26.875 center at -30.000 36.934
\circulararc 16.210 degrees from -19.019 21.314 center at -30.000
36.934 \circulararc 11.774 degrees from -19.753 24.523 center at
-30.000 36.934 \circulararc 11.774 degrees from -40.247 -24.523
center at -30.000 -36.934 \circulararc 16.210 degrees from -40.981
-21.314 center at -30.000 -36.934 \circulararc 62.915 degrees from
-17.437 -26.875 center at -30.000 -36.934 \circulararc 60.263
degrees from -15.095 -25.000 center at -30.000 -36.934
\circulararc 60.263 degrees from -37.453 -37.908 center at -40.335
-19.033 \circulararc 62.915 degrees from -37.905 -34.943 center at
-40.335 -19.033 \circulararc 16.210 degrees from -21.317 -17.333
center at -40.335 -19.033 \circulararc 11.774 degrees from -24.464
-16.365 center at -40.335 -19.033 \circulararc 11.774 degrees from
-25.290 -34.112 center at -19.665 -19.033 \circulararc 16.210
degrees from -27.702 -36.353 center at -19.665 -19.033
\circulararc 62.915 degrees from -34.658 -13.182 center at -19.665
-19.033 \circulararc 60.263 degrees from -37.453 -12.092 center at
-19.665 -19.033 \circulararc 137.364 degrees from -42.563 -26.875
center at -38.660 -30.000 \circulararc 137.364 degrees from
-22.095 -34.943 center at -21.340 -30.000 \circulararc 137.364
degrees from -25.342 -13.182 center at -30.000 -15.000
\circulararc 137.364 degrees from -44.905 -25.000 center at
-38.660 -30.000 \circulararc 137.364 degrees from -22.547 -37.908
center at -21.340 -30.000 \circulararc 137.364 degrees from
-22.547 -12.092 center at -30.000 -15.000 \setdashes \plot -45.000
-10.000 -36.245 -10.000 / \plot -23.755 -10.000 -15.000 -10.000 /
\plot -23.755 10.000 -15.000 10.000 / \plot -45.000 10.000 -36.245
10.000 / \put {$+$} at -30.000 0.000 \put {$=$} at 0.000 0.000
\put {Рисунок 5} at 0.000 -45.000
\endpicture
 & 
\beginpicture 
\setcoordinatesystem units <0.700mm,0.700mm>
\setplotarea x from -45 to 40, y from -40 to 35


\setsolid \setplotsymbol ({\fiverm.}) \arrow <2mm> [.2,.4] from
20.000 0.000 to 30.000 0.000 \arrow <2mm> [.2,.4] from -15.000
0.000 to -5.000 0.000 \plot 35.000 30.000 35.000 -30.000 / \plot
14.571 0.000 0.774 -14.001 / \plot 0.697 13.874 14.571 0.000 /
\plot -39.324 13.895 -32.773 7.344 / \circulararc 45.000 degrees
from -39.324 -13.895 center at -37.690 -15.528 \plot -30.519 5.090
-21.464 -3.964 / \plot -40.000 -15.528 -40.000 -30.000 / \plot
-39.902 15.421 -39.902 29.894 / \circulararc 12.223 degrees from
-30.524 -8.555 center at -17.402 -0.106 \circulararc 59.403
degrees from -28.438 10.929 center at -17.402 -0.106 \plot -27.357
-2.132 -39.226 -13.801 / \plot -21.367 3.858 -23.277 1.948 /
\circulararc 180.000 degrees from -21.367 3.858 center at -24.902
7.394 \circulararc 45.000 degrees from -39.902 15.421 center at
-37.592 15.421 \circulararc 180.000 degrees from -28.438 -11.142
center at -24.902 -7.606 \circulararc 45.000 degrees from 0.774
-14.001 center at 2.408 -15.634 \plot 0.000 -15.528 0.000 -30.000
/ \plot 0.000 15.528 0.000 30.000 / \circulararc 45.000 degrees
from 0.000 15.528 center at 2.310 15.528 \put {$\bullet$} at
14.571 0.000 \put {Рисунок 6} at 0.000 -38.000
\endpicture

\end{tabular}
\bigskip

\subsection{Почти вложения} Говоря неформально, {\it почти вложение} --- это отображение, имеющее только ``локальные'' самопересечения; см.~рисунок~4. Чтобы дать формальное определение,
зафиксируем отмеченную точку $*\in S^p$ и шар $B ^ {p+q} \subset S^p\times S^q$ коразмерности $0$, такой
что $B ^ {p+q} \cap (*\times S^q) = \emptyset$;
см.~рисунок~1b. Отображение $f:S^p\times S^q\to S^m$ является {\it почти вложением},
если выполнены следующие два условия:

\smallskip
\noindent  (i) $f$ является вложением вне шара $B ^ {p+q} $;

\noindent  (ii) $fB ^ {p+q} \cap f (S^p\times S^q-B ^ {p+q}) = \emptyset$.
\smallskip

\noindent  {\it Почти конкордантность} определяется аналогично, с той разницей, что шар $B ^ {p+q} $ заменяется на $B ^ {p+q} \times I$.

\subsection{Структура коммутативной группы} Операция па\-ра\-метр\-ической связной
суммы задает естественную структуру коммутативной группы на множестве вложений $S^p\times S^q\to S^m$ с точностью до конкордантности; см.~рисунок~5.
Эта структура корректно определена при $m> 2p+q+2$ \cite{Sko08}. Мы приводим ее формальное определение в \S5.

\subsection{Действие узлов на заузленных торах}
При $m> p+q+2$ множество вложений
$S ^ {p+q} \to S^m$ с точностью до конкордантности является группой относительно  операции связной суммы \cite{Hae66A}.
Эта же операция задает действие этой группы на множестве вложений
$S^p\times S^q\to S^m$ с точностью до конкордантности. В \S5 мы докажем, что это действие инъективно при $m> 2p+q+2$.
Заметим, что множество орбит данного действия находится во взаимно-однозначном соответствии с множеством вложений $S^p\times S^q\to S^m$
с точностью до конкордантности, гладкой вне некоторого конечного множества;
см.~рисунок~6.

\subsection{Обозначения}

\noindent  (a) $ {E} ^m (S^p\times S^q)/E^m (S ^ {p+q}) $ обозначает группу всех гладких вложений $S^p\times S^q\to S^m$ с точностью до конкордантности и связного суммирования с вложениями
$S ^ {p+q} \to S^m$.

\noindent  (b) $\overline {E} ^m (S^p\times S^q) $ обозначает группу всех
почти вложений $S^p\times S^q\to S^m$ с точностью до почти конкордантности;
структура группы определяется с помощью параметрической связной суммы.

\noindent  (c) $\Omega^m _ {p, q}: = \pi _ {p+2q-m+1} (V _ {N+m-p-q-1, N}) $,
где $V _ {i, j} $ --- {\it
многообразие Штифеля} $j$-осна\-ще\-ний начала координат в пространстве $\Bbb {R} ^i$ и $N$ --- достаточно большое число. Эквивалентно, при $m\ge p +\frac {4} {3} q+2$ группа
$\Omega _ {p, q} ^m$ --- это {\it группа нормальных бордизмов} $\Omega _ {2p+3q-2m+2} (P ^\infty, (m-p-q-1) \lambda) $.
Многие из групп $\Omega _ {p, q} ^m$ известны \cite{Kos88L, Pae56}.

\subsection{} Сформулируем основной ``теоретический'' результат работы.

\begin{theorem}\label{tII}
Для каждого $m\ge p +\frac {4} {3} q+2$
и $m> 2p+q+2$ существует точная последовательность
\begin{multline*}
{E} ^m (S^p\times S^q)/E^m (S ^ {p+q}) \to
\overline{E} ^m (S^p\times S^q) \xrightarrow{\beta}
\Omega _ {p, q} ^m\to\\
\to {E} ^ {m-1} (S^p\times S ^ {q-1})/E ^ {m-1} (S ^ {p+q-1}) \to \overline{E} ^{m-1} (S^p\times S^{q-1})\to\Omega _ {p, q-1} ^{m-1}\to\dots
\end{multline*}
\end{theorem}

Данная теорема имеет ряд непосредственных следствий. Во-первых,
она позволяет оценить порядок или ранг группы ${E} ^m (S^p\times S^q)/E^m (S ^ {p+q})$, если известны соответствующие оценки для группы
$\overline {E} ^m (S^p\times S^q) $. Из теоремы~\ref{tII} также несложно вывести формулу Хефлигера \cite{Hae66C} для группы зацеплений $S^q\sqcup S^q\to S^m$ в $2$-метастабильной размерности. Короткое доказательство этого
классического результата, вместе с доказательством теоремы~\ref{tII} при $p=0$,
приводится в \cite{Sko09}. Теорема~\ref{tII} при $p=0$ аналогична точной последовательности Кошорке: сравни
теоремы~3.1 и~3.5 в \cite{Sko09}, и теорему~3.1 в \cite{Kos90}.

\subsection{Бета-инвариант} Отображение
$\beta:\overline {E} ^m (S^p\times S^q) \to\Omega^m _ {p, q} $ в теореме~\ref{tII} --- это
новый инвариант, являющийся главным инструментом настоящей работы. Он обобщает:

\noindent  (a) $\beta$-инвариант сингулярных зацеплений со значениями в группе нормальных бордизмов
\cite{HaKa98, Kos88, Ki90};

\noindent  (b) $\beta$-инвариант заузленных торов \cite{Sko07F, Sko08}, сравни с \cite{Hae66C, Hud63}.

Идея этого инварианта такова. Для любого почти вложения
$f:S^p\times S^q\to S^m$ по определению $fB ^ {p+q} \cap f (*\times S^q) = \emptyset$. Наш $\beta$-инвариант измеряет ``зацепленность'' сферы
$f (*\times S^q) $ и (подходящего ретракта) полиэдра $fB ^ {p+q} $.

Основная идея работы состоит в том, что с этой точки зрения исследование почти вложений $S^p\times S^q\to S^m$ аналогично исследованию сингулярных зацеплений $S^q\sqcup S ^ {p+q} \to S^m$. Таким образом, мы можем использовать все преимущества метода Хабеггера--Кайзера теории сингулярных зацеплений для нашей задачи. Теперь сконцентрируемся на том, что сделано в работе в дополнение к методу статьи Хабеггера--Кайзера \cite{HaKa98}.

\subsection{Набросок доказательства Теоремы~\ref{tII}} Наметим доказательство точности в члене $\overline {E} ^m (S^p\times S^q) $. Нам нужно доказать, что любое почти вложение $f:S^p\times S^q\to S^m$, такое что
$\beta (f) =0$, почти конкордантно вложению.

Достаточно построить шар $B^m\subset S^m$, такой что
$f ^ {-1} B^m=B ^ {p+q} $. Действительно, тогда узел $f:\partial B ^ {p+q} \to \partial B^m$ тривиален согласно теории сглаживания.
Таким образом, можно продеформировать ограничение $f:B ^ {p+q}\to B^m$ во вложение, и получить нужное нам вложение $S^p\times S^q\to S^m$.

Для построения шара $B^m$ достаточно заклеить меридианы $f (*\times S^q) $ и
$f (S^p\times *) $ двумя дисками $D ^ {q+1} $ и $D ^ {p+1} $ (назовем их {\it паутинами}),
внутренности которых не пересекают друг друга и образ
$f (S^p\times S^q) $. Тогда нужный нам шар $B^m$ будет дополнением к малой трубчатой окрестности объединения $D ^ {q+1} \cup D ^ {p+1} $ в сфере $S^m$.

Существование паутины $D ^ {p+1} $ гарантируется неравенством $m> 2p+q+2$ и соображениями общего положения.
Покажем, как строится паутина $D ^ {q+1} $ при некоторых дополнительных предположениях. Это наиболее трудный шаг доказательства, который использует предположения $\beta (f) =0$ и $m\ge p +\frac {4} {3} q+2$. В силу результатов Хабеггера--Кайзера \cite {HaKa98}, мы можем
считать, что сфероид $f\left | _ {*\times S^q} \right.$  стягиваем вне $fB ^ {p+q} $. Таким образом, мы можем заклеить сфероид $f (*\times S^q) $ (не
обязательно вложенным) диском $D ^ {q+1} $, расположенным в пространстве $S^m-fB ^ {p+q} $. После этого мы можем устранить самопересечения этого диска, пользуясь
теоремой вложения Хадсона.

Это заканчивает доказательство при предположении, что внутренность диска $D ^ {q+1} $ не пересекает множества $f (S^p\times S^q-B ^ {p+q}) $ и этот диск --- гладкий. Остальные случаи доказываются при помощи подходящей относительной версии данного рассуждения; см.~\S4.
Подчеркнем, что в формальном доказательстве Теоремы~\ref{tII} используется не буквально приведенная конструкция, а некоторая ее относительная версия.

\section{Бета-инвариант}

В этом разделе мы приводим подробную конструкцию
$\beta$-инварианта почти вложений $S^p\times S^q\to S^m$.

\subsection{Идея $\beta$-инварианта}
Наглядно представить идею инварианта можно, обратившись к аналогии с почти вложениями $f\sqcup g:S^1\sqcup S^0\to S^2$ малой размерности.

Зафиксируем дугу $B^1\subset S^1$. Отображение
$f\sqcup g:S^1\sqcup S^0\to S^2$ назовем
{\it почти вложением}, если оно является вложением вне дуги $B^1$ и $fB^1\cap(f(S^1-B^1)\cup gS^0)=\emptyset$.
{\it Почти изотопия} $f_t\sqcup g_t:S^1\sqcup S^0\to S^2$ --- это гомотопия в классе почти вложений.

Простейший инвариант почти изотопии почти
вложения $f\sqcup g:S^1\sqcup S^0\to S^2$ --- это коэффициент зацепления $\operatorname {lk} (f, g) $, который принимает значения в группе $\Bbb Z_2$. Этот инвариант не является полным: например, почти вложение на рисунке~7 нельзя ``расцепить'' с помощью почти изотопии, хотя $\operatorname {lk} (f, g) =0$.

\smallskip

\beginpicture 
\setcoordinatesystem units <0.600mm,0.600mm>
\setplotarea x from -50 to 45, y from -35 to 35


\setsolid
\setplotsymbol ({\fiverm.})
\arrow <2mm> [.2,.4] from 5.000 0.000 to 15.000 0.000
\circulararc 360 degrees from -0.000 0.000 center at -20.000 0.000
\circulararc 263.621 degrees from 22.546 -8.333 center at 30.000 -15.000
\plot 22.546 8.333 37.454 -8.333 /
\plot 22.546 -8.333 37.454 8.333 /
\setplotsymbol ({\rm.})
\circulararc 180.000 degrees from -20.000 20.000 center at -20.000 0.000
\circulararc 263.621 degrees from 37.454 8.333 center at 30.000 15.000
\plot 30.000 -0.000 37.454 8.333 /
\plot 22.546 8.333 30.000 -0.000 /
\put {$fC_{\Delta}$} at 30.000 30.000
\put {$\Delta$} at 35.000 0.000
\put {$gS^0$} at 30.000 -20.000
\put {$gS^0$} at 30.000 20.000
\put {$f$} at 10.000 5.000
\put {$C_{\Delta}$} at -45.000 0.000
\put {$f^{-1} \Delta$} at -20.000 -25.000
\put {$f^{-1} \Delta$} at -20.000 25.000
\put {$\bullet$} at -20.000 -20.000
\put {$\bullet$} at -20.000 20.000
\put {$\bullet$} at 30.000 -15.000
\put {$\bullet$} at 30.000 15.000
\put {$\bullet$} at 30.000 0.000
\put {Рисунок 7} at 10 -33.000
\endpicture

\bigskip

В ситуации на рисунке~7 полезен следующий {\it $\beta$-инвариант} почти изотопии. Возьмем двойную точку $\Delta$
отображения $f:S^1\to S^2$. При условии общего положения множество $f ^ {-1} \Delta$ состоит из двух точек.
Соединим эти две точки дугой $C_\Delta\subset S^1$. Образ
$fC_\Delta$ --- некоторый цикл, и число
$\beta (f, g) = \sum _ {\Delta} \operatorname {lk} (fC_\Delta, g) \pmod2$,
где суммирование ведется по всем двойным точкам отображения $f$, является инвариантом почти
изотопии. Он определен корректно, только если
$\operatorname {lk} (f, g) =0$. Например, для почти вложения на рисунке~7 имеем $\beta(f,g)=1$; и это доказывает, что указанное почти вложение действительно нельзя ``расцепить''.

\subsection{Построение} Чтобы реализовать данную идею в более высоких размерностях, мы строим:

\noindent  (i) аналог цикла $fC_\Delta$, см.~определение цикла $\tilde f$ ниже;

\noindent  (ii) обобщение коэффициента зацепления
$\operatorname {lk} (fC_\Delta, f (*\times S^q)) $, см.~определение инварианта~$\beta (f) $ ниже.

\begin{definition}[данных, связанных с двойными точками] \cite{Kos88L}
Пусть $f:S^p\times S^q\to S^m$ --- некоторое почти вложение.
С помощью почти изотопии сделаем $f$ погружением общего положения. Рассмотрим диаграмму
$$
\begin{CD}
\tilde\Delta @>{\tilde i}>> B ^ {p+q} \\
@V{2:1}VV @VfVV \\
\Delta @> i>> S^m
\end{CD}
$$
Здесь $\tilde\Delta =\operatorname{Cl} \{\, (x, y) \in B ^ {p+q} \times B ^ {p+q} \, | \,
x\ne y, fx=fy \, \} $ и $\Delta =\tilde\Delta/\Bbb Z_2$ --- {\it многообразия двойных точек}, группа $\mathbb{Z}_2$ действует на $\tilde\Delta$ перестановкой сомножителей. Погружения $\tilde i:\tilde\Delta\to B ^ {p+q} $ и $i:\Delta\to S^m$ определяются формулами
$\tilde i (x, y) =x$ и $i \{x, y \} = fx$. Обозначим через
$\Sigma (f) = \operatorname {Im} \tilde i$ {\it сингулярное множество} отображения $f$.
\end{definition}

\begin{definition}[цикла $\tilde f:T (\lambda_\Delta) \to S^m$] \cite{HaKa98} Обозначим через $T (\lambda_\Delta) $ конус накрытия $\tilde\Delta\to\Delta$. Возьмем любое продолжение $\bar i:C\tilde\Delta\to B$ отображения $\tilde i:\tilde\Delta\to B$. Определим $\tilde f:T (\lambda_\Delta) \to S^m$ как фактор-отображение композиции $f\bar i:C\tilde\Delta\to S^m$.
\end{definition}

Теперь мы собираемся определить некоторое обобщение коэффициента зацепления. Это --- заявленный выше шаг (ii) построения инварианта. Продолжим вложение
$f\left | _ {*\times S^q} \right.$ до погружения
$\bar f:D ^ {q+1} \to S^m$ общего положения ({\it паутины}). Говоря неформально, нужный нам ``коэффициент зацепления'' ---
это класс бордизма ``пересечения''
$\operatorname {Im} \tilde f\cap\operatorname {Im} \bar f$ с естественным ``косым оснащением''. Он принимает значения в {\it группе нормальных бордизмов} $\Omega^m _ {p, q} $, обобщающей группу оснащенных зацеплений.

Построение всех оснащений ниже очевидно, и читатель может пропустить громоздкие шаги (2) и (3) во всех последующих определениях. Указанные шаги не используются в работе, за исключением доказательства предложения~\ref{p4.3}, по существу заимствованного из \cite{HaKa98}.

\begin{definition}[группы нормальных бордизмов $\Omega_s (P ^\infty, l\lambda) $]  \cite{Kos88L} Пусть $l$ и $s$ --- некоторые целые числа.
Назовем {\it $l\lambda$-многообразием} тройку, состоящую из

\noindent  (1) многообразия $M$ (``зацепления'');

\noindent  (2) линейного расслоения $\lambda_M$ на $M$; и

\noindent  (3) стабильного изоморфизма $\bar g_M:\nu (M) \cong
l\lambda_M: =\underbrace {\lambda_M\oplus\dots\oplus \lambda_M} _ {l} $
(``косого оснащения''), где $\nu (M) $ --- стабильное нормальное расслоение к многообразию~$M$.

\noindent  {\it Группа нормальных бордизмов} $\Omega_s (P ^\infty, l\lambda)$ --- это множество всех $l\lambda$-много\-образ\-ий размерности $s$ с точностью до
бордизма (имеющего аналогичное ``косое оснащение'').
На этом множестве операция несвязного объединения
определяет структуру коммутативной группы.
\end{definition}

Начиная с этого момента, положим $s=2p+3q-2m+2$, $l=m-p-q-1$ и $n=p+q$. Согласно \cite{Kos88L} мы имеем
$\Omega_s (P ^\infty, l\lambda) \cong\pi _ {l+s} (V _ {N+l, N}) = \Omega^m _ {p, q} $ при $s<l$ для достаточно большого целого числа $N$.
Обозначим через $N (X, Y) $ нормальное расслоение к многообразию $X$, погруженному в многообразие $Y$, и через $\epsilon$ --- тривиальное одномерное линейное расслоение.

\begin{definition}[$(m-n)\lambda$-многообразия двойных точек]
\cite{Kos88L}
Назовем {\it $(m-n)\lambda$-многообразием двойных точек} тройку $ (\Delta, \lambda_\Delta, \bar g_\Delta) $, в которой:

\noindent  (1) $\Delta$ является многообразием двойных точек;

\noindent  (2) $\lambda_\Delta$ является линейным расслоением, ассоциированным с накрытием $\tilde \Delta\to \Delta$;

\noindent  (3) $\bar g_\Delta:N (\Delta, S^m) \cong
(m-n) \lambda_\Delta\oplus\epsilon ^ {m-n} $ строится
следующим образом.

Для каждой  точки $\{x, y \}\in \Delta$ мы имеем канонические изоморфизмы
$$
N (\Delta, S^m) _ {\{x, y \}}\cong N (\tilde \Delta, B) _ {(x, y)} \oplus
N (\tilde \Delta, B) _ {(y, x)} \cong N (B, S^m) _y\oplus N (B, S^m) _x.
$$
Пусть векторы $ \{e^k_x \} _ {k=1} ^ {m-n}$ образуют тривиализацию расслоения
$N (B, S^m) $ в точке $x\in B$. Тогда векторы $ \{e^k_x,
e^k_y \} $ образуют ``косое оснащение'' многообразия $\Delta$. Если мы поменяем местами точки $x$ и $y$ в паре
$(x,y)\in\tilde\Delta$, то векторы $e^k_x$
и $e^k_y$ поменяются местами. Таким образом, расслоение $N (\Delta, S^m) $ распадается в сумму
всех линейных расслоений $\langle e^k_x+e^k_y\rangle\cong\epsilon$ и
линейных расслоений $\langle e^k_x-e^k_y\rangle\cong\lambda_\Delta$, где $k=1, \dots, m-n$. Данное разложение определяет нужный нам
изоморфизм $\bar g_\Delta:N (\Delta, S^m) \cong
(m-n) \lambda_\Delta\oplus \epsilon ^ {m-n} $.
\end{definition}

Заметим, что пространство расслоения $\lambda_\Delta$ можно отождествить с подмножеством конуса
$T (\lambda_\Delta) $ накрытия $\tilde\Delta\to\Delta$.
Будем обозначать ограничение цикла $\tilde f:T (\lambda_\Delta) \to S^m$ на данное подмножество также через $\tilde f:\lambda_\Delta\to S^m$. Обозначим через $\check f:\lambda_\Delta\to S^m$ гладкое отображение общего положения (не обязательно погружение), достаточно близкое к $\tilde f:\lambda_\Delta\to S^m$.

\begin{definition}[бета-инварианта $\beta (f) $]
Назовем бета-инвариантом $\beta (f) $ почти вложения $f:S^p\times S^q\to S^m$
класс бордизма $l\lambda$-многообразия $(\beta,\lambda_\beta,\bar g_\beta)$, определяемого следующим образом:

\noindent  (1) {\it Многообразие $\beta =\beta (f) $.} Положим
$$
\beta = \{\, (x, y) \in D ^ {q+1} \times \lambda_\Delta \, : \, \bar fx =\check fy
\, \}.
$$

\noindent  (2) {\it Линейное расслоение $\lambda_\beta$ на многообразии $\beta$.}
Обозначим $\operatornamewithlimits {pr}:\beta\to \Delta$
очевидную композицию $\beta\to D ^ {q+1} \times \lambda_\Delta\to
0\times \Delta =\Delta$. Положим
$\lambda_\beta =\operatornamewithlimits {pr ^ *} (\lambda_\Delta) $.

\noindent  (3) {\it Стабильный изоморфизм $\bar g_\beta:\nu (\beta) \cong l\lambda_\beta$.} Ограничим стабильный изоморфизм $\bar g_\Delta:\nu (\Delta) \cong (l+1) \lambda_\Delta$, построенный в шаге (3) определения
$(m-n) \lambda$-многообразия двойных точек, на расслоение
$ <e^1_x-e^1_y> ^ \perp$. Мы получим изоморфизм
$\nu (\lambda_\Delta) \cong l\lambda_\Delta\oplus \epsilon ^ {l+1} $.
Отождествим $\nu (\lambda_\Delta) $ и $\nu (D ^ {q+1} \times\lambda_\Delta) $.
Ограничивая предыдущий изоморфизм на многообразие $\beta$, мы получаем изоморфизм
$g_1:\nu (D ^ {q+1} \times\lambda_\Delta) \left | _ {\beta} \right. \cong
l\lambda_\beta\oplus \epsilon ^ {l+1} $. Возьмем тривиализацию нормального расслоения $N (D ^ {q+1}, S^m) $. Эта тривиализация задает
изоморфизм $g_2:\nu (\beta) \cong
\nu (D ^ {q+1} \times\lambda_\Delta) \left | _ {\beta} \right.\oplus \epsilon ^ {m-q-1} $. Положим $\bar g_\beta = (g_1\oplus \operatorname {id}) \circ
g_2$.
\end{definition}

\begin{propos}\label{p3.1} Формула $f\mapsto\beta (f) $
корректно определяет некоторое отображение $\overline {E} ^m (S^p\times S^q) \to \Omega^m _ {p, q} $.
\end{propos}

\begin{proof} Нам нужно проверить следующее:

(1) {\it Класс бордизма $\beta (f)$ не зависит от выборов в его построении.} Действительно, в процессе построения мы делали следующие четыре выбора. В определении цикла $\tilde f$ мы выбирали продолжение $\bar
i:C\tilde \Delta\to D ^ {p+q} $. В шаге (3) определения
$ (m-n) \lambda$-многообразия двойных точек мы выбирали тривиализацию расслоения $N (D ^ {p+q}, S^m) $. В шаге
(1) определения инварианта $\beta (f) $  мы выбирали продолжение $\bar f:D ^ {q+1} \to S^m$. И в шаге (3)
того же определения мы выбирали тривиализацию расслоения
$N (D ^ {q+1}, S^m) $. Ясно, что все рассматриваемые продолжения и тривиализации единственны с точностью до гомотопии. Следовательно, класс $\beta (f) $ определен корректно.

(2) {\it Если $f_1$ и $f_2$ почти конкордантны, то
$\beta (f_1) = \beta (f_2) $.} Действительно, пусть $f:S^p\times S^q\times I\to S^m\times I$ --- почти конкордантность общего положения между $f_1$
и $f_2$. Аналогично построению выше можно построить
$\beta$-инвариант $\beta (f)$ этой почти конкордантности. Он будет бордизмом между $\beta (f_1) $ и $\beta (f_2) $. Таким образом,
$\beta (f_1) = \beta (f_2) $.
\end{proof}

Приведем относительную версию предыдущего построения; см.~рисунок~8. Отображение $f:X\to Y$ назовем {\it собственным}, если $f^{-1}\partial Y=\partial X$.

\begin{definition}[собственного почти вложения $S^p\times D^q\to D^m$]
Зафиксируем шар $B ^ {p+q} \subset S^p\times \operatorname {Int} D^q$ такой, что $B ^ {p+q} \cap (*\times D^q) = \emptyset$.
Собственное отображение $f:S^p\times
D^q\to D^m$ будем называть {\it собственным почти вложением}, если выполнены следующие условия:

\smallskip

\noindent  (i) $f$ является вложением вне $B ^ {p+q} $; и

\noindent  (ii) $fB ^ {p+q} \cap f (S^p\times D^q-B ^ {p+q}) = \emptyset$.

\smallskip

\noindent  {\it Собственная почти конкордантность} определяется аналогично, с той разницей, что шар $B ^ {p+q} $ заменяется $B ^ {p+q} \times I$.
{\it Стандартное} вложение $S^p\times D^q\to D^m$ определяется как композиция $S^p\times D^q\subset D ^ {p+1} \times D^q\cong D ^ {p+q+1} \subset D^m$.
\end{definition}

Далее $f:S^p\times D^q\to D^m$ будет обозначать некоторое собственное почти вложение, если явно не оговорено обратное.

\smallskip

\beginpicture 
\setcoordinatesystem units <0.700mm,0.700mm>
\setplotarea x from -30 to 30, y from -35 to 30


\setsolid \setplotsymbol ({\fiverm.}) \circulararc 13.818 degrees
from -2.502 13.124 center at 12.205 -20.000 \plot -7.352 10.526
-8.325 8.059 / \plot -5.500 11.635 -7.165 8.001 / \plot -3.839
12.443 -5.800 8.288 / \circulararc 120.000 degrees from -10.000
8.660 center at -7.500 12.990 \circulararc 120.000 degrees from
-15.000 -0.000 center at -7.500 12.990 \circulararc 60.000 degrees
from -12.500 8.660 center at -8.750 2.165 \circulararc 120.000
degrees from -8.750 6.495 center at -11.250 6.495 \circulararc
60.000 degrees from -12.500 -0.000 center at -16.250 6.495
\circulararc 120.000 degrees from -16.250 2.165 center at -13.750
2.165 \circulararc 120.000 degrees from 6.250 10.825 center at
5.000 12.990 \circulararc 60.000 degrees from -1.250 10.825 center
at 2.500 17.321 \circulararc 120.000 degrees from -1.250 15.155
center at 0.000 12.990 \circulararc 60.000 degrees from 6.250
15.155 center at 2.500 8.660 \circulararc 360 degrees from 25.000
0.000 center at 0.000 0.000 \circulararc 42.261 degrees from
-7.221 -22.420 center at -3.770 -21.110 \circulararc 42.261
degrees from -2.167 -23.454 center at 0.917 -21.424 \circulararc
42.261 degrees from -11.929 -20.310 center at -8.276 -19.783
\circulararc 42.261 degrees from 8.007 -22.151 center at 9.935
-19.003 \circulararc 42.261 degrees from 2.992 -23.363 center at
5.559 -20.711 \circulararc 42.261 degrees from 24.400 -5.443
center at 20.711 -5.559 \circulararc 42.261 degrees from 22.630
-10.625 center at 19.003 -9.935 \circulararc 42.261 degrees from
22.630 10.625 center at 19.783 8.276 \circulararc 42.261 degrees
from 25.000 0.000 center at 21.424 -0.917 \circulararc 42.261
degrees from 24.400 5.443 center at 21.110 3.770 \circulararc
42.467 degrees from 1.717 16.119 center at -0.431 15.019
\circulararc 106.543 degrees from -1.251 10.825 center at -1.366
10.264 \circulararc 102.259 degrees from 0.915 9.990 center at
0.804 9.131 \circulararc 102.259 degrees from 5.059 9.058 center
at 4.196 9.131 \circulararc 106.543 degrees from 6.937 10.214
center at 6.366 10.264 \circulararc 42.467 degrees from 4.589
17.281 center at 5.431 15.019 \circulararc 72.093 degrees from
8.982 13.426 center at 8.513 12.210 \circulararc 39.567 degrees
from 7.691 16.124 center at 8.082 13.431 \circulararc 39.567
degrees from -1.065 15.258 center at -3.082 13.431 \circulararc
39.567 degrees from -12.619 10.392 center at -10.091 9.384
\circulararc 39.567 degrees from -16.247 2.376 center at -15.673
-0.284 \circulararc 72.093 degrees from -15.035 0.020 center at
-14.830 -1.268 \circulararc 42.467 degrees from -15.601 5.216
center at -15.723 2.806 \circulararc 106.543 degrees from -12.501
-0.001 center at -12.072 -0.382 \circulararc 102.259 degrees from
-10.694 1.457 center at -10.006 0.931 \circulararc 102.259 degrees
from -7.815 4.580 center at -8.310 3.869 \circulararc 106.543
degrees from -7.877 6.785 center at -8.206 6.315 \circulararc
42.467 degrees from -15.171 8.284 center at -12.792 7.883
\circulararc 24.761 degrees from -9.093 -1.925 center at -11.752
2.006 \circulararc 46.460 degrees from -3.879 -1.566 center at
-5.004 -0.186 \circulararc 39.363 degrees from 1.677 1.125 center
at -0.364 2.636 \put {$f(S^p \times D^q)$} at 0.000 -5.000 \put
{$\overline{f}D^{q+1}$} at -10.000 15.000 \put {$D^m$} at 23.000
20.000 \put {Рисунок 8} at 0.000 -30.000
\endpicture

\bigskip

\begin{definition}[собственного почти вложения $D ^ {p+q} \to D^m$]
{\it Собственное почти вложение} $D ^ {p+q} \to D^m$ --- это собственное отображение, ограничение которого на край является вложением.
\end{definition}

\begin{definition}[паутины $\bar f$]
Зафиксируем стандартное экваториальное разбиение $\partial D ^ {q+1} =D^q _ +\cup_{\partial D^q_+=\partial D^q_-=\partial D^q} D^q_-$. {\it Паутина} собственного почти вложения
$f:S^p\times D^q\to D^m$ --- это отображение $\bar
f:D ^ {q+1} \to D^m$, удовлетворяющее следующим двум условиям, см.~рисунок~8:

\smallskip
\noindent  (i) $\bar f\left | _ {D^q _ +}\right. = f\left | _ {*\times D^q} \right.$; и

\noindent  (ii) $\bar fD^q_-\subset \partial D^m$.
\end{definition}

\begin{definition}[относительных бета-инвариантов $\beta (f) $ и $\beta (f, g) $] Определение (относительного)
{\it $\beta$-инварианта $\beta (f) $} собственного почти вложения $f:S^p\times D^q\to D^m$ полностью аналогично определению инварианта $\beta (f) $ выше, с той разницей, что отображение $f\left | _ {*\times S^q} \right.$
заменяется на $f\left | _ {*\times D^q} \right.$.

Если дано собственное отображение
$g:D^q\to D^m$, образ которого не пересекает $f (S^p\times D^q) $, то можно определить {\it
$\beta$-инвариант $\beta (f, g) $} аналогично инварианту~ $\beta (f) $, с той разницей, что
отображение $f\left | _ {*\times D^q} \right.$ заменяется отображением~$g$.
\end{definition}

В качестве непосредственного следствия определений получаем:

\begin{propos}\label{p3.2} Отображение $g\mapsto \beta (f, g) $ индуцирует
гомоморфизм групп $\pi_q (D^m-\operatorname {Im} f, \partial D^m-\operatorname {Im} f) \to\Omega^m _ {p, q} $. Кроме того, если собственное отображение $g:D^q\to D^m$ достаточно близко к отображению $f\left | _ {*\times D^q} \right.$, то $\beta (f, g) = \beta (f) $.
\end{propos}

\section{Полнота бета-инварианта}\label{s4}

В этом разделе мы доказываем полноту относительного
$\beta$-инварианта:

\begin{theorem}\label{t4.1} Пусть $m\ge p+\frac {4}{3} q+2$. Тогда любое собственное почти вложение $f:S^p\times D^q\to D^m$, для которого $\beta (f) =0$, собственно почти конкордантно связной сумме \textup{(}относительно края\textup{)}
стандартного вложения $S^p\times D^q\to D^m$ и некоторого собственного почти вложения $D ^ {p+q} \to D^m$.
\end{theorem}

Сначала сформулируем нашу центральную лемму, которая описывает
гомотопические группы дополнения к собственному почти вложению; сравни с \cite[следствие~4.4]{HaKa98}:

\begin{lemma}[о дополнении]\label{l3}
Пусть
$m\ge p +\frac {4} {3} q+2$. Тогда любое собственное почти вложение $f:S^p\times D^q\to D^m$
собственно почти конкордантно собственному почти вложению $f ':S^p\times D^q\to D^m$ общего положения, такому что
$$
\pi _ {q} (D^m-\operatorname {Im} f', \partial D^m-\operatorname {Im} f') \cong\Omega^m _ {p, q}.
$$
Последний изоморфизм задается формулой $g\mapsto\beta (f', g) $.
\end{lemma}

В частности, если $\beta (f) =0$, то любое собственное отображение $g:D^q\to D^m-\operatorname {Im} f '$, достаточно близкое к $f '\left | _ {*\times D^q} \right.$, собственно гомотопно постоянному отображению (в классе собственных отображений $D^q\to D^m-\operatorname {Im} f '$). Это замечание следует из предложения~\ref{p3.2} выше
и составляет основу следующего рассуждения.

\begin{proof}[теоремы \ref{t4.1}, исходя из Леммы~\ref{l3}] Доказательство проходит в $3$ этапа:

(1) {\it Построение паутины,
внутренность которой не пересекает $\operatorname {Im} f$.} Возьмем собственное почти вложение $f:S^p\times D^q\to D^m$. Пусть $f ':S^p\times D^q\to D^m$ --- собственное почти вложение (собственно почти конкордантное отображению $f$), такое что выполнен изоморфизм из леммы~\ref{l3} о дополнении.

Без потери общности можно считать, что
$f' (S^p\times D^q) $ ортогонально границе $\partial D^m$. Так как ограничение нормального расслоения $N (f' (S^p\times D^q), D^m) $ на диск $f'(*\times D^q)$
тривиально, то существует единичное векторное поле на диске $f'(*\times D^q)$, ортогональное образу $f '(S^p\times D^q)$.
Приклеим к диску $f' (*\times D^q) $ воротниковую окрестность в направлении этого векторного поля.
Компонента края этой воротниковой окрестности
образует собственное отображение $g:D^q\to D^m-\operatorname {Im} f '$.

Согласно лемме~\ref{l3} о дополнении (и абзацу после ее формулировки) отображение $g:D^q\to D^m-\operatorname {Im} f '$ собственно гомотопно постоянному отображению.
Значит, к рассматриваемой компоненте края воротниковой окрестности можно приклеить диск (возможно, самопересекающийся), не пересекающий $\operatorname {Im} f '$, край которого лежит в $\operatorname {Im}g\cup \partial D^m$. Объединение диска и воротниковой окрестности --- образ искомой паутины $\bar f:D ^ {q+1} \to D^m$.

(2) {\it Устранение самопересечений паутины.} Согласно \cite[Theorem~2.1]{Ha67} существует
кусочно-гладкий гомеоморфизм $h:D^m\to D^m$, такой что отображения $f' _ {PL} :=hf' $
и $\bar f _ {PL} :=h\bar f$ кусочно линейны. Обозначим через $M^m$
дополнение к регулярной окрестности образа $\operatorname {Im} f' _ {PL} $ в $D^m$. Тогда пара $(M^m, \partial M^m\cap\partial D^m) $
достаточно высоко связна (см.~предложение~\ref{p4.5} ниже).
Следующая теорема позволяет устранить самопересечения
паутины:

\begin{theorem}[о вложении, допускающая движение части края]\label{t4.2}
Пусть $ (M^m, M^ {m-1}) $ --- пара
кусочно линейных многообразий, такая что $M ^ {m-1} \subset \partial M^m$. Предположим, что эта пара
$ (2q-m+3)$-связна и $ m \ge q+4$. Пусть $\bar
f _ {PL}: (D ^ {q+1}, D ^ {q} _-)\to (M^m, M^ {m-1}) $ --- кусочно линейное отображение,
которое вкладывает $D ^ {q} _ +$ в $\partial M^m-\operatorname {Int} M ^ {m-1} $. Тогда
$\bar f _ {PL} $ гомотопно $\operatorname {rel} D ^ {q} _ +$ кусочно линейному вложению
$\bar f _ {Emb}: (D ^ {q+1}, D ^ {q} _-)\to (M^m, M^ {m-1}) $.
\end{theorem}

Эта теорема доказывается полностью аналогично результату \cite[теорема~9.2.1]{Hud69}.
Так как все препятствия к сглаживаемости вложения $h ^ {-1} \bar f _ {Emb} $ лежат в нулевых группах $H^k(D^{q+1},D^{q}_-;C^{m-q}_{q-k})$, то
согласно теории сглаживания \cite{Ha67}, кусочно линейное вложение $h ^ {-1} \bar f _ {Emb} $ собственно гомотопно $\operatorname {rel} D ^ {q} _ +$
паутине $\bar f ':D ^ {q+1} \to D^m$, которая уже является гладким вложением.

(3) {\it Разложение $f '$ в связную сумму.} Пусть $B^m$ --- дополнение к объединению трубчатых окрестностей многообразия $f' (S^p\times D^q-\operatorname {Int} B ^ {p+q}) $ и паутины $\bar f'D ^ {q+1} $ в $D^m$.
Тогда $B^m$ --- гладкий шар, причем $ (f') ^ {-1} B^m=B ^ {p+q} $.
Обозначим ограничение $f ':B ^ {p+q} \to B^m$ через $g':B^{p+q}\to B^m$. Легко видеть, что
$f '$ собственно почти конкордантно связной сумме (относительно края) отображения $g':B^{p+q}\to B^m$ и некоторого вложения $S^p\times D^q\to D^m$.
Последнее собственно объемлемо конкордантно стандартному
(согласно несложному предложению~\ref{p5.7}(a) ниже).
Теорема~\ref{t4.1} доказана, исходя из леммы~\ref{l3}.
\end{proof}

Оставшаяся часть \S4 посвящена доказательству леммы~\ref{l3} о дополнении. Наше рассуждение параллельно
\cite[\S\S3--4]{HaKa98}.

Начиная с этого момента, отождествим диск $D^m$ с верхней полусферой сферы $S^m$
и зафиксируем разбиение $S^m=D^m\cup CS ^ {m-1} $.
Обозначим через $ \{X, Y \} $ множество стабильных гомотопических классов отображений $X\to Y$.

\begin{example}
$\{ (D^q; S^{q-1})\,,\,(D^m-\operatorname{Im}f,\partial)\}\cong\{S^q,S^m-f(S^p\times D^q)\}$
(потому что пара $ (X, Y) $ стабильно гомотопически эквивалентна пространству $X\cup CY$).
\end{example}

Следующее предложение показывает, как выразить $\beta$-инвариант через гомотопический класс собственного отображения $g:D^q\to D^m-\operatorname {Im} f$
в группе $\pi_q (D^m-\operatorname {Im} f, \partial) $.
Данное предложение можно рассматривать как альтернативное определение $\beta$-инварианта $\beta (f, g) $.

\begin{propos} \label{p4.3} \textup{(сравни с \cite[предложение~3.2]{HaKa98})}
При композиции
\begin{multline}
\pi_q (D^m-\operatorname {Im} f, \partial)
\xrightarrow{\Sigma^\infty}
\{S^q, S^m-\operatorname {Im} f \}
\xrightarrow{SW} \\
\xrightarrow{SW}
\{\operatorname {Im} f, S ^ {m-q-1} \}
\xrightarrow{\tilde f_ *}
\{T (\lambda_\Delta), S ^ {m-q-1} \}
\xrightarrow{PT}
\Omega_s (P ^ {\infty}, l\lambda)
\end{multline}
гомотопический класс собственного отображения $g:D^q\to D^m-\operatorname {Im} f$ переходит в $\beta (f, g) $.
\end{propos}

Здесь первая стрелка --- итерация отображения надстройки. Второе отображение --- это двойственность Спениера--Уайтхеда. Третья стрелка индуцируется отображением
$\tilde f$, определенным в \S3. Четвертая стрелка задается конструкцией Понтрягина--Тома (см.~детали в \cite{Kos88L}). Данное предложение доказывается непосредственной проверкой. Фактически, его достаточно доказать для отображения $f:B ^ {p+q} \to S^m$, что делается полностью аналогично \cite[предложение~3.2]{HaKa98}.

Данное предложение подсказывает, что полезно найти гомотопический тип образа $\operatorname {Im} f$:

\begin{propos} \label{p4.4} \textup{(сравни с \cite[\S4]{HaKa98})} Обозначим через $C$ конус отображения $f:\Sigma (f) \to f\Sigma (f) $. Тогда $\operatorname {Im} f\simeq C\vee S^p$.
\end{propos}

Доказательство этого предложения немедленно следует из следующего наблюдения: оба рассматриваемых пространства могут быть получены из полиэдра $\operatorname {Cyl} (\Sigma (f) \to f\Sigma (f)) \cup _ {\Sigma (f) \subset S^p\times D^q} (S^p\times D^q) $ с помощью стягивания подходящих подполиэдров; см.~рисунок~9.

\beginpicture 
\setcoordinatesystem units <0.600mm,0.600mm>
\setplotarea x from -60 to 60, y from -30 to 30


\arrow <2mm> [.2,.4] from -5.000 0.000 to 5.000 0.000
\put {$f$} at 0.000 5.000
\put {Рисунок 9} at 0.000 -25.000

\setcoordinatesystem units <0.600mm,0.600mm> point at 30 0
\setplotarea x from -30 to 30, y from -30 to 30


\setsolid
\setplotsymbol ({\fiverm.})
\circulararc 360 degrees from 20.000 0.000 center at 0.000 0.000
\setplotsymbol ({\rm.})
\plot -16.454 9.500 -18.187 10.500 /
\plot 16.454 9.500 18.187 10.500 /
\circulararc 120.000 degrees from 17.321 10.000 center at -0.000 0.000
\setdashes
\setplotsymbol ({\fiverm.})
\plot 0.000 0.000 -14.142 14.142 /
\plot 0.000 0.000 0.000 20.000 /
\plot 0.000 0.000 14.142 14.142 /
\put {$B^{p+q}$} at -20.000 18.000
\put {$\Sigma(f)$} at 0.000 25.000
\put {Cyl} at 0.000 -5.000
\put {$\bullet$} at -14.142 14.142
\put {$\bullet$} at 14.142 14.142
\put {$\bullet$} at 0.000 20.000
\put {$\bullet$} at 0.000 0.000

\setcoordinatesystem units <0.600mm,0.600mm> point at -30 0
\setplotarea x from -30 to 30, y from -30 to 30


\setsolid
\setplotsymbol ({\fiverm.})
\plot -4.330 7.500 4.330 -7.500 /
\plot 4.330 7.500 -4.330 -7.500 /
\plot -8.660 0.000 8.660 0.000 /
\circulararc 240.000 degrees from -8.660 -0.000 center at -8.660 -5.000
\circulararc 240.000 degrees from 4.330 -7.500 center at 8.660 -5.000
\circulararc 240.000 degrees from 4.330 7.500 center at 0.000 10.000
\put {Im $f$} at 0.000 -15.000
\put {$\bullet$} at 0.000 -0.000

\endpicture

\bigskip

Начнем теперь исследование гомотопического типа пары
$ (D^m-\operatorname {Im} f, \partial) $.

\begin{propos}\label{p4.5} \textup{(сравни с \cite[лемма~4.2]{HaKa98})} Для любого погружения $f$ общего положения пара
$ (D^m-\operatorname {Im} f, \partial) $ является $c$-связной, где
$c =\min \{m-p-2,2m-2p-2q-3 \} $.
\end{propos}

\begin{proof} При условии коразмерности по крайней мере 3 пара
$ (D^m-\operatorname {Im} f, \partial) $ по общему положению односвязна.
Согласно теореме вырезания имеем $H_i (D^m-\operatorname {Im} f, \partial)
\cong H_i (S^m-\operatorname {Im} f) $. По двойственности Александера
$H_i (S^m-\operatorname {Im} f) \cong H ^ {m-i-1} (\operatorname {Im} f) $.
Согласно предложению~\ref{p4.4} мы имеем $H ^ {m-i-1} (\operatorname {Im} f) =0$ для $i\le c$,
так как по общему положению $\operatorname {dim} C=2p+2q-m+1$. По теореме Гуревича пара $(D^m-\operatorname {Im} f, \partial) $ является
$c$-связной.
\end{proof}

\begin{proof}[леммы~\ref{l3} о дополнении] Достаточно показать, что при подходящем выборе собственного почти вложения $f '$ все стрелки в композиции из предложения~\ref{p4.3} биективны.
Нам потребуется изменить исходное собственное почти вложение $f$ только в шаге (4) ниже.

(1) {\it Первое отображение биективно} согласно предложению~\ref{p4.5} и теореме о надстройке, так как предположение $m\ge p +\frac {4} {3} q+2$ влечет
неравенство $q/2\le c$.

(2) {\it Второе отображение биективно} по двойственности Спениера--Уайт\-хеда.

(3) {\it Третье отображение биективно.} Согласно предложению~\ref{p4.4} получаем
$\{\operatorname {Im} f, S^{m-q-1} \} 
\cong \{C, S ^ {m-q-1} \} $. Остается проверить, что
$\{C, S ^ {m-q-1} \} \cong \{T (\lambda_\Delta), S ^ {m-q-1} \} $. Обозначим через
$\lambda_T$ ограничение расслоения $\lambda_\Delta$ на тройные
точки, через $C_T$ --- образ пространства $T (\lambda_T) $ при отображении
$T (\lambda_\Delta) \to C$, и через $R$ --- конус отображения $T (\lambda_T) \to C_T$. Тогда $R$ --- деформационный ретракт конуса отображения $T (\lambda_\Delta) \to C$. Рассмотрим точную последовательность Пуппе пары
$ (\operatorname {Cyl} (T (\lambda_\Delta) \to C), T (\lambda_\Delta)) $:
$$
\{R, S ^ {m-q-1} \} \to \{C, S ^ {m-q-1} \}
\to \{T (\lambda_\Delta), S ^ {m-q-1} \} \to \{R, S ^ {m-q} \}.
$$
Так как $R$ имеет размерность не более $3q+3p-2m+2$, то ввиду предположения $m\ge
p +\frac {4} {3} q+2$ получаем, что $ \{C, S ^ {m-q-1} \} \cong \{T (\lambda_\Delta), S ^ {m-q-1} \} $. Ясно, что полученный изоморфизм  $\{\operatorname {Im} f, S^{m-q-1} \} \cong \{T (\lambda_\Delta), S ^ {m-q-1} \} $
индуцируется отображением $\tilde f:T (\lambda_\Delta)\to\operatorname {Im} f$.

(4) {\it Превращение четвертого отображения в биективное.} Легко видеть, что четвертое отображение биективно, если отображение $\Delta\to P ^\infty$, классифицирующее
расслоение $\lambda_\Delta$, является $ (s+1)$-связным. Действительно, согласно конструкции Понтрягина--Тома $\{T (\lambda_\Delta), S ^ {m-q-1} \} \cong
\Omega_s (\Delta, l\lambda_\Delta) $, и по теореме о накрывающей гомотопии индуцированное отображение
$\Omega_s (\Delta, l\lambda_\Delta) \to\Omega_s (P ^\infty, l\lambda) $ будет изоморфизмом. (Определение группы $\Omega_s (\Delta, l\lambda_\Delta)$ и детали данного рассуждения приводятся в
\cite[\S3]{HaKa98} и \cite {Kos88L}).

Таким образом, нам остается сделать классифицирующее отображение $\Delta\to P ^\infty$
$ (s+1)$-связным при помощи подходящей собственной почти конкордантности почти вложения $f$.
Это возможно в силу следующей теоремы (которая применяется здесь при $s=2p+3q-2m+2$ и $n=p+q$):

\begin{theorem}[о хирургии]\label{t4.6} \textup{(сравни с \cite[теорема~4.5]{HaKa98})}
Рассмотрим $(s+1)$-связное многообразие $M^m$ и собственное погружение $f:B^n\to M^m$, ограничение $f\left | _ {\partial B^n} \right.$ которого на край --- вложение.
Предположим, что $2s \le 2n-m-2$ и $0\le s\le m-n-3$.
Тогда с помощью некоторой регулярной
гомотопии  $\operatorname {rel} \partial B^n$ погружения $f:B^n\to M^m$ классифицирующее отображение
$\Delta\to P ^\infty$ накрытия $\tilde\Delta\to\Delta$ можно сделать $ (s+1)$-связным.
\end{theorem}

Доказательство теоремы~\ref{t4.6} о хирургии полностью аналогично доказательству \cite[теорема~4.5]{HaKa98}. Мы приводим его в Приложении для удобства читателя.
\end{proof}

Таким образом, мы доказали полноту $\beta$-инварианта, считая теорему~\ref{t4.6} о хирургии доказанной.

\section{Точная последовательность}\label{s5}

В данном разделе мы выводим теорему~\ref{tII} из утверждения о полноте $\beta$-ин\-ва\-ри\-ан\-та. Формально, теорема~\ref{tII} следует из утверждений~\ref{t4.1}, \ref{t5.1}, \ref{t5.2} и~\ref{p5.5}.

\begin{theorem}\label{t5.1} Пусть $m\ge p +\frac {4} {3} q+2$. Тогда для каждого $x\in\Omega^m _ {p, q} $ существует собственное почти
вложение $\omega_x:S^p\times D^q\to D^m$, такое что $\beta (\omega_x) =x$.
\end{theorem}

\begin{proof}[теоремы~\ref{t5.1}, исходя из леммы~\ref{l3}]
Построение элемента
$\omega_x$ состоит из $3$ шагов:

(1) {\it Построение отображения $f '\sqcup g:S^p\times D^q\sqcup D^q\to D^m$, такого что $\beta (f', g) =x$.} Начнем со стандартного
вложения $f:S^p\times D^q\to D^m$. По лемме~\ref{l3} о дополнении из \S4
отображение $f$ собственно почти конкордантно собственному почти вложению $f ':S^p\times D^q\to D^m$, такому что отображение
$\beta (f', {\cdot}):\pi_q (D^m-\operatorname {Im} f', \partial) \to
\Omega^m _ {p, q} $ является изоморфизмом.
Возьмем собственное отображение $g:D^q\to
D^m-\operatorname {Im} f '$, такое что $\beta (f', g) =x$.

(2) {\it Гомотопическое продолжение отображения $g:D^q\to D^m-\operatorname {Im} f '$ до собственного вложения
$ g':S^p\times D^q\to
D^m-\operatorname {Im} f '$.}
Так как пара
$ (D^m-\operatorname {Im} f', \partial) $ достаточно высоко связна (см.~ предложение~\ref{p4.5}), то по теореме вложения \cite[теорема~8.2.1]{Hud69} самопересечения
отображения $g:D^q\to
D^m-\operatorname {Im} f '$ можно устранить.
Согласно теории Хирша можно сделать $g:D^q\to
D^m-\operatorname {Im} f '$ гладким вложением с тривиальным нормальным расслоением. Значит, можно продолжить вложение $g:D^q\to D^m-\operatorname {Im} f '$ до собственного вложения
$g':S^p\times D^q \to D^m-\operatorname {Im} f '$, такого что образ
$\operatorname {Im} g '$ содержится в трубчатой окрестности образа $\operatorname {Im} g$.

(3) {\it $S^p$-параметрическая связное суммирование отображений $f '$ и $g' $.}
Зафиксируем отмеченную точку $*\in \partial D ^ {q} $.
По построению
$f '\left | _ {\partial} \right.$ конкондартно стандартному вложению.

Значит, сферу $f' (S^p\times *) $ можно заклеить оснащенным диском
$D ^ {p+1} \subset \partial D^m$ ({\it паутиной}), таким что $\operatorname {Int} D ^ {p+1} \cap\operatorname {Im} f' = \emptyset$ и первые $q$ векторных полей оснащения диска $\partial D ^ {p+1} $ касаются образа $f' (S^p\times S^q) $. По общему положению данная паутина не пересекается с $\operatorname {Im} g$, а значит, и с $\operatorname {Im} g '$.
Ясно, что сфера $ g' (S^p\times *) $ также может  быть заклеена паутиной $\bar D ^ {p+1} \subset \partial D^m$, внутренность которой не пересекается с
$\operatorname {Im} (f '\sqcup g') $ и с первой паутиной $D ^ {p+1} $.

Соединим центры паутин $D ^ {p+1} $ и $\bar D ^ {p+1} $
путем $I$ в $\partial D^m$. По общему положению путь $I$ пересекает $\operatorname {Im} (f '\sqcup g') \cup D ^ {p+1} \cup \bar D ^ {p+1} $ только по своей границе
$\partial I$. Пусть $\bar D^m$ --- объединение малых трубчатых окрестностей многообразий
$D ^ {p+1} $, $I$ и $\bar D ^ {p+1} $ в шаре $D^m$. Ясно, что
пересечение $\operatorname {Im} (f '\sqcup g') \cap \bar D^m$ стандартно. Делая $S^p$-параметрическую связную сумму отображений $f ':S^p\times D^q\to
D^m$ и $ g' :S^p\times D^q\to
D^m$ относительно края в шаре $\bar D ^ {m} $, мы получаем искомое собственное почти вложение $\omega_x:S^p\times D^q\to D^m$.
\end{proof}

Пусть $\overline {E} ^m (S^p\times D^q, S^p\times S ^ {q-1}) $ --- группа
собственных почти вложений $S^p\times D^q\to D^m$ с точностью до собственной
почти конкордантности и связного суммирования с собственными почти вложениями $D ^ {p+q} \to D^m$; структура группы задается операцией параметрической связной суммы относительно края.

\begin{remark} Тем самым относительный $\beta$-инвариант является
отображением $\overline {E} ^m (S^p\times D^q, S^p\times S ^ {q-1}) \to \Omega^m _ {p, q} $.
Теоремы~\ref{t4.1} и~\ref{t5.1} вместе утверждают, что это отображение биективно при $m\ge p +\frac {4} {3} q+2$.
\end{remark}

\noindent 
\beginpicture 
\setcoordinatesystem units <0.850mm,0.850mm>
\setplotarea x from -80 to 70, y from -30 to 30
\put {Рисунок 10} at 0.000 -29.000


\arrow <2mm> [.2,.4] from -43.000 0.000 to -37.000 0.000
\put {$e$} at -40.000 5.000
\arrow <2mm> [.2,.4] from -3.000 0.000 to 3.000 0.000
\put {$h$} at 0.000 5.000
\arrow <2mm> [.2,.4] from 32.000 0.000 to 38.000 0.000
\put {$p$} at 35.000 5.000

\setcoordinatesystem units <0.850mm,0.850mm> point at 60 0
\setplotarea x from -20 to 20, y from -25 to 25


\setsolid
\setplotsymbol ({\fiverm.})
\circulararc 120.000 degrees from -1.121 6.424 center at -3.536 7.071
\circulararc 60.000 degrees from -4.182 9.486 center at -2.241 2.241
\circulararc 120.000 degrees from -9.486 4.183 center at -7.071 3.536
\circulararc 60.000 degrees from -6.424 1.121 center at -8.365 8.365
\circulararc 57.312 degrees from -2.433 -13.738 center at -0.428 -12.776
\circulararc 57.312 degrees from 6.452 -12.370 center at 7.455 -10.385
\circulararc 57.312 degrees from 1.825 -13.832 center at 3.447 -12.310
\circulararc 57.312 degrees from 14.497 -3.852 center at 12.310 -3.447
\circulararc 57.312 degrees from 12.301 -8.584 center at 10.385 -7.455
\circulararc 57.312 degrees from 14.984 0.702 center at 12.776 0.428
\circulararc 360 degrees from 15.000 0.000 center at 0.000 0.000
\circulararc 60.000 degrees from 6.424 -1.121 center at 8.365 -8.365
\circulararc 120.000 degrees from 1.121 -6.424 center at 3.536 -7.071
\circulararc 60.000 degrees from 4.183 -9.486 center at 2.241 -2.241
\circulararc 120.000 degrees from 9.486 -4.183 center at 7.071 -3.536

\put {$S^1 \sqcup S^1 \to S^2$} at 0.000 20.000
\put {$\frac{{E}^2(S^0\times S^1)}{E^2(S^1)}$} at 0.000 -20.000

\setcoordinatesystem units <0.850mm,0.850mm> point at 20 0
\setplotarea x from -20 to 20, y from -25 to 25


\setsolid
\setplotsymbol ({\fiverm.})
\circulararc 57.312 degrees from -2.433 -13.738 center at -0.428 -12.776
\circulararc 57.312 degrees from 6.452 -12.370 center at 7.455 -10.385
\circulararc 57.312 degrees from 1.825 -13.832 center at 3.447 -12.310
\circulararc 57.312 degrees from 14.497 -3.852 center at 12.310 -3.447
\circulararc 57.312 degrees from 12.301 -8.584 center at 10.385 -7.455
\circulararc 57.312 degrees from 14.984 0.702 center at 12.776 0.428
\circulararc 360 degrees from 15.000 0.000 center at 0.000 0.000
\plot -2.500 5.159 -4.183 -1.121 /
\plot 1.121 4.183 -5.314 2.458 /
\circulararc 60.000 degrees from 6.424 -1.121 center at 8.365 -8.365
\circulararc 120.000 degrees from 1.121 -6.424 center at 3.536 -7.071
\circulararc 60.000 degrees from 4.183 -9.486 center at 2.241 -2.241
\circulararc 120.000 degrees from 9.486 -4.183 center at 7.071 -3.536
\circulararc 60.000 degrees from -4.022 6.674 center at -3.052 3.052
\circulararc 120.000 degrees from -6.674 4.022 center at -5.467 3.699
\circulararc 120.000 degrees from -2.492 5.143 center at -3.699 5.467
\circulararc 120.000 degrees from 4.183 1.121 center at 1.768 1.768
\circulararc 60.000 degrees from -1.121 -4.183 center at -3.062 3.062
\circulararc 120.000 degrees from -4.183 -1.121 center at -1.768 -1.768
\setdashes
\circulararc 120.004 degrees from 8.365 2.241 center at 3.536 3.536
\circulararc 59.990 degrees from -2.241 -8.365 center at -6.125 6.125
\circulararc 120.004 degrees from -8.365 -2.242 center at -3.536 -3.536
\circulararc 60.001 degrees from 2.242 8.365 center at 6.124 -6.124
\put {$\bullet$} at -3.062 3.062
\put {$S^1 \sqcup S^1 \to S^2$} at 0.000 20.000
\put {\small $\overline{E}^2(S^0\times S^1)$} at 0.000 -20.000

\setcoordinatesystem units <0.680mm,0.680mm> point at -22 0
\setplotarea x from -20 to 20, y from -30 to 30


\setsolid \setplotsymbol ({\fiverm.}) \circulararc 93.135 degrees
from 7.262 -13.125 center at 0.000 -20.000 \circulararc 70.529
degrees from 5.000 0.000 center at -10.000 0.000 \circulararc
180.000 degrees from -5.000 -0.000 center at -0.000 0.000
\circulararc 70.529 degrees from 5.000 14.142 center at 10.000
0.000 \circulararc 360 degrees from 15.000 0.000 center at 0.000
0.000 \put {$\bullet$} at -5.000 14.142 \put {$\bullet$} at 5.000
14.142 \put {$\bullet$} at -7.262 -13.125 \put {$\bullet$} at
7.262 -13.125 \put {$D^1 \sqcup D^1 \to D^2$} at 0.000 25.000 \put
{\footnotesize $\overline{E}^2(S^0\times D^1,S^0\times S^{0})$} at 0.000 -25.000

\setcoordinatesystem units <0.680mm,0.680mm> point at -67 0
\setplotarea x from -20 to 20, y from -30 to 30


\setsolid
\setplotsymbol ({\fiverm.})
\circulararc 360 degrees from 15.000 0.000 center at 0.000 0.000
\put {$\bullet$} at -5.000 14.142
\put {$\bullet$} at 5.000 14.142
\put {$\bullet$} at -7.262 -13.125
\put {$\bullet$} at 7.262 -13.125
\put {$S^0 \sqcup S^0 \to S^1$} at 0.000 25.000
\put {$\frac{{E}^1(S^0\times S^0)}{E^{1}(S^{0})}$} at 0.000 -25.000

\endpicture

\bigskip

\begin{lemma}\label{t5.2} \textup{(сравни с \cite{Ne84} и \cite[теорема~1]{Sko09}, см.~рисунок~10)}
Для каждого $m> 2p+q+2$ существует точная последовательность
$$
\dfrac {{E} ^m (S^p\times S^q)} {E ^ {m} (S ^ {p+q})}
\xrightarrow{e} \overline {E} ^m (S^p\times S^q) \xrightarrow{h}
\overline {E} ^m (S^p\times D^q, S^p\times S ^ {q-1}) \xrightarrow{p}
\dfrac {{E} ^ {m-1} (S^p\times S ^ {q-1})} {E ^ {m-1} (S ^ {p+q-1})} \to\dots
$$
\end{lemma}

В оставшейся части \S5 мы доказываем лемму~\ref{t5.2}.
Начнем с построения структуры коммутативной группы
на множестве заузленных торов
(а также почти вложений $S^p\times S^q\to S^m$).
Данное построение эквивалентно приведенному в \cite[\S2]{Sko08}.

\begin{definition}[паутины $D ^ {p+1} $]
Зафиксируем отмеченную точку $*\in S^q$.
Не ограничивая общности, будем считать, что $*\times S^q\cap B ^ {p+q} = \emptyset$.
{\it Паутина} почти вложения $f:S^p\times S^q\to S^m$ --- это оснащенный диск $D ^ {p+1} \subset S^m$, удовлетворяющий следующим $3$ условиям:

\smallskip

\noindent  (i) $\partial D ^ {p+1} =f (S^p\times *) $;

\noindent  (ii) $\operatorname {Int}
D ^ {p+1} \cap\operatorname {Im} f =\emptyset$;

\noindent  (iii) первые $q$ векторных полей оснащения края $\partial D ^ {p+1} $ совпадают с очевидным оснащением ``меридиана'' $f(S^p\times *)$ в ``торе'' $f(S^p\times S^q)$.

\smallskip
\noindent  {\it Паутины}
почти конкордантности $f:S^p\times S^q\times I\to S^m\times I$ и
собственного почти вложения $f:S^p\times D^q\to D^m$
определяются аналогично (в качестве отмеченной точки $*\in D^q$ берется центр диска  $D^q$).
\end{definition}

Следующее утверждение эквивалентно лемме о стандартизации~2.1 из \cite{Sko08}.

\begin{propos}\label{p5.3}
\cite[лемма о стандартизации 2.1]{Sko08} Если $m> 2p+q+1$, то для любого
почти вложения $f:S^p\times S^q\to S^m$ существует паутина.
Если $m> 2p+q+2$, то для любой почти конкордантности между почти вложениями $f_1, f_2:S^p\times S^q\to
S^m$ существует паутина, продолжающая данные паутины почти вложений $f_1$ и $f_2$.
\end{propos}

\begin{proof} \cite{Sko08}
Расслоение $N (f | _ {S^p\times *}; S^m) $
стабильно тривиально и $m-p-q\ge p$, следовательно, это расслоение тривиально.
Возьмем $(m-p-q)$-оснащение $\xi$ этого расслоения.

Рассмотрим сечение, образованное первыми векторами оснащения $\xi$. Так как $m\ge2p+q+2\ge2p+2$, то вложение $f | _ {S^p\times *} $ незаузлено
в $S^m$. Таким образом, существует вложение $\bar f:D ^ {p+1} \subset S^m$,
удовлетворяющее условию (i) из определения паутины. Так как $m\ge2p+q+2$, то по общему положению мы можем
добиться также выполнения условия (ii).

Выбрасывая первый вектор из оснащения $\xi$, мы получаем некоторое $(m-p-q-1)$-оснащение $\xi_1$ края $\bar f(\partial D^{p+1})$, ортогональное диску $\bar f(D^{p+1})$. Обозначим через $\eta$ очевидное нормальное $q$-оснащение сферы $f (S^p\times *) $ в торе $f (S^p\times S^q) $. Тогда $ (\xi_1, \eta) $ будет
$(m-p-1)$-оснащением края $\bar f(\partial D^{p+1})$, ортогональным диску $\bar f(D^{p+1})$. Так как  $p <m-p-q-1$, то отображение
$\pi_p (SO_ {m-p-q-1}) \to\pi_p (SO _ {m-p-1}) $ сюръективно. Следовательно, мы
можем изменить оснащение $\xi_1$ (и вместе с ним оснащение $\xi$) так, чтобы оснащение $ (\xi_1, \eta) $ продолжалось до оснащения всего диска $\bar f (D ^ {p+1}) $. Ясно, что полученное оснащение будет удовлетворять условию (iii).

Второе утверждение предложения доказывается аналогично, см.~детали в \cite[Доказательство леммы о стандартизации 2.1 в \S3]{Sko08}.
\end{proof}

Следующее определение эквивалентно приведенному в \cite{Sko08}.

\begin{definition}[параметрической связной суммы] (см.~рисунок~5)
Пусть $f_1, f_2:S^p\times S^q\to S^m$
--- пара почти вложений. Не теряя общности, можно считать, что $\operatorname {Im} f_1\subset D^m _ + $, $\operatorname {Im} f_2\subset D^m_-$.
Возьмем паутины $D ^ {p+1} \subset D^m _ + $ и $\bar D ^ {p+1} \subset D^m_-$ этих почти вложений.
Соединим центры паутин $D ^ {p+1} $ и $\bar D ^ {p+1} $
путем $I$ в $S^m$. По общему положению путь $I$ пересекает
$\operatorname {Im} f_1\cup\operatorname {Im} f_2\cup D ^ {p+1} \cup \bar D ^ {p+1} $ только по своей границе
$\partial I$. Пусть $\bar D^m$ --- объединение небольших трубчатых окрестностей многообразий $D ^ {p+1} $, $I$ и $\bar D ^ {p+1} $ в $S^m$. Ясно, что
пересечение $(\operatorname {Im} f_1 \cup\operatorname {Im} f_2) \cap \bar D^m$
стандартно. {\it Параметрической связной суммой} почти вложений $f_1$ и $f_2$ мы назовем почти вложение $f_1+f_2:S^p\times S^q\to S^m$,
полученное с помощью $S^p$-параметрической связной суммы почти вложений $f_1$ и $f_2$ внутри шара $\bar D^m$.
\end{definition}

\begin{definition}[противоположного элемента] Пусть $f:S^p\times S^q\to S^m$ ---
почти вложение. {\it Противоположным элементом} для $f$ мы называем почти вложение, заданное формулой $ (-f) (x, y) = \sigma_mf (x, \sigma_qy) $,
где $\sigma_k$ --- симметрия сферы $S^k$ относительно гиперплоскости $x_1=0$.
\end{definition}

\begin{definition}[нейтрального элемента]
{\it Нейтральным элементом} мы назовем {\it стандартное вложение}
$S^q\times S^p\to D ^ {q+1} \times D ^ {p+1}\cong D^{p+q+2} \subset D^m\subset S^m$.
\end{definition}

Следующий важный результат доказывается прямой проверкой (см.~детали доказательства в \cite{Sko08}).

\begin{theorem}[о групповой структуре]\label{t5.4}
\cite{Sko08}
Пусть $m> 2p+q+2$. Тогда приведенное выше построение корректно определяет структуру коммутативной группы на множестве почти вложений $S^p\times S^q\to S^m$ с точностью до почти конкордантности.
\end{theorem}

Аналоги этой теоремы для собственных почти вложений и гладких вложений также справедливы и доказываются аналогично.

\begin{propos}\label{p5.5}
При $m>2p+q+2$ относительный $\beta$-инвариант $\beta:\overline {E} ^m (S^p\times D^q, S^p\times S ^ {q-1}) \to \Omega^m_{p, q}$ является гомоморфизмом.
\end{propos}

\begin{proof} Рассмотрим тройку
собственных почти вложений $f_1$, $f_2$ и $f_3=f_1+f_2$, где ``$+$'' обозначает $S^p$-параметрическую связную сумму
относительно края. Легко видеть, что с точностью до гомотопии $\tilde f_3 =\tilde f_1\vee\tilde f_2$ и
$\bar f_3 =\bar f_1 \natural \bar f_2$, где ``$\natural$'' обозначает связную сумму относительно края. Значит, $\beta (f_3) = \beta (f_1) \sqcup\beta (f_2) $, и, следовательно,
$\beta (f_1+f_2) = \beta (f_1) + \beta (f_2) $.
\end{proof}

Теперь мы собираемся доказать, что действие группы вложений $S ^ {p+q} \to S^m$ на множестве вложений $S^p\times S^q\to S^m$ инъективно.


Введем обозначения:

\noindent  (a) $E^m (S ^ {p+q}) $ --- группа всех вложений $S ^ {p+q} \to S^m$ с точностью до конкордантности;

\noindent  (b) $E^m (S^p\times S^q) $ --- группа всех вложений $S^p\times S^q\to S^m$ с точностью до конкордантности;
структура группы задается операцией параметрической связной суммы.

\noindent  (c) $\kappa^*:E ^ {m} (S ^ {p+q}) \to E ^ {m} (S^p\times S^q) $ --- отображение, сопоставляющее вложению $g:S ^ {p+q} \to S^m$ связную сумму этого вложения $g$ и стандартного вложения $S^p\times S^q\to S^m$. (Связное суммирование производится вдоль некоторого пути $I$, соединяющего образы этих вложений; предполагается, что образы разделены гиперплоскостью).

\begin{propos}\label{p5.6}
Отображение $\kappa^*:E ^ {m} (S ^ {p+q}) \to E ^ {m} (S^p\times S^q) $ инъективно при $m> 2p+q+2$.
\end{propos}

Данное предложение немедленно влечет случай ``$p+q+1$ делится на $4$'' в теореме~\ref{tI}, согласно теореме~\ref{t0} из \S1.

\begin{proof}[предложения \ref{p5.6}]
Достаточно построить левое обратное
$\bar\kappa^*:E ^ {m} (S^p\times S^q) \to E ^ {m} (S ^ {p+q}) $ для отображения $\kappa^*$.

Отображение $\bar\kappa^*:E ^ {m} (S^p\times S^q) \to E ^ {m} (S ^ {p+q}) $ определяется следующим образом.
Возьмем вложение $f:S^p\times S^q\to S^m$. Согласно предложению~\ref{p5.3} у него есть паутина $D ^ {p+1} $.
Произведем вложенную хирургию многообразия $S^p\times S^q$ вдоль оснащенного диска $D ^ {p+1} $.
Пусть $\bar\kappa^* (f) $ ---  изотопический класс вложения $S ^ {p+q} \to S^m$, полученного с помощью этой хирургии.

Элемент $\bar\kappa^* (f) $ корректно определен согласно второму утверждению предложения~\ref{p5.3}.
Мы имеем $\bar\kappa^*\kappa^* =\mathop {id} $,
потому что
$\bar\kappa^*\kappa^* (f) = \bar\kappa^*(f \# s) =f \#\bar\kappa^* (s) =f \# 0=f$ для любого $f\in E^m (S ^ {p+q}) $.
\end{proof}

\begin{propos}\label{p5.7} \textup{(a)} Для каждого $m> 2p+q+2$ все собственные вложения $S^p\times D^q\to D^m$
собственно конкордантны друг другу.

\textup{(b)} \cite[критерий тривиальности]{Sko08} При $m> 2p+q+2$ вложение $S^p\times S ^ {q-1} \to S ^ {m-1} $
конкордантно стандартному вложению, если и только если оно продолжается до собственного вложения $S^p\times D^q\to D^m$.
\end{propos}

\begin{proof} (a) Возьмем собственное вложение $f:S^p\times D^q\to D^m$. Аналогично
предложению~\ref{p5.3} можно показать, что у него есть паутина $D ^ {p+1} $. Пусть $\bar D^m$ --- трубчатая окрестность паутины $D ^ {p+1} $. Ясно, что ограничение $f:f ^ {-1} \bar D^m\to \bar D^m$
конкордантно стандартному вложению $S^p\times D^q\to D^m$.
Остается доказать, что $f$ конкордантно этому ограничению.
Пусть $c:S^p\times D^q\times I\to D^m\times I$ ---
неподвижная конкордантность, совпадающая с вложением $f:S^p\times D^q\to D^m$ при каждом $t\in I$.
Пусть $h_m:D^m\times I\to D^m\times I$ --- диффеоморфизм, неподвижный на $D^m\times 0 $
и переводящий $\bar D^m\times 1 $ в $D^m\times 1 $. Пусть $h _ {p+q}:S^p\times D^q\times I\to S^p\times D^q\times I$ --- диффеоморфизм, неподвижный на $S^p\times S^q\times 0 $ и переводящий $S^p\times D^q\times 1 $ в $f ^ {-1} \bar D^m$.
Тогда композиция
$h_m c h _ {p+q} $ является собственной конкордантностью между $f:S^p\times D^q\to D^m$ и ограничением $f:f ^ {-1} \bar D^m\to \bar D^m$.

(b) следует непосредственно из (a).
\end{proof}

Аналог данного предложения для собственных почти вложений также справедлив, с той разницей, что для этого шар $B ^ {p+q} $ в определении собственного почти вложения нужно заменить на шар $\bar B ^ {p+q} \subset S^p\times D^q$, пересекающий край по своей грани. {\it Грань} шара $\bar B ^ {p+q} $ --- это шар, содержащийся в $\partial \bar B ^ {p+q} $.

\begin{proof}[леммы \ref{t5.2}] (сравни c
\cite[доказательство теоремы 3.1]{Sko09}) (1) {\it
Построение гомоморфизмов.} Пусть $e$ --- очевидное отображение.
Пусть $p$ --- отображение ограничения на край. Гомоморфизм $h$ --- это ``гомоморфизм вырезания'', определяемый следующим образом. Возьмем почти вложение
$f:S^p\times S^q\to S^m$. Согласно предложению~\ref{p5.3} у него есть паутина $D ^ {p+1} \subset S^m$. Пусть $\bar D^m$ --- трубчатая окрестность паутины
$D ^ {p+1} $. Положим $h (f) $ равным ограничению $f:
(S^p\times S^q-f ^ {-1} \operatorname {Int} \bar D^m) \to
S^m-\operatorname {Int} \bar D^m$.

(2) {\it Точность в члене $E^m (S^p\times S^q)/E^m (S ^ {p+q}) $.}
Аналогично предложению~\ref{p5.7}~(b) можно показать, что
вложение $f:S^p\times S^q\to S^m$ продолжается
до собственного почти вложения $S ^ {p} \times D ^ {q+1} \to D ^ {m+1} $, если и
только если $f:S^p\times S^q\to S^m$ почти конкордантно стандартному вложению.
Это означает, что наша последовательность точна в члене
$E^m (S^p\times S^q)/E^m (S ^ {p+q}) $.

(3) {\it Точность в члене $\overline {E} ^m (S^p\times D^q, S^p\times S ^ {q-1}) $.} Включение $\operatorname {Im} h\subset\operatorname {ker} p$
следует из предложения~\ref{p5.7} (b).
Чтобы доказать $\operatorname {ker} p\subset\operatorname {Im} h$,
возьмем собственное почти вложение $f:S^p\times D^q\to D^m$, такое что $p (f) =0$. Пусть $f\left |_\partial\right.:S^p\times \partial D^q\to \partial D^m$
--- ограничение отображения $f$ на край.
По определению, существует гладкое вложение $g:S ^ {p+q-1} \to S ^ {m-1} $, такое что связная сумма $f\left |_\partial\right. + g$ конкордантна стандартному вложению. Продолжим $g:S ^ {p+q-1} \to S ^ {m-1} $
до собственного почти вложения  $g ':D ^ {p+q} \to D^m$.
Пусть $f+g '$ --- связная сумма $f$ и $g' $ относительно края. По предложению~\ref{p5.7} (b)
отображение $f+g ':S^p\times D^q\to D^m$ продолжается до почти вложения $f ':S^p\times S^q\to S^m$. Таким образом, $f=h (f') $.

(4) {\it Точность в члене $\overline {E} ^m (S^p\times S^q) $.} Включение
$\operatorname {Im} e\subset \operatorname {ker} h$ следует из предложения~\ref{p5.7} (a). Чтобы доказать $\operatorname {ker} h\subset\operatorname {Im} e$, возьмем почти вложение $f:S^p\times S^q\to S^m$,
такое что $h (f) =0$. Тогда
существует собственная почти конкордантность $c$ между $h (f) $ и связной суммой
стандартного вложения $S^p\times D^q\to D^m$ с некоторым собственным почти вложением $g:D ^ {p+q} \to D^m$.
Так как ограничение почти конкордантности $c$ на край является гладкой конкордантностью, то по предложению~\ref{p5.6}
ограничение почти вложения $g$ на край незаузленно.
Значит, можно считать, что $g$ --- гладкое
вложение. 
По теореме о продолжении конкордантности
ограничение конкордантности $c$ на край продолжается до объемлемой конкордантности шара $S^m-D^m$. Тем самым $c$ может быть продолжена до почти конкордантности
почти вложения $f$ без добавления новых самопересечений. Последняя почти конкордантность является почти конкордантностью между почти вложением $f$ и некоторым вложением $f ':S^p\times S^q\to S^m$.
Следовательно, $f=e (f') $.
\end{proof}

Доказательство теоремы~\ref{tII} завершено.

\section{Критерий конечности}\label{s6}

В данном разделе мы выводим теорему~\ref{tI} из теоремы~\ref{tII}. Для этого нам потребуется классификация почти вложений $S^p\times S^q\to S^m$. Суммируем данную классификацию и полученные выше результаты в следующей теореме.

Пусть
$E^m (D^p\times S^q) $ обозначает группу гладких
вложений $D^p\times S^q\to S^m$ с точностью до конкордантности; структура группы задается операцией $D^p$-параметрической связной суммы.

\begin{theorem}\label{t6.1}
При $m> p +\frac {4} {3} q+2$ и $m> 2p+q+2$ существуют
следующие точные последовательности

\smallskip
\noindent  $\begin{matrix}
(1) &0 &\to &E^m (S ^ {p+q}) &\to & {E} ^m (S^p\times S^q) &\to &\frac {{E} ^m (S^p\times S^q)} {E ^ {m} (S ^ {p+q})}
\\
(2) &\overline {E} ^ {m+1} (S^p\times S ^ {q+1}) &\to &\Omega ^ {m+1} _ {p, q+1} &\to &\frac {{E} ^m (S^p\times S^q)} {E ^ {m} (S ^ {p+q})} &\to
&\overline {E} ^m (S^p\times S^q)
\\
(3) &E ^ {m+1} (D^p\times S ^ {q+1}) &\to&\pi _ {p+q} (S ^ {m-q-1}) &\to &\overline {E} ^m (S^p\times S^q) &\to
&E^m (D^p\times S^q)
\\
(4) &E ^ {m+1} (S ^ {q+1}) &\to&\pi _ {q} (V _ {m-q, p}) &\to &E^m (D^p\times S^q) &\to &E^m (S ^ {q}).
\end{matrix}
$
\end{theorem}


Группы во втором столбце теоремы~\ref{t6.1}
известны {\it рационально}:

\begin{theorem}\label{t6.2}
Пусть $p +\frac {4} {3} q+2\le m <
p +\frac {3} {2} q+2$, $m> 2p+q+2$ и $m> n+2$. Тогда

\smallskip
\noindent  \textup{(1)}
$E ^ {m} (S^n) $ бесконечна, если и только если $m\le\frac {3} {2} n +\frac {3} {2} $ и
$n+1$ делится на $4$.

\noindent  \textup{(2)}  $\Omega _ {p, q+1} ^ {m+1} $ бесконечна, если и только если
$m=p +\frac {3} {2} q +\frac {3} {2} $ и $q+1$ делится на $4$.

\noindent  \textup{(3)}  $\pi _ {p+q} (S ^ {m-q-1}) $ бесконечна, если и только если
$m =\frac {1} {2} p +\frac {3} {2} q +\frac {3} {2} $ и $p+q+1$ делится на~$4$.

\noindent  \textup{(4)} $\pi_q (V _ {m-q, p}) $ бесконечна, если и только если
$p\ge1$, $\frac {3} {2} q +\frac {3} {2} \le m\le p +\frac {3} {2} q +\frac {1} {2} $
и $q+1$ делится на $4$.
\end{theorem}

План наших дальнейших рассуждений таков.
Сначала мы докажем теорему~\ref{tI}, считая теоремы~\ref{t6.1} и~\ref{t6.2} доказанными. Потом мы докажем сами теоремы~\ref{t6.1} и~\ref{t6.2}, опираясь на некоторые известные результаты.

\begin{proof}[теоремы~\ref{tI}, исходя из теорем~ \ref{t6.1} и \ref{t6.2}] (1) {\it Случай, когда $q+1, p+q+1$ не делятся на $4$.}
Напомним, что если $X\to Y\to Z$ --- точная
последовательность с конечными группами $X$ и $Z$, то группа $Y$ также конечна.
Применим это наблюдение $4$ раза к последним $3$ столбцам таблицы в теореме~\ref{t6.1},
начиная снизу. Мы получаем требуемый результат, так как по теореме~\ref{t6.2} группы во втором столбце теоремы~\ref{t6.1} конечны, когда $q+1, p+q+1$ не делятся на $4$.

\smallskip

(2) {\it Случай, когда $p+q+1$ делится на $4$.} По теореме~\ref{t6.1} (1) и~\ref{t6.2} (1) непосредственно получаем, что
группа $E^m (S^p\times S^q) $ в этом случае бесконечна.

\smallskip

(3) {\it Случай, когда $q+1$ делится на $4$, $m\le
\frac {3} {2} q +\frac {3} {2} $.}
По теореме~\ref{t6.2} (1)
группа $E^m (S^q) $ в этом случае бесконечна. Возьмем ее
элемент $x$ бесконечного порядка. Препятствие к существованию $ (p+1) $-оснащения
вложения $x:S ^ {q} \to S^m$ лежит в группе
$\pi _ {q-1} (V _ {m-q, p+1}) $. По теореме~\ref{t6.2} (4)
эта группа в нашем случае конечна. Значит,
для некоторого целого положительного числа $N$ вложение $Nx$ продолжается до
гладкого вложения $H:S^p\times S^q\to S^m$. Так как ограничение этого вложения на сферу $*\times S^q$ имеет бесконечный порядок, то и само вложение $H$ имеет бесконечный порядок.

\smallskip
(4) {\it Случай, когда $q+1$ делится на $4$, $\frac {3} {2} q +\frac {3} {2} <m\le p +\frac {3} {2} q +\frac {1} {2} $, $p\ge 1 $.} Достаточно построить вложение $T:S^p\times S^q\to S^m$, имеющее бесконечный порядок в группе
$E^m (S^p\times S^q) $.

{\it Построение вложения $T$.} По теореме~\ref{t6.2} (4) группа
$\pi_q (V _ {m-q, p}) $ в этом случае бесконечна. Возьмем элемент $x$ бесконечного порядка в этой группе. Рассмотрим отображение
$\tau:\pi_q (V _ {m-q, p}) \to E^m (D^p\times S^q) $ из теоремы~\ref{t6.1} (4). Это отображение
переводит элемент $x$ в каноническое $p$-оснащение $D^p\times S^q\to S^m$ стандартной сферы $S^q\subset S^m$. Полное
препятствие к продолжению этого $p$-оснащения до $ (p+1) $-оснащения
лежит в группе $\pi _ {q-1} (S ^ {m-p-q-1}) $. Последняя группа конечна в нашем случае. Значит, для некоторого целого положительного числа $N$ элемент $N\tau x $
продолжается до искомого гладкого вложения $T:S^p\times S^q\to S^m$.

{\it Доказательство того, что вложение $T$ имеет бесконечный порядок.} Достаточно доказать, что элемент $\tau (x) \in E^m (D^p\times S^q) $, который является
ограничением вложения $T$ на $D^p\times S^q$, имеет бесконечный порядок.
Предположим противное. Тогда $N\tau (x) =0$ для некоторого целого положительного числа $N$. Значит, по
теореме~\ref{t6.1} (4) $Nx$ принадлежит образу отображения
$E ^ {m+1} (S ^ {q+1}) \to\pi _ {q} (V _ {m-q, p}) $. Но группа $E ^ {m+1} (S ^ {q+1}) $ в нашем случае конечна. Значит, $x$ имеет конечный порядок, вопреки нашему выбору. Данное противоречие показывает, что $T$ имеет
бесконечный порядок.

\smallskip
(5) {\it Случай, когда $q+1$ делится на $4$, $m=p +\frac {3} {2} q +\frac {3} {2} $.}

{\it Построение вложения $W$.} По теореме~\ref{t6.2} (2)
группа $\Omega ^ {m+1} _ {p, q+1} $ в этом случае бесконечна. Возьмем в ней элемент $x$ бесконечного порядка. Пусть $W:S^p\times S^q\to S^m$ ---
вложение, реализующее образ $x$ при отображении
$\Omega ^ {m+1} _ {p, q+1} \to {{E} ^m (S^p\times S^q)} / {E ^ {m} (S ^ {p+q})} $
из теоремы~\ref{t6.1} (2).

{\it Доказательство того, что $W$ имеет бесконечный порядок.} Рассмотрим точную последовательность (2) из теоремы~\ref{t6.1}.  Достаточно доказать, что $\overline {E} ^ {m+1} (S^p\times S ^ {q+1}) $ конечна в нашем случае. Так как $q+1$ делится на $4$, то по теореме~\ref{t6.2} группы $\pi _ {q+1} (V _ {m-q, p}) $ и $E ^ {m+1} (S ^ {q+1}) $ конечны. Как показывает несложное вычисление, $\pi _ {p+q+1} (S ^ {m-q-1}) $ также конечна в нашем случае.
По теореме~\ref{t6.1} (3)--(4) получаем, что $\overline {E} ^ {m+1} (S^p\times S ^ {q+1}) $ конечна.
\end{proof}

В оставшейся части раздела мы докажем теоремы~\ref{t6.1} и~\ref{t6.2}.

Утверждения~(1) и~(2) теоремы~\ref{t6.1} являются переформулировками предложений~\ref{p5.6} и теоремы~\ref{tII}, которые были доказаны в ~\S5.
Теорема~\ref{t6.1} ~ (4) доказывается непосредственно, аналогично \cite[следствие~5.9] {Hae66A}. Для удобства читателя мы приводим набросок доказательства ниже.
Теорема~\ref{t6.1} (3) доказана ранее в \cite[лемма об ограничении~5.2]{Sko08}
при $p\ge 1 $, $m> 2p+q+2$ и $m\ge \frac {1} {2} p +\frac {3} {2} q+2$. Дополнительное ограничение
$m\ge \frac {1} {2} p +\frac {3} {2} q+2$ в доказательстве этого утверждения в \cite{Sko08}
не используется. Мы приводим набросок альтернативного доказательства для удобства читателя.

\begin{proof}[утверждения (3) в теореме~\ref{t6.1} (набросок)] (a) {\it Определение
групп $\tilde E^m (S^p\times S^q) $ и $\tilde {E} ^m (S^p\times D^q, S ^ {p} \times S ^ {q-1}) $.} Отображение
$f:S^p\times S^q\to S^m$, назовем {\it слабым почти
вложением}, если оно является вложением вне фиксированного шара $B ^ {p+q} \subset S^p\times S^q$ (пересечению $fB ^ {p+q} \cap f (S^p\times S^q-B ^ {p+q}) $ разрешается быть непустым). {\it Слабая почти конкордантность} определяется аналогично. Обозначим $\tilde E^m (S^p\times S^q) $ группу слабых почти
вложений с точностью до слабой почти конкордантности.
Отождествим группы $\tilde E^m (S^p\times S^q) $ и $E^m (D^p\times S^q) $. Ясно, что эти группы изоморфны при $m> 2p+q+2$.

Зафиксируем шар $\bar B ^ {p+q} \subset (S^p-*)\times D^q$, пересекающий край по своей грани.
Собственное отображение $f:S^p\times D^q\to D^m$ назовем
{\it слабым собственным почти вложением}, если выполнены следующие два условия:

\noindent  (i) $f$ является вложением вне шара $\bar B ^ {p+q} $; и

\noindent  (ii) $f (S^p\times\partial D^q\cap \bar B ^ {p+q}) \cap
f (S^p\times\partial D^q - \bar B ^ {p+q}) = \emptyset$.

\noindent  {\it Слабая собственная почти конкордантность} определяется аналогично. Обозначим через $\tilde {E} ^m (S^p\times D^q, S ^ {p} \times S ^ {q-1}) $
группу слабых собственных почти вложений с точностью до слабой собственной почти конкордантности.

(b) {\it При каждом $m> 2p+q+2$ существуют точная последовательность:}
$$
\overline {E} ^m (S ^ {p} \times S ^ {q})
\xrightarrow{e}
\tilde {E} ^m (S^p\times S^q) \xrightarrow{h}
\tilde {E} ^m (S^p\times D^q, S ^ {p} \times S ^ {q-1}) \xrightarrow{p}
\overline {E} ^ {m-1} (S ^ {p} \times S ^ {q-1}) \to\dots
$$
Здесь $e$, $h$ и $p$ --- очевидно определяемые гомоморфизмы забывания, вырезания и ограничения, соответственно. Данное утверждение доказывается
полностью аналогично лемме~\ref{t5.2}.

(c) {\it Определение гомоморфизма
$\lambda:\tilde {E} ^m (S^p\times D^q, S ^ {p} \times S ^ {q-1}) \to\pi _ {p+q-1} (S ^ {m-q-1}) $.}
Возьмем слабое собственное почти вложение $f:S^p\times D^q\to D^m$. По определению получаем
$f\partial \bar B ^ {p+q} \cap f (*\times D^q) = \emptyset$. Заметим, что
$D^m-f (*\times D^q) \simeq S ^ {m-q-1} $. Положим $\lambda (f) $ равным
гомотопическому классу ограничения $f:\partial \bar B ^ {p+q} \to D^m-f (*\times D^q) $.

(d) {\it $\lambda$ инъективно.} Возьмем слабое собственное почти вложение $f:S^p\times D^q\to D^m$, такое что $\lambda (f) =0$. Тогда
$f\left | _ {\partial \bar B ^ {p+q}} \right.$ продолжается до отображения $g:\bar B ^ {p+q} \to D^m$, не пересекающего $f (*\times D^q) $. Так как $f\left | _ {S^p\times D^q-\bar B ^ {p+q}} \right.$ ---
вложение, то
пересечение $g\bar B ^ {p+q} $ с $f (S^p\times D^q-\bar B ^ {p+q}) $ может быть устранено с помощью подходящей
гомотопии относительно края.
Значит, мы можем считать, что $g$ не пересекает
$f (S^p\times D^q-\bar B ^ {p+q})$. Произведем слабую собственную почти конкордантность, заменяющую
$f\left | _ {\bar B ^ {p+q}} \right.$ на $g$. Аналогично предложению~\ref{p5.7} (a) можно показать, что
полученное отображение слабо собственно почти конкордантно стандартному вложению $S^p\times D^q\to D^m$.

(e) {\it $\lambda$ сюръективно.} Возьмем элемент $x\in
\pi _ {p+q-1} (S ^ {m-q-1}) $. Возьмем стандартное вложение
$f:S ^ {p} \times D^q\to D^m$. Реализуем элемент $x$ отображением
$g:S ^ {p+q-1} \to \partial D^m-f (*\times D^q) $.
Так как $f\left | _ {S^p\times D^q-\bar B ^ {p+q}} \right.$ ---
вложение, то
пересечение $g\bar B ^ {p+q} $ с $f (S^p\times D^q-\bar B ^ {p+q}) $ может быть устранено с помощью подходящей
гомотопии относительно края.
Таким образом, мы можем считать, что $g$ не пересекает
$f (S^p\times D^q-\bar B ^ {p+q}) $.
Продолжим отображение $g:S ^ {p+q-1} \to \partial D^m$ до собственного отображения $g':D ^ {p+q} \to D^m$. Пусть $\mu_x$ --- связная сумма
(относительно края) отображений $ g'$ и $f$.
Ясно, что $\lambda (\mu_x) =x$. Это заканчивает доказательство утверждения ~ (3).
\end{proof}

\begin{proof}[утверждения (4) в теореме~\ref{t6.1} (набросок)]
(a) {\it Определение гомоморфизмов.} Отображение $i ^*: {E} ^m (D^p\times S^q) \to E ^ {m} (S^q) $ индуцировано ограничением. Здесь произведение $0\times S^q$ отождествлено с $S^q$ очевидным способом.

Отображение $Ob:E ^ {m} (S^q) \to \pi _ {q-1} (V _ {m-q, p}) $ --- это полное препятствие к существованию $p$-оснащения вложения $S^q\to S^m$.
Напомним определение этого препятствия. Возьмем вложение $f:S^q\to S^m$.
Возьмем (единственное с точностью до гомотопии) $p$-оснащение диска $fD^q_-$. Возьмем (единственное с точностью до гомотопии) $(m-q)$-оснащение диска $fD^q_+$. Тогда сфера $fS^{q-1}$ будет снабжена как $p$-оснащением,
так и $(m-q)$-оснащением. Отождествим каждый слой ограничения нормального расслоения к $fS^q$ на сферу $fS^{q-1}$ с пространством $\mathbb{R}^{m-q}$, пользуясь данным $(m-q)$-оснащением. Тем самым $p$-оснащение определит некоторое отображение $S^{q-1}\to V _ {m-q, p}$. Гомотопический класс этого отображения и будет искомым препятствием $Ob(f)$.

Отображение $\tau:\pi _ {q} (V _ {m-q, p}) \to {E} ^m (D^p\times S^q) $ определяется следующим образом.
Реализуем $f\in \pi _ {q} (V _ {m-q, p}) $ как гладкое отображение $f:D^p\times S^q\to D ^ {m-q} $, линейное на каждом слое $D^p\times *$.
Определим $\tau (f) $ как композицию $D^p\times S^q\to D ^ {m-q} \times S^q\to S^m$ вложения $f\times \operatorname {pr} _2$ и стандартного
вложения $s$.

(b){\it Точность в членах $ {E} ^m (D^p\times S^q) $ и $E ^ {m} (S^q) $} проверяется непосредственно.

(c){\it Точность в члене $\pi _ {q} (V _ {m-q, p}) $.} Пусть $f:S ^ {q+1} \to S ^ {m+1} $ --- вложение. Тогда $f$ изотопно {\it стандартизованному} вложению $f ':S ^ {q+1} \to S^m$, то есть такому вложению, что
\begin{itemize}
\item $f ':D ^ {q+1} _-\to D ^ {m+1} _-$ является ограничением стандартного вложения $S ^ {q+1} \to S ^ {m+1} $;
\item $f' (D ^ {q+1}) _ + \subset D ^ {m+1} _ +$.
\end{itemize}
Возьмем $p$-оснащение диска $f' (D ^ {q+1}) $. Ясно, что вложение $\tau Ob f ':D^p\times S^q\to S^m$ продолжается до вложения $D^p\times D ^ {q+1} \to D ^ {m+1} $, заданного $p$-оснащением. Поэтому $\tau Ob f' $ изотопно стандартному вложению $D^p\times S^q\to S^m$. Таким образом, $\operatorname{Im}\tau\subset \operatorname{ker} Ob$. Аналогично, $\operatorname{Im}\tau\supset \operatorname{ker} Ob$.
\end{proof}

Теорема~\ref{t6.2} легко сводится к известным результатам. Утверждение (1) доказано в
\cite[следствие~6.7]{Hae66A}. Утверждение (3) следует из теоремы Серра.

\begin{proof}[утверждения (4) в теореме~\ref{t6.2}]
Из предположений $m> 2p+q+2$ и $m <p +\frac {3} {2} q+2$
следует, что $m\le 2q$. Докажем требуемое утверждение, в котором данная пара ограничений заменена на единственное ограничение $m\le 2q$,
индукцией по $p$.

{\it (a) Случай, когда $q+1$ не делится на 4.} Так как $m\le 2q$, то группа $\pi_q (V _ {m-q, 1}) \cong \pi_q (S ^ {m-q-1}) $ конечна.
Пользуясь точной гомотопической последовательностью расслоения ``ограничения'' $S ^ {m-p-q} \to V _ {m-q, p} \to
V _ {m-q, p-1} $, тензорно помноженной на $\Bbb Q$, мы получаем по индукции, что группа $\pi_q (V _ {m-q, p}) $ конечна.

{\it (b) Случай, когда $q+1$ делится на 4 и выполнено одно из двух условий: $ m<\frac {3} {2} q +\frac {3} {2} $ или
$m> p +\frac {3} {2} q +\frac {1} {2} $.} В этом случае группы $\pi_q (S ^ {m-q-i}) $ по-прежнему
конечны для каждого $i=1,2, \dots, p$. Аналогично предыдущему
мы получаем, что группа $\pi_q (V _ {m-q, p})$ конечна.

{\it (c) Случай, когда $q+1$ делится на 4 и
$\frac {3} {2} q +\frac {3} {2} \le m\le p +\frac {3} {2} q +\frac {1} {2} $.}
Возьмем такое число $i$, что $m=i +\frac {3} {2} q +\frac {1} {2} $. Рассмотрим указанную выше точную гомотопическую последовательность при $p=i$. Аналогично предыдущему
можно показать, что если $q+1$ делится на $4$ и $m\le 2q$, то группа $\pi _ {q+1} (V _ {m-q, i-1}) $ конечна. Таким образом, группа
$\pi_q (V _ {m-q, i}) $ бесконечна. По индукции получаем, что группа $\pi_q (V _ {m-q, p}) $ также бесконечна.
\end{proof}

\begin{proof}[утверждения (2) в теореме~\ref{t6.2}]
Будем использовать обозначения $s=2p+3q-2m+2$ и $l=m-p-q-1$. Тогда рассматриваемая группа будет изоморфна $\pi_{s+l+1}(V_{M+l,M})$. Наше ограничение $p+\frac{4}{3}q+2\le m < p+\frac{3}{2}q+2$ эквивалентно ограничению $-1\le s\le l-3$.

(a) {\it Случай  $s =-1$}. Согласно таблицам
в~\cite{Pae56} группа $\pi _ {l} (V _ {M+l, M}) $
бесконечна, если и только если $l$ делится на $2$.
Вместе с условием $s =-1$ это эквивалентно условиям
$m=p +\frac {3} {2} q +\frac {3} {2} $, $q+1$ делится на $4$.

(b) {\it  Случай $0\le s\le l-3$}.
Докажем индукцией по $s$, что
группа $\pi _ {s+l+1} (V _ {M+l, M}) $ конечна. База $s=0$ следует из таблиц в~\cite{Pae56}. При $s> 0 $
рассмотрим точную гомотопическую последовательность расслоения ``ограничения''
$ S ^ {l} \to V _ {M+l, M} \to V _ {M+l, M-1} $, тензорно помноженную на $\mathbb{Q}$. В этой последовательности $\pi _ {s+l+1} (S^l) $ конечна, потому что $0\le s\le l-3$. Согласно предположению индукции, группа
$\pi _ {s+l+1} (V _ {M+l, M-1}) \cong\pi _ {(s-1) + (l+1) +1} (V _ {M+l+1, M}) $ конечна. Следовательно, и группа $\pi _ {s+l+1} (V _ {M+l, M}) $ конечна.
\end{proof}

\begin{remark}\label{r6.3} Аналогично, используя расслоение $SO_{m-q-p-1}\to SO_{m-q}\to \pi_{q}(V_{m-q,p+1})$, можно доказать следующее утверждение: {\it предположим, что $m\ge p +\frac32 q +2$; тогда группа $\pi_q (V _ {m-q, p+1}) $
бесконечна, если и только если либо $m=2q+1$, $q$ нечетно, либо $m=p+2q+1$, $q$ четно} \cite{CFS11}.
\end{remark}

Доказательство теоремы~\ref{tI} завершено.

\section{Заключительные замечания}\label{s7}

Приведем пример, показывающий, что теорема~\ref{tII} и лемма~\ref{l3} о дополнении не могут быть доказаны с помощью стандартного трюка Уитни или подхода статьи \cite{HaQu74}.

\begin{example} Найдутся такие $p$, $q$ и $m$, удовлетворяющие неравенствам из условия теоремы~\ref{tII}, и такое собственное почти вложение $f:S^p\times D^q\to D^m$, что $\beta (f) =0$, но у $f$ нет паутин. В частности, изоморфизм из леммы~\ref{l3} о дополнении не выполнен при $f' =f$.
\end{example}

\begin{proof} Положим $m=p +\frac {3} {2} q +\frac {3} {2} $, и выберем $p,q$ так, что число $l=m-p-q-1$ нечетно. Возьмем собственное почти вложение $f:S^p\times
D^q\to D^m$ общего положения, такое что многообразие $\Delta$ связно и многообразие $\tilde\Delta$ несвязно (например, начнем со стандартного вложения и сделаем пальцевое движение Уитни). Наше доказательство леммы~\ref{l3} о дополнении в \S4
фактически показывает, что  $\pi_q (D^m-\operatorname {Im} F, \partial) \cong
\Omega_s (\Delta; l\lambda_\Delta) $. Последняя группа изоморфна $\Bbb Z$ при $s=0$ и $l$ нечетном \cite[конец \S4]{HaKa98}. С другой стороны,
$\Omega_0 (P ^\infty; l\lambda) = \Bbb Z_2$. Значит, отображение $\beta (f, {\cdot}):\Bbb Z\to\Bbb Z_2$ не инъективно.

Значит, существует собственное отображение $g:D^q\to D^m-\operatorname {Im} f$,
такое что $\beta (f, g) =0$, но $g$ не (собственно) гомотопно постоянному отображению. Выполним построение шагов (2) и (3) доказательства теоремы~\ref{t5.1} из \S5. Мы получим новое собственное почти вложение $f$, такое
что по-прежнему $\beta (f) =0$. С другой стороны, отображение $g:D^q\to D^m-\operatorname {Im} f$ близко к $f\left | _ {*\times D^q} \right.$,
но не гомотопно постоянному отображению. Значит,
у ~ $f$ нет паутины.
\end{proof}

\subsection{Открытые вопросы} Есть несколько направлений для дальнейшего изучения заузленных торов:

\noindent  (i) {\it Явные классификационные результаты.} Сколько вложений $S^1\times S^5\to S ^ {10} $ существует с точностью до изотопии?

\noindent  (ii) {\it Ослабление размерностных ограничений.} Можно ли отбросить ограничения $m> 2p+q+2$ или $m> p +\frac {4} {3} q+2$ в теореме~\ref{tI}? (сравни с \cite{Hae66A})

\noindent  (iii) {\it Произвольные многообразия.} Было бы интересно обобщить $\beta$-инвариант и теорему~\ref{tII} на вложения произвольных многообразий \cite{ReSk99, Sko08}.

\noindent  (iv) {\it Рациональная классификация вложений.} Для данного многообразия $N$ и числа $m$ определить, является ли множество вложений $N\to S^m$
с точностью до изотопии конечным.

\appendix

\section{Приложение.
Хирургия многообразия двойных точек}

В данном разделе производится хирургия многообразия двойных точек $\Delta$, делающая классифицирующее отображение  $\Delta\to P ^ {\infty} $
накрытия $\tilde\Delta\to\Delta$ достаточно высоко связным. Это необходимо для доказательства леммы~\ref{l3} о дополнении, сформулированной в \S4.
Наше изложение полностью аналогично \cite[приложение ~A]{HaKa98}, но производится в большей общности и более подробно.

Пусть $f:D^n\to M^m$ собственное погружение общего положения, такое что $f\left | _ {\partial D^n} \right.:\partial D^n\to \partial M^m$ --- вложение.
Теорема вложения \cite{Hud69}, сравни с теоремой~\ref{t4.2} выше, позволяет устранить самопересечения погружения $f$ с помощью гомотопии $\operatorname {rel}
\partial D^n$ при некоторых условиях. В диапазоне размерностей, где теорема вложения не верна, мы предлагаем подход к ``упрощению'' двойных точек погружения $f$.

Цель раздела состоит в том, чтобы показать, что классифицирующее отображение
$\Delta\to P ^\infty$ можно сделать $ (s+1)$-связным с помощью гомотопии $\operatorname {rel} \partial D^n$
погружения $f:D^n\to M^m$ при следующих условиях:

\noindent  (i) $M^m$ является $ (s+1)$-связным;

\noindent  (ii) $2s\le 2n-m-2$; и

\noindent  (iii) $0\le s\le m-n-3$.

В этом состоит утверждение теоремы~\ref{t4.6} о хирургии.

\begin{proof}[теоремы~\ref{t4.6} о хирургии]
Превращение отображения $\Delta\to P ^\infty$
в $ (s+1)$-связное происходит в 2 этапа:

\noindent  Шаг 1. Превращение многообразия $\Delta$ в связное и отображения $\pi_1 (\Delta) \to\pi_1 (P ^\infty) $ --- в сюръективное (иными словами, многообразия $\tilde \Delta$ --- в связное).

\noindent  Шаг 2. Убивание элементов ядра
$\operatorname {ker} (\pi_i (\Delta) \to\pi_i (P ^\infty)) $ для каждого $i=1,\dots, s$.

Эти два шага достаточны для доказательства теоремы, потому что отображение
$\pi _ {i+1} (\Delta) \to\pi _ {i+1} (P ^\infty) =0$ сюръективно для каждого
$1\le i\le s$.

В каждом из шагов 1 и 2 мы производим следующий трюк Уитни--Хефлигера, производящий хирургию многообразия $\Delta$.
Сначала построим {\it стандартную модель Хабеггера--Кайзера}
для хирургии оснащенной $i$-мерной сферы, лежащей в двойных точек $n$-мерного диска, погруженного в $S^m$.

\subsection{Стандартная модель для хирургии} \cite{HaKa98} Будем использовать ``модельное'' многообразие $\Bbb R^m =\Bbb R\times \Bbb R ^ {i+1} \times \Bbb R ^ {2n-m-i} \times \Bbb R ^ {m-n-1} \times \Bbb R^{m-n-1}$
и два вложения, $g_+$ и $g_-$, пространства $\Bbb
R^n =\Bbb R ^ {i+1} \times \Bbb R ^ {2n-m-i} \times \Bbb R ^ {n-m-1} $ в $\Bbb R^m$,
пересекающиеся трансверсально по $0\times S^i\times
\Bbb R ^ {2n-m-i} \times 0\times 0 $. Например, можно взять
$g_-(x, y, z) = (|x | ^ 2-1, x, y, 0, z) $ и $ g_ + (x, y, z) = (1-| x | ^ 2, x, y, z, 0) $.
Сфера $S^i$ ограничивает шар $D ^ {i+1} \subset\Bbb
R ^ {i+1} \subset\Bbb R^n$. Далее, сфера
$S ^ {i+1} =D ^ {i+1} _ + \cup D ^ {i+1} _-$, где $D ^ {i+1} _\pm=g_\pm (D ^ {i+1}\times 0\times 0) $, ограничивает шар $D ^ {i+2} \subset\Bbb R\times\Bbb
R ^ {i+1} \subset\Bbb R^m$ с углами вдоль $S^i$. Движение одной из двух шапочек $D ^ {i+2}_\pm $ через $D ^ {i+2} $ производит нужную нам хирургию. Более
точно, двойные точки получающейся регулярной гомотопии
образуют след нужной хирургии.

Теперь сделаем некоторые приготовления к выполнению хирургии, которые немного отличаются в Шагах 1 и 2.

\subsection{Шаг 1: превращение многообразий $\Delta$ и $\tilde \Delta$ в связные}
Если $\Delta =\emptyset$, то мы сначала создаем непустое самопересечение (например, с помощью пальцевого движения Уитни). Пусть $\Delta\ne\emptyset$. Возьмем пару точек $ (a, b), (c, d) $, принадлежащих разным компонентам многообразия $\tilde
\Delta$. Можно считать, что они расположены вне множества тройных точек. Рассмотрим сферы $S^0 = \{\{a, b \}, \{c, d \} \} $, $S^0 _ + = \{
(a, b), (c, d) \} $, и $S^0_-= \{(b, a), (d, c) \} $. Пусть $\eta$ --- тривиализация нормального расслоения $N (S^0, \Delta) $. Хирургия по сфере $S^0$
(см.~Завершение доказательства ниже) соединит различные
компоненты многообразия $\tilde \Delta$, потому что
$\operatorname{dim}\Delta=2n-m\ge 2s+2\ge 2 $.

\subsection{Шаг 2: убивание элементов
$\operatorname {ker}(\pi_i (\Delta) \to\pi_i (P ^\infty)) $} Пусть отображение $g:S^i\to \Delta$ представляет элемент ядра $\operatorname {ker} (\pi_i (\Delta) \to\pi_i (P ^\infty) $, где $1\le i\le s$. Так как $2i\le\operatorname {dim} \Delta-1$ (потому что $2s\le
2n-m-2$), то $g$ можно считать вложением.
По общему положению множество тройных точек имеет размерность $\le 3n-2m$. Так как $s\le m-n-3$, то
$i+3n-2m\le\operatorname {dim} \Delta-1$, так что в общем положении $\operatorname {Im} g$ не содержит тройных точек.

Так как композиция $S^i\to \Delta\to P ^\infty$ тривиальна, то расслоение $ (m-n) \lambda_\Delta\left | _ {S^i} \right.$ тривиально.
Значит, расслоение $\eta=N (S^i, \Delta) $ стабильно тривиально, и стало быть, тривиально. Так как сфероид $S^i\to P ^ {\infty} $ тривиален, то сфера $S^i$ тривиально накрывается в $\tilde \Delta$.
Обозначим $S^i _ + $ и $S^i_-$ две копии $S^i$ в $\tilde
\Delta$.

\subsection{Завершение доказательства: хирургия по сфере $S^i$}
Сначала заклеим сферы $S^i _ + $ и $S^i_-$ непересекающимися
дисками $D ^ {i+1} _ +$ и $D ^ {i+1} _-$ в шаре $D^n$.
Для этого возьмем некоторую тривиализацию расслоения
$N (\tilde \Delta,D^n) $ (например, построенную в определении $(m-n) \lambda$-многообразия двойных точек в \S3). Сдвинем сферы $S^i _ + $
и $S^i_-$ вдоль первого векторного поля этой тривиализации. Так как $s\le m-n-3$ и $2s\le 2n-m-1$, то
$2s\le 2n-m-1$. Следовательно, $i+1 +\operatorname {dim} \Delta \le n-1$ и 2 $ (i+1) \le n-1$. Значит, по соображениям общего положения сдвинутые сферы можно заклеить вложенными дисками в шаре $fD^n$, не пересекающими образ $\tilde i\tilde\Delta$.

Рассмотрим очевидное разложение
$N (D ^ {i+1} _ +, D^n) \left | _ {S^i _ +}\right. =\eta\oplus\epsilon ^ {m-n-1} $,
где $\eta=N (S^i, \Delta) $.
Это разложение задает $(m-n-1)$-оснащение сферы $S^i _ +$. Мы собираемся продолжить это $(m-n-1)$-оснащение на диск $D ^ {i+1} _ +$. Полное препятствие к продолжению лежит в группе
$\pi_i (V _ {n-i-1, m-n-1})$. Последняя группа тривиальна при условии $2i\le 2n-m-1$, которое следует из $2s\le 2n-m-2$. Таким образом, мы получаем разложение
$N (D ^ {i+1} _ +, D^n) = \eta _ +\oplus\epsilon ^ {m-n-1} $, где $\eta _ + $ ---
продолжение расслоения $\eta$ на диск $D ^ {i+1} _ +$,
дополнительное к нашему $(m-n-1)$-оснащению. Определим расслоение $\eta_-$ аналогично.

Теперь продолжим вложение сферы $S ^ {i+1} =D ^ {i+1} _ + \cup D ^ {i+1} _-$ в $S^m$
до вложения диска $D ^ {i+2} $ в $S^m$. Для этого сдвинем полусферы $D ^ {i+1} _ +$ и $D ^ {i+1} _-$ вдоль первого векторного поля некоторой тривиализации расслоения $N (D^n, S^m) $. В итоге получим вложение воротниковой окрестности
сферы $S ^ {i+1} $ в $S^m$.
Так как $n+s+3\le m$ и $s\le n-1$, то $i+2+n\le m-1$ и
$ 2(i+2) \le m-1$. Поэтому по общему положению воротниковая окрестность продолжается до вложенного диска $D ^ {i+2} $ в $S^m$, внутренность которого
не пересекает $fD^n$.

Наконец, рассмотрим следующее $(m-n-1)$-оснащение сферы $S ^ {i+1} $.
На диске $D^{i+1}_+$ возьмем $(m-n-1)$-оснащение, дополнительное к $\eta_+$. На диске
$D ^ {i+1} _-$ возьмем $(m-n-1)$-оснащение, полученное
из тривиализации расслоения $N (D^n, S^m) $ отбрасыванием первого векторного поля.
По построению эти два оснащения совпадают на $S ^ {i} $. Тем самым $(m-n-1)$-оснащение сферы $S ^ {i+1} $ построено. Продолжим его на $D ^ {i+2} $.
Полное препятствие к этому лежит в группе $\pi _ {i+1} (V _ {m-i-2, m-n-1}) $.
Эта группа тривиальна при условии $2i+2\le n-1$, которое следует из $3s\le n-4$.
Пусть $\eta '$ --- дополнительное расслоение  к
полученному $ (m-n-1) $-оснащению диска $D ^ {i+2} $.
Мы имеем разложение $\eta' = \eta_-\oplus\epsilon ^ {m-n-1} $ на диске $D ^ {i+1} _-$.
Продолжая его на диск $D ^ {i+2}$, получаем разложение
$N (D ^ {i+2}, S^m) = \eta '' \oplus\epsilon ^ {m-n-1} \oplus\epsilon ^ {m-n-1} $
для некоторого расслоения $\eta ''$.

Таким образом, наше погружение в окрестности диска $D ^ {i+2} $ устроено так же, как в стандартной модели
(то есть, переводится в стандартную модель
подходящим диффеоморфизмом). Мы можем произвести описанную хирургию и убить сфероид $g:S^i\to \Delta$. Теорема доказана.
\end{proof}

\subsection{Благодарности}
Авторы благодарны А.\,Б.~Скопенкову за постоянное внимание к данной работе, а также П.\,М.~Ахметьеву, У.~Кайзеру, У.~Кошорке, Г.~Лауресу, С.\,А.~Мелихову, А.\,С.~Мищенко, В.\,М.~Нежинскому и Е.\,В.~Щепину за полезные обсуждения.

\end{fulltext}